\documentclass[11pt]{article}
\usepackage{amssymb}
\usepackage{multirow}
\usepackage{graphicx}
\usepackage{amsmath,amssymb,verbatim}
\usepackage{amsthm}
\usepackage{epsfig}
\usepackage[latin1]{inputenc}

\numberwithin{equation}{section}

\usepackage[main=english]{babel}
\newtheorem{definition}{Definition}[section]
\newtheorem{example}{Example}[section]
\newtheorem{theorem}{Theorem}[section]

\newtheorem{lemma}{Lemma}[section]
\newtheorem{corollary}{Corollary}[section]
\newtheorem{proposition}{Proposition}[section]

\begin{document}

%%%%%%%%%%%%%%%%%%%%%%%%%%%%%%%%%%%%%%%%%%%%%%%%%%%%%%%%%%%%%%
% Authors' definitions
\newcommand{\mmmm}{\,\mbox{(mod $m_0$)}}
\newcommand{\nnnn}{\,\mbox{(mod $n_0$)}}
\newcommand{\DDDD}{D_{\rho,\delta}}
\newcommand{\DDDA}{\DDD(A_0^{\gamma})}
\newcommand{\UPP}{\overline{u}}
\newcommand{\LOP}{\underline{u}}
\newcommand{\weight}{e^{2s\varphi}}
\newcommand{\tilweight}{e^{2s\widetilde{\va}}}
\newcommand{\ep}{\varepsilon}
\newcommand{\la}{\lambda}
\newcommand{\va}{\varphi}
\newcommand{\ppp}{\partial}
\newcommand{\izt}{\int^t_0}
\newcommand{\ddd}{\mathcal{D}}
\newcommand{\www}{\widetilde}
\newcommand{\whwh}{\widehat}
\newcommand{\pppg}{\partial_{t,\gamma}^{\alpha}}
\newcommand{\pppa}{\partial_t^{\alpha}}
\newcommand{\lat}{-\lambda_n t^{\alpha}}
\newcommand{\ddda}{d_t^{\alpha}}
\newcommand{\DDDa}{D_t^{\alpha}}
\newcommand{\HOLST}{C^{2+\theta}(\ooo{\OOO})}
\newcommand{\HOLSZ}{C^{\theta}(\ooo{\OOO})}
\newcommand{\CONE}{C^1(\ooo{\OOO})}
\newcommand{\lla}{L^{\frac{3}{4}}}
\newcommand{\llb}{L^{-\frac{1}{4}}}
\newcommand{\UPPY}{\ooo{y}}
\newcommand{\LOWY}{\underline{y}}
\newcommand{\lowa}{\underline{a}}
\newcommand{\Ahalf}{A_0^{\frac{1}{2}}}
\newcommand{\QQQQ}{\Omega \times (0,T)}
\newcommand{\HHMIN}{H^{-1}(\Omega)}
\newcommand{\HHONE}{H^1_0(\Omega)}
\newcommand{\uab}{u_{a,b}}
\newcommand{\Halp}{H_{\alpha}(0,T)}
\newcommand{\Hal}{H^{\alpha}(0,T)}
\newcommand{\lllb}{L^{-\frac{3}{4}}}
\newcommand{\pdif}[2]{\frac{\partial #1}{\partial #2}}
\newcommand{\ppdif}[2]{\frac{\partial^2 #1}{{\partial #2}^2}}
\newcommand{\uen}{u_N^{\ep}}
\newcommand{\fen}{F_N^{\ep}}
\newcommand{\penk}{p^{\ep}_{N,k}}
\newcommand{\penl}{p^{\ep}_{N,\ell}}
\newcommand{\R}{\mathbb{R}}
\newcommand{\Q}{\mathbb{Q}}
\newcommand{\Z}{\mathbb{Z}}
\newcommand{\C}{\mathbb{C}}
\newcommand{\N}{\mathbb{N}}
\newcommand{\uu}{\mathbf{u}}
\renewcommand{\v}{\mathbf{v}}
\newcommand{\y}{\mathbf{y}}
\newcommand{\RR}{\mathbf{R}}
\newcommand{\Y}{\mathbf{Y}}
\newcommand{\w}{\mathbf{w}}
\newcommand{\z}{\mathbf{z}}
\newcommand{\G}{\mathbf{G}}
\newcommand{\ooo}{\overline}
\newcommand{\OOO}{\Omega}
\newcommand{\MLONE}{E_{\alpha,1}}
\newcommand{\MLTWO}{E_{\alpha,\alpha}}
\newcommand{\MLAA}{E_{\alpha,\alpha}}
\newcommand{\oo}{\omega}
\newcommand{\sumjn}{\sum_{j=1}^{\ell}}
\newcommand{\sumjm}{\sum_{j=1}^m}
\newcommand{\sumij}{\sum_{i,j=1}^d}
\newcommand{\suml}{\sum_{\ell=1}^p}
\newcommand{\sumN}{\sum_{k=1}^N}
\newcommand{\nene}{\not\equiv}
\newcommand{\NUNU}{\ppp_{\nu_A}}
\newcommand{\CCO}{{_{0}C^1[0,T]}}
\newcommand{\CC}{{_{0}C^1[0,T]}}
\newcommand{\dddaO}{{_{0}\ddda}}
\newcommand{\WWW}{{_{0}W^{1,1}(0,T)}}
\newcommand{\CCOO}{{^{0}C^1[0,T]}}
\newcommand{\HH}{H_{\alpha}}
\newcommand{\HHH}{_{0}{H^{\alpha}}(0,T)}
\newcommand{\HHO}{_{0}H^1(0,T)}
\newcommand{\hhalf}{\frac{1}{2}}
\newcommand{\DDD}{\mathcal{D}}
\newcommand{\RRR}{\mathcal{R}}
\newcommand{\RRRR}{\longrightarrow}
\newcommand{\lbra}{_{H^{-1}(\Omega)}{\langle}}
\newcommand{\wdel}{\www{\delta}}
\newcommand{\weps}{\www{\varepsilon}}
\newcommand{\sumn}{\sum_{n=1}^{\infty}}
\newcommand{\pfpf}{(P_nf)(x_0)}
\newcommand{\pgpg}{(P_ng)(x_0)}
\newcommand{\pkpk}{p\in \{p_k\}_{k\in \N}}
\newcommand{\rbra}{\rangle_{H^1_0(\Omega)}}
\newcommand{\lbran}{_{(H^{-1}(\Omega))^N}{\langle}}
\newcommand{\rbran}{\rangle_{(H^1_0(\Omega))^N}}
\newcommand{\llll}{L^{\infty}(\Omega\times (0,t_1))}
\renewcommand{\div}{\mathrm{div}\,}  %div
\newcommand{\grad}{\mathrm{grad}\,}  %grad
\newcommand{\rot}{\mathrm{rot}\,}  %rot
\newcommand{\urur}{\left( \frac{1}{r}, \theta\right)}

%%% ........... %%%%%%%%%%%%%%%%%%%%%%%%%%%%%%%%%%

%%%%%%%%%%%%%%%%%%%%%%%%%%%%%%%%%%%%%%%%%%%%%%%%%%

%%%%%%%%%%%%%%% begin make title %%%%%%%%%%%%%%%%%

%\journalname{Fract. Calc. Appl. Anal.} %%%%% and FCAA-logo will appear on the right %%

%\begin{document}

 %%% Pls. specify the kind of the article:
 %%% This is: %% ORIGINAL PAPER %%%% 

%%%% Title of article for FCAA %%%%%%%%%%%%%%%%%%%%%%%%
%\title{Comparison principles for the linear and semiliniar time-fractional
 %       diffusion equations with the Robin boundary condition}

%\titlerunning{Comparison principles for the linear and semiliniar time-fractional equations \dots}
%% if too long for running head - use the text from 1st line !%

\centerline{{\Large{\bf Comparison principles for the linear and semiliniar time-fractional
}}}
\vspace{0.2cm}

\centerline{{\Large{\bf 
        diffusion equations with the Robin boundary condition}}}

\vspace{0.3cm}

%%%%% authors:

\vspace{0.3cm}
\centerline{{\bf Yuri Luchko$^{1,*}$ and  Masahiro Yamamoto$^2$}}
\vspace{0.3cm}

\centerline{$^{1}${Department of Mathematics, Physics, and Chemistry}}

\centerline{{Berlin University of Applied Sciences and Technology}}

\centerline{{Luxemburger Str. 10, 13353 Berlin,\ Germany}}

\centerline{{e-mail: luchko@bht-berlin.de}$^*$ corresponding author }

%\centerline{{telephone number: +49-30-45045295}}

\vspace{0.2cm}

\centerline{$^{2}$

Department of Mathematical Sciences}

\centerline{The University of Tokyo}

\centerline{Komaba, Meguro, 
153 Tokyo, Japan}

\centerline{Honorary Member of Academy of Romanian Scientists}

\centerline{Ilfov, nr. 3, Bucuresti, Romania}

\centerline{Correspondence member of Accademia Peloritana dei Pericolanti}

\centerline{Palazzo Universit\`a, Piazza S. Pugliatti 1,
98122 Messina, Italy}

\centerline{e-mail: myama@ms.u-tokyo.ac.jp}

%%% Dates:
%\date{Received: XX 2022 / Revised: ....  / Accepted: ......}
% These dates will be entered by the Editor (EiC, V. Kiryakova) or authors/ system %

%%%%% For Production Dept.: variants of COPYRIGHT notice to appear: %%%%%%
% if Open Access option chosen, it is as:
   %% "\copyright The Author(s) 2022" % (or next year..) %
% if No Open Access, the TCA form signed to FCAA journal is used, and it appears:
   %% "\copyright Diogenes Co. Ltd. 2022" %% or 2021, or next year %%%%%
%%%%%%%%%%%%%%%%%%%%%%%%%%%%%%%%%%%%%%%%%%%%%%%%%%%%%%%%%%%%%%%%%%%%%%%%%%

%\maketitle

%%%%%%%%%%%%%%%%%%%%%%%%%%%%%%%%%%%%%%%%%%%%%%%%%%%%%%
\begin{abstract}
{The main objective of this paper is analysis of the initial-boundary value problems for the linear and semilinear time-fractional diffusion equations
with a uniformly elliptic spatial differential operator of the second order and
the Caputo type fractional derivative  acting in the fractional Sobolev spaces.
The boundary conditions are formulated in form of  the homogeneous Neumann or 
Robin conditions. First we deal with the uniqueness and existence of the 
solutions to these initial-boundary value problems. Then we show a positivity 
property for the solutions and derive the corresponding comparison principles.
In the case of the semilinear time-fractional diffusion equation,
we also apply 
the monotonicity method by upper and lower solutions.
As an application of our results, we present some a priori estimates 
for solutions to the semilinear time-fractional diffusion equations.}

\end{abstract} 

%%%%%%%%%% Enter suitable key words and phrases %%% examples:

\vspace{0.2cm}


\noindent
{\sl Key Words}: fractional calculus; 
fractional diffusion equation; positivity of solutions; 
comparison principle; monotonicity method; upper and lower solutions
%%% These are examples only. Pls. use MSC 2020 for suitable topic numbers %%%%%%

%%%%%%%%%%%%%

%%%%%%%% begin papers' body %%%%%%%%%%%%%%%%%%%%%%%%%%%%%

%%%%%%%%%%%% Section 1 %%%%%%%%%%%%%%%%%%%%%%%%%%
\section{Introduction} \label{sec:1}

\setcounter{section}{1} 
%\setcounter{equation}{0} %% to have proper 2-digits numbers of eqs
%% Note that this style produces 1-digit numbering of definitons, statements, exmaples, etc.

In this paper, we first deal with a linear time-fractional 
diffusion equation
$$
 \pppa (u(x,t)-a(x)) = \sum_{i,j=1}^d \ppp_i(a_{ij}(x)\ppp_j u(x,t))
$$
\begin{equation}
\label{(1.1)}
+ \sum_{j=1}^d b_j(x,t)\ppp_ju(x,t) + c(x,t)u(x,t) + F(x,t),\quad 
 x \in \Omega,\, 0<t<T,                                         
\end{equation}
where $\pppa$ is the Caputo fractional derivative of order 
$\alpha\in (0,1)$ defined in the fractional Sobolev spaces 
(see Section \ref{sec2}) and $\OOO \subset \R^d, \ d=1,2,3$ is
a bounded domain with smooth boundary $\ppp\OOO$.   
All functions under consideration are assumed to be real-valued.

In what follows, we always assume that the following conditions are satisfied:
\begin{equation}
\label{(1.2)}
\left\{ \begin{array}{rl}
& a_{ij} = a_{ji} \in C^1(\ooo{\OOO}), \quad 1\le i,j \le d, \\
& b_j,\, c \in C^1([0,T]; C^1(\ooo{\OOO})) \cap C([0,T];C^2(\ooo{\OOO})),
  \quad 1\le j \le d, \\
& \mbox{and there exists a constant $\kappa>0$ such that}\\
& \sumij a_{ij}(x)\xi_i\xi_j \ge \kappa \sum_{j=1}^d \xi_j^2, \quad
x\in \OOO, \, \xi_1, ..., \xi_d \in \R.
\end{array}\right.
\end{equation}

Using the notations $\ppp_j = \frac{\ppp}{\ppp x_j}$, $j=1, 2, ..., d$ and 
$\ppp_t = \frac{\ppp}{\ppp t}$, we define the conormal derivative
$\NUNU w$ with respect to $\sumij \ppp_j(a_{ij}\ppp_i)$ by 
\begin{equation}\label{(1.3)}
\NUNU w(x) =  \sumij a_{ij}(x)\ppp_jw(x)\nu_i(x), \quad x\in \ppp\OOO,
\end{equation}
where $\nu = \nu(x) =: (\nu_1(x), ..., \nu_d(x))$ is the unit outward normal
vector to $\ppp\OOO$ at the point $x := (x_1,..., x_d) \in \ppp\OOO$. 

For the equation \eqref{(1.1)},  we consider the initial-boundary value 
problems with the Neumann boundary condition
\begin{equation}
\label{(1.3a)}
\ppp_{\nu_A}u = 0 \quad \mbox{on $\ppp\OOO \times (0,T)$}   
\end{equation}
or with the Robin boundary condition
\begin{equation}
\label{(1.4)}
\ppp_{\nu_A}u + \sigma(x)u = 0 \quad \mbox{on $\ppp\OOO \times (0,T)$},   
\end{equation}
where  $\sigma$ is a smooth function on $\ppp\OOO$ that satisfies the condition $\sigma \ge 0$.

In the second part of the article, we treat a semilinear time-fractional 
diffusion equation
$$
\pppa (u(x,t)-a(x)) = \sum_{i,j=1}^d \ppp_i(a_{ij}(x)\ppp_j u(x,t))
$$
\begin{equation}
\label{(1.5)}
+ \sum_{j=1}^d b_j(x,t)\ppp_ju(x,t) + c(x,t)u(x,t) 
+ F(x,u(x,t), \nabla u(x,t)),\  x \in \Omega,\, 0<t<T \\
\end{equation}
with the boundary condition \eqref{(1.3a)} or \eqref{(1.4)}.

For the parabolic partial differential equations which 
correspond to the case $\alpha=1$ in the equations \eqref{(1.1)} and 
\eqref{(1.5)},  several 
important qualitative properties of solutions to the initial-boundary 
value problems are known.  In particular, 
we mention the maximum principle and the comparison principle
(\cite{PW}, \cite{RR}).

The main purpose of this article is in establishing the comparison principles 
and the monotonicity method for the linear and semilinear time-fractional  
diffusion equations \eqref{(1.1)} and \eqref{(1.5)} with the Neumann or
Robin boundary conditions.

For the linear time-fractional diffusion equations of type \eqref{(1.1)} 
with the Dirichlet boundary condition,  the maximum principle was  studied and 
used in  \cite{Lu1}, \cite{Lu2},
 \cite{LY1},  \cite{LY2}, \cite{LY3}, \cite{Za}. For derivation of the
maximum principle for the time-fractional transport equations
we refer to \cite{LSY}.

Because the maximum principle involves the Dirichlet boundary values, it is
difficult to formulate it in the case of the Neumann or Robin boundary 
conditions. However, for this kind of the boundary conditions, the comparison 
principle and the positivity of solutions can be proved that is the main 
objective of this article.  One typical result of this sort says
that the solution $u$ to the equation \eqref{(1.1)} with the boundary 
condition \eqref{(1.4)} and an appropriately formulated initial condition is 
non-negative in $\OOO\times (0,T)$ 
if the initial value $a$ and the non-homogeneous term $F$ are non-negative 
in $\OOO$ and in $\QQQQ$, respectively.
Such positivity properties and their applications were intensively derived 
and used for the parabolic partial differential equations ($\alpha=1$ 
in the equations \eqref{(1.1)} and \eqref{(1.5)}), see e.g., \cite{E}, 
 \cite{Fr}, \cite{Pao2} or \cite{RR}.
Moreover, we refer to  \cite{Am},  \cite{Ke},  \cite{Pao2}
for a discussion of the monotonicity method by upper and lower solutions for 
studying the semilinear parabolic partial differential equations. 
This method is based on the corresponding comparison principle.  

However, to the best knowledge of the authors, there are no such works for 
the time-fractional diffusion equations in the case of the Neumann or Robin
boundary conditions. The main purpose of this article is to fill this gap in 
the theory of the linear and semilinear time-fractional  diffusion equations.  
The positivity property and the comparison principle for the linear equation 
\eqref{(1.1)} with the boundary condition
\eqref{(1.3a)} or \eqref{(1.4)} and an appropriately formulated initial 
condition are the subject of the first part of the article.
In the second part, these result are extended to the semilinear 
time-fractional diffusion 
equation \eqref{(1.5)}.  The arguments employed in this article 
rely on an operator theoretical approach to the fractional integrals and 
derivatives in the fractional Sobolev spaces that is an extension of the 
theory developed for the case $\alpha=1$ (\cite{He}, \cite{Pa}, \cite{Ta}).

The rest of this article is organized as follows. In Section \ref{sec2},  
some important results regarding the unique existence of solutions 
to the initial-boundary value problems for the linear time-fractional  
diffusion equations are presented. Section \ref{sec3} is devoted to a proof
of a key lemma that is a basis for the proofs of the comparison principles 
discussed in Sections \ref{sec4} and \ref{sec6}. The lemma asserts that the
solution mapping $\{a, F\} \longrightarrow
u_{a,F}$ preserves the ordering of $a$ and $F$, 
where $a$ and $F$ denote an initial condition 
and a source function of the 
problem under consideration, respectively, and $u_{a,F}$ denotes its  
solution in a certain class.  In Section \ref{sec4}, 
the key lemma and the 
fixed point theorem are employed to prove the comparison principles and some 
of their corollaries for the initial-boundary value problems for the linear 
time-fractional  diffusion equations.  
In Section \ref{sec5}, the unique existence of solutions to the 
initial-boundary value problems for the semilinear time-fractional  
diffusion equations are discussed. Section \ref{sec6} presents the comparison 
principles and some of their corollaries for the initial-boundary value 
problems for the semilinear time-fractional  diffusion equations.
Finally, some conclusions and remarks are formulated in the last section.

%%%%%%%%%%% Section 2 %%%%%%%%%%%%%%%
\section{Well-posedness for initial-boundary value problems for 
linear time-fractional diffusion equations} 
\label{sec2}

\setcounter{section}{2}

For $x \in \Omega, \thinspace 0<t<T$,  we define the operator
\begin{equation}
\label{(2.1)}
-Au(x,t) = \sum_{i,j=1}^d \ppp_i(a_{ij}(x)\ppp_jv(x,t)
+ \sum_{j=1}^d b_j(x,t)\ppp_ju(x,t) + c(x,t)u(x,t)            
\end{equation}
and assume that the conditions \eqref{(1.2)} for the coefficients 
$a_{ij}(x), b_j(x,t), c(x,t)$ are satisfied. 
 
In this section, we deal with the following initial-boundary value problem for 
the linear  time-fractional diffusion equation with the fractional derivative 
of order $\alpha\in (0,1)$
\begin{equation}
\label{(2.2)}
\left\{ \begin{array}{rl}
& \pppa (u(x,t)-a(x)) + Au(x,t) = F(x,t), 
\quad x \in \Omega, \thinspace 0<t<T, \\
& \NUNU u + \sigma(x)u(x,t) = 0, \quad x\in \ppp\OOO, \, 0<t<T,
\end{array}\right.
\end{equation}
along with the initial condition \eqref{incon} formulated below. 
To appropriately define the  Caputo derivative $\ddda w(t)$, $0<\alpha<1$, 
we start with its definition on the space
$$
\CC := \{ u \in C^1[0,T];\thinspace u(0) = 0\}
$$
that reads as follows:  
$$
\ddda w(t) = \frac{1}{\Gamma(1-\alpha)}\int^t_0
(t-s)^{-\alpha}\frac{dw}{ds}(s) ds,\ w\in 
\CC.
$$
Then we extend this operator from the domain $\mathcal{D}(\ddda)
:= \CC$ into $L^2(0,T)$ taking into account the closability of the operator 
(\cite{Yo}).  As have been shown in \cite{KRY}, there exists a unique minimum 
closed 
extension of $\ddda$ with $\mathcal{D}(\ddda) = \CC$.
Moreover, the domain of this 
extension is the closure of $\CC$ in the Sobolev-Slobodecki space 
$H^{\alpha}(0,T)$.  Let us recall that the norm $\Vert \cdot\Vert_{H^{\alpha}
(0,T)}$ is defined as follows (\cite{Ad}):
$$
\Vert v\Vert_{H^{\alpha}(0,T)}:=
\left( \Vert v\Vert^2_{L^2(0,T)}
+ \int^T_0\int^T_0 \frac{\vert v(t)-v(s)\vert^2}{\vert t-s\vert^{1+2\alpha}}
dtds \right)^{\hhalf}.
$$
By setting 
$$
H_{\alpha}(0,T):= \ooo{\CC}^{H^{\alpha}(0,T)},
$$
we see that 
$$
H_{\alpha}(0,T) = 
\left\{ \begin{array}{rl}
&H^{\alpha}(0,T), \quad  0<\alpha<\hhalf, \\
&\left\{ v \in H^{\hhalf}(0,T);\, \int^T_0 \frac{\vert v(t)\vert^2}{t}
dt < \infty \right\}, \quad  \alpha=\hhalf, \\
& \{ v \in H^{\alpha}(0,T);\, v(0) = 0\}, \quad \hhalf < \alpha < 1,
\end{array}\right.
$$
and 
$$
\Vert v\Vert_{H_{\alpha}(0,T)} = 
\left\{ \begin{array}{rl}
&\Vert v\Vert_{H^{\alpha}(0,T)}, \quad  \alpha \ne \hhalf, \\
&\left( \Vert v\Vert_{H^{\hhalf}(0,T)}^2
+ \int^T_0 \frac{\vert v(t)\vert^2}{t}dt\right)^{\hhalf}, \quad
\alpha=\hhalf.
\end{array}\right.
$$
In what follows, we also use the Riemann-Liouville fractional
integral operator $J^{\beta}$, $\beta > 0$ defined by 
$$
J^{\beta}f(t) := \frac{1}{\Gamma(\beta)}\int^t_0 (t-s)^{\beta-1}f(s) ds,
\quad 0<t<T.
$$
Then, according to \cite{GLY} and \cite{KRY}, 
$$
H_{\alpha}(0,T) = J^{\alpha}L^2(0,T),\quad 0<\alpha<1.
$$
We define 
$$
\pppa = (J^{\alpha})^{-1} \quad \mbox{with $\mathcal{D}(\pppa)
= H_{\alpha}(0,T)$}.
$$
Then there exists a constant $C>0$ depending only on $\alpha$ such that 
$$
C^{-1}\Vert v\Vert_{H_{\alpha}(0,T)} \le \Vert \pppa v\Vert_{L^2(0,T)}
\le C\Vert v\Vert_{H_{\alpha}(0,T)} \quad \mbox{for all } v\in H_{\alpha}(0,T).
$$

Now we can introduce a suitable form for an initial condition for the problem 
\eqref{(2.2)} as follows
\begin{equation}
\label{incon}
u(x, \cdot) - a(x) \in \HH(0,T) \quad \mbox{for almost all } x\in \OOO
\end{equation}
and write down the complete formulation of the initial-boundary value problem 
for the linear  time-fractional diffusion equation:  
\begin{equation}
\label{(2.3)}
\left\{ \begin{array}{rl}
& \pppa (u(x,t)-a(x)) + Au(x,t) = F(x,t), 
\quad x \in \Omega, \thinspace 0<t<T, \\
& \NUNU u(x,t) + \sigma(x)u(x,t) = 0, \quad x\in \ppp\OOO, \, 0<t<T,\\
& u(x, \cdot) - a(x) \in \HH(0,T) \quad \mbox{for almost all }x\in \OOO.
\end{array}\right.
\end{equation}
We note that the term $\pppa (u(x,t) - a(x))$  in \eqref{(2.3)} is 
well-defined via the third condition in \eqref{(2.3)}.
In particular, for $\frac{1}{2} < \alpha < 1$, the Sobolev embedding leads 
to the inclusions
$\HH(0,T) \subset H^{\alpha}(0,T)\subset C[0,T]$. 
Thus, $u\in \HH(0,T;L^2(\OOO))$ implies  $u \in C([0,T];L^2(\OOO))$  
and in this case, we can see that the initial condition is formulated as 
$u(\cdot,0) = a$  in $L^2$-sense.

In the following theorem, a fundamental result regarding the unique existence 
of the solution to 
the initial-boundary value problem \eqref{(2.3)} is formulated.  

\begin{theorem}
\label{t2.1}
For $a\in H^1(\OOO)$ and $F \in L^2(0,T;L^2(\OOO))$, there exists a unique 
solution 
$u(F,a) = u(F,a)(x,t) \in L^2(0,T;H^2(\OOO))$ to the initial-boundary value 
problem \eqref{(2.3)} such that 
$u(F,a)-a \in \HH(0,T;L^2(\OOO))$.
Moreover, there exists a constant $C>0$ such that 
\begin{align*}
& \Vert u(F,a)-a\Vert_{\HH(0,T;L^2(\OOO))} 
+ \Vert u(F,a)\Vert_{L^2(0,T;H^2(\OOO))} \\
\le &C(\Vert a\Vert_{H^1(\OOO)} + \Vert F\Vert_{L^2(0,T;L^2(\OOO))}).
\end{align*}
\end{theorem}

\begin{proof}
For the proof, we introduce an elliptic operator $A_0$ in $L^2(\OOO)$
with an arbitrary $c_0>0$ as follows:
\begin{equation}
\label{(3.2)}
\left\{ \begin{array}{rl}
& (-A_0v)(x) = \sumij \ppp_i(a_{ij}(x)\ppp_jv(x)) - c_0v(x), \quad
x\in \OOO, \\
& \mathcal{D}(A_0) = \left\{ v\in H^2(\OOO);\,
\NUNU v + \sigma v = 0 \quad \mbox{on } \ppp\OOO \right\}.
\end{array}\right.
\end{equation}
We recall that $\sigma$ is smooth and $\sigma\ge 0$ on $\ppp\OOO$ 
and the coefficients $a_{ij}$, $b_j$ and $c$ satisfy the conditions 
\eqref{(1.2)}. 

Henceforth, by $\Vert \cdot\Vert$ and $(\cdot,\cdot)$ we denote the 
standard norm and
the scalar product  in $L^2(\OOO)$, respectively. 
It is well-known that the operator $A_0$ is self-adjoint and its resolvent 
is a compact operator.  Moreover, for sufficiently large constant
$c_0>0$, by Lemma \ref{lem1} in Section \ref{sec8}, we can verify that $A_0$ is 
positive definite.
Therefore, by choosing a large constant $c_0>0$,
its spectrum consists entirely of discrete positive 
eigenvalues $0 < \la_1 \le \la_2 \le \cdots$,
which are numbered according to their multiplicities and  
$\la_n \to \infty$ as $n\to \infty$.
Let $\va_n$ be an eigenvector corresponding to the eigenvalue $\la_n$ such that
$(\va_n, \va_m) = 0$ if $n \ne m$ and $(\va_n,\va_n) = 1$.
Note that $A\va_n = \la_n\va_n$ for $n\in\N$.
Then the system $\{ \va_n\}_{n\in \N}$ of the 
eigenvectors forms an
orthonormal basis in $L^2(\OOO)$ and for any $\gamma\ge 0$ we can define
the fractional powers $A_0^{\gamma}$ of the operator $A_0$ 
by the following relation (see, e.g., \cite{Pa}):
$$
A_0^{\gamma}v = \sum_{n=1}^{\infty} \la_n^{\gamma} (v,\va_n)\va_n,
$$
where
$$
v \in \mathcal{D}(A_0^{\gamma})
:= \left\{ v\in L^2(\OOO): \thinspace
\sum_{n=1}^{\infty} \la_n^{2\gamma} (v,\va_n)^2 < \infty\right\}
$$
and
$$
\Vert A_0^{\gamma}v\Vert = \left( \sum_{n=1}^{\infty}
\la_n^{2\gamma} (v,\va_n)^2 \right)^{\frac{1}{2}}.
$$
We note that $\mathcal{D}(A_0^{\gamma}) \subset H^{2\gamma}(\OOO)$.

The proof of Theorem \ref{t2.1} is similar to case of the homogeneous Dirichlet 
boundary condition (\cite{GLY}, \cite{KRY}), and we 
define two operators $S(t)$ and $K(t)$ by
\begin{equation}
\label{(5.1)}
S(t)a = \sum_{n=1}^{\infty} E_{\alpha,1}(-\la_n t^{\alpha})
(a,\va_n)\va_n, \quad a\in L^2(\OOO), \thinspace t>0  
\end{equation}
and
\begin{equation}
\label{(5.2)}
K(t)a = -A_0^{-1}S'(t)a 
= \sum_{n=1}^{\infty} t^{\alpha-1}E_{\alpha,\alpha}(-\la_n t^{\alpha})
(a,\va_n)\va_n, \quad a\in L^2(\OOO), \thinspace t>0.
\end{equation}
Here $E_{\alpha,\beta}(z)$ denotes the Mittag-Leffler function defined by 
a convergent series as follows:
$$
E_{\alpha,\beta}(z) = \sum_{k=0}^\infty \frac{z^k}{\Gamma(\alpha\, k + \beta)},
\ \alpha >0,\ \beta \in \mathbb{C},\ z \in \mathbb{C}.
$$
It follows directly from the definitions given above that
$A_0^{\gamma}K(t)a = K(t)A_0^{\gamma}a$
and $A_0^{\gamma}S(t)a = S(t)A_0^{\gamma}a$ for $a \in \mathcal{D}
(A_0^{\gamma})$.
Moreover, since 
$$
\vert E_{\alpha,1}(-\la_nt^{\alpha})\vert, \, 
\vert E_{\alpha,\alpha}(-\la_nt^{\alpha})\vert \le \frac{C}{1+\la_nt^{\alpha}}
\quad \mbox{for all $t>0$}
$$
(e.g., Theorem 1.6 (p.35) in \cite{Po}), we can prove 
\begin{equation}
\label{(5.3)}
\left\{ \begin{array}{l}
\Vert A_0^{\gamma}S(t)a\Vert \le Ct^{-\alpha\gamma}\Vert a\Vert, \\
\Vert A_0^{\gamma}K(t)a\Vert \le Ct^{\alpha(1-\gamma)-1}
\Vert a\Vert, \quad a \in L^2(\OOO), \thinspace t > 0, \thinspace
0 \le \gamma \le 1
\end{array}\right.                    
\end{equation}
(\cite{GLY}).
In order to shorten the notations and focus on the dependence of the time 
variable $t$, henceforth we sometimes omit the variable  $x$ in the functions 
of two variables $x$ and $t$, and write simply
$u(t) = u(\cdot,t)$, $F(t) = F(\cdot,t)$, $a = a(\cdot)$, etc.

Because of \eqref{(5.3)}, the remaining estimation can be carried out analog 
to the one in
the case of the fractional powers of generators of analytic semigroups 
(\cite{He}). 

To do this, we first formulate and prove the following lemma:

\begin{lemma}
\label{l5.1}
Under the conditions formulated above, the following estimates hold true for 
$F\in L^2(0,T;L^2(\OOO))$ and  $a \in L^2(\OOO)$:

\noindent
(i) 
$$
\left\Vert \int^t_0 A_0K(t-s)F(s) ds \right\Vert_{L^2(0,T;L^2(\OOO))}
\le C\Vert F\Vert_{L^2(0,T;L^2(\OOO))},
$$
\noindent
(ii)
$$
\left\Vert \int^t_0 K(t-s)F(s) ds \right\Vert_{\HH(0,T;L^2(\OOO))}
\le C\Vert F\Vert_{L^2(0,T;L^2(\OOO))},
$$
\noindent
(iii) 
$$
\Vert S(t)a - a\Vert_{\HH(0,T;L^2(\OOO))}
+ \Vert S(t)a\Vert_{L^2(0,T;H^2(\OOO))} \le C\Vert a\Vert.
$$
\end{lemma}

\begin{proof} 
We start with the proof of the estimate (i). By \eqref{(5.2)}, we have
\begin{align*}
& \left\Vert \int^t_0 A_0 K(t-s)F(s) ds \right\Vert^2\\
=& \left\Vert \sumn \left(\int^t_0 \la_n(t-s)^{\alpha-1}
\MLTWO(-\la_n(t-s)^{\alpha})
(F(s), \va_n) ds\right) \va_n\right\Vert^2\\
=& \sumn \left\vert \int^t_0 \la_n(t-s)^{\alpha-1}\MLTWO(-\la_n(t-s)^{\alpha})
(F(s), \va_n) ds \right\vert^2.
\end{align*}
Therefore, using the Parseval equality and the Young inequality for the 
convolution, we obtain
\begin{align*} 
& \left\Vert \int^t_0 A_0K(t-s)F(s) ds \right\Vert^2_{L^2(0,T;L^2(\OOO))}\\
= & \sumn \int^T_0 \vert (\la_ns^{\alpha-1}\MLTWO(-\la_ns^{\alpha}) \, *\,
(F(s), \va_n) \vert^2 ds\\
= & \sumn \Vert \la_ns^{\alpha-1}\MLTWO(-\la_ns^{\alpha}) \, * \,
(F(s), \va_n) \Vert^2_{L^2(0,T)}\\
\le & \sumn \left( \la_n\int^t_0 \vert t^{\alpha-1}\MLTWO(-\la_nt^{\alpha})
\vert dt \right)^2 \Vert (F(t),\va_n)\Vert^2_{L^2(0,T)}.
\end{align*} 
Then we employ the representation
\begin{equation}\label{(2.9a)}
\frac{d}{dt}\MLONE(-\la_nt^{\alpha}) 
= -\la_nt^{\alpha-1}\MLTWO(-\la_nt^{\alpha}),
\end{equation}
and the complete monotonicity of the Mittag-Leffler function (\cite{GKMR})
$$
\MLONE(-\la_nt^{\alpha}) > 0, \quad
\frac{d}{dt}\MLONE(-\la_nt^{\alpha}) \le 0, \quad t\ge 0, \quad 0<\alpha\le 1
$$
to reach the inequality
\begin{equation}
\label{(5.4)}
\int^T_0 \vert \la_nt^{\alpha-1}\MLTWO(-\la_nt^{\alpha})\vert dt 
= -\int^T_0 \frac{d}{dt}\MLONE(-\la_nt^{\alpha})dt
\end{equation}
$$
= 1 - \MLONE(-\la_nT^{\alpha}) \le 1 \quad \mbox{for all $n\in \N$}.
$$
Hence,
\begin{align*}
& \left\Vert \int^t_0 A_0K(t-s)F(s) ds \right\Vert^2_{L^2(0,T;L^2(\OOO))}
\le \sumn \Vert (F(t), \va_n)\Vert^2_{L^2(0,T)}\\
=& \int^T_0 \sumn \vert (F(t), \va_n) \vert^2 dt
= \int^T_0 \Vert F(\cdot,t)\Vert^2 dt 
= \Vert F\Vert_{L^2(0,T;L^2(\OOO))}^2.
\end{align*}

Now we proceed with the proof of the estimate (ii). For $0<t<T, \, n\in \N$
and $f\in L^2(0,T)$,
we set 
$$
(L_nf)(t) := \int^t_0 (t-s)^{\alpha-1}\MLTWO(-\la_n(t-s)^{\alpha}) 
f(s) ds.
$$
Then 
$$
\int^t_0 K(t-s)F(s) ds = \sumn (L_nf)(t)\va_n
$$
in $L^2(\OOO)$ for any fixed $t \in [0,T]$.

First we prove that 
\begin{equation}
\label{(5.5)}
\left\{ \begin{array}{rl}
& L_nf \in \HH(0,T), \\
& \pppa(L_nf)(t) = -\la_nL_nf(t) + f(t), \quad 0<t<T, \\
& \Vert L_nf\Vert_{\HH(0,T)} \le C\Vert f\Vert_{L^2(0,T)},
\quad n\in \N \quad \mbox{for each } f\in L^2(0,T).            
\end{array}\right.
\end{equation}
In order to prove this, we apply to $L_nf$ the Riemann-Liouville fractional 
integral operator $J^{\alpha}$:
\begin{align*}
& J^{\alpha}(L_nf)(t) 
= \frac{1}{\Gamma(\alpha)}\int^t_0 (t-s)^{\alpha-1}(L_nf)(s) ds\\
=& \frac{1}{\Gamma(\alpha)} \int^t_0 (t-s)^{\alpha-1}
\left( \int^s_0  (s-\xi)^{\alpha-1}\MLTWO(-\la_n (s-\xi)^{\alpha})f(\xi) d\xi
\right) ds\\
=& \frac{1}{\Gamma(\alpha)}\int^t_0 f(\xi) \left( 
\int^t_{\xi} (t-s)^{\alpha-1}(s-\xi)^{\alpha-1}\MLTWO(-\la_n (s-\xi)^{\alpha})
ds \right) d\xi.
\end{align*}
By direct calculation, using \eqref{(2.9a)}, we obtain the formula
\begin{align*}
& \frac{1}{\Gamma(\alpha)}\int^t_{\xi} (t-s)^{\alpha-1}(s-\xi)^{\alpha-1}
\MLTWO(-\la_n (s-\xi)^{\alpha}) ds\\
=&  -\frac{1}{\la_n}(t-\xi)^{\alpha-1}\left(\MLTWO(-\la_n t^{\alpha}) 
- \frac{1}{\Gamma(\alpha)}\right).
\end{align*}
% \begin{align*}
% & \frac{1}{\Gamma(\alpha)}\int^t_{\xi} (t-s)^{\alpha-1}(s-\xi)^{\alpha-1}
% \MLTWO(-\la_n (s-\xi)^{\alpha}) ds\\
% = &\int^t_{\xi} \frac{1}{\Gamma(\alpha)}(t-s)^{\alpha-1}(s-\xi)^{\alpha-1}
% \sum_{k=0}^{\infty} \frac{(-\la_n)^k(s-\xi)^{\alpha k}}
% {\Gamma(\alpha k + \alpha)} ds\\
% =& \frac{1}{\Gamma(\alpha)}\sum_{k=0}^{\infty} 
% \frac{(-\la_n)^k}{\Gamma(\alpha k + \alpha)}
% \int^t_{\xi} (t-s)^{\alpha-1}(s-\xi)^{\alpha k+\alpha -1} ds\\
% = & \frac{1}{\Gamma(\alpha)}\sum_{k=0}^{\infty} 
% \frac{(-\la_n)^k}{\Gamma(\alpha k + \alpha)}
% \frac{\Gamma(\alpha)\Gamma(\alpha k + \alpha)}{\Gamma(\alpha k+2\alpha)}
% (t-\xi)^{\alpha k+2\alpha}\\
% =& -\frac{1}{\la_n}(t-\xi)^{\alpha-1} \sum_{j=1}^{\infty} 
% \frac{(-\la_n t^{\alpha})^j}{\Gamma(\alpha j + \alpha)}
% = -\frac{1}{\la_n}(t-\xi)^{\alpha-1}\left(\MLTWO(-\la_n t^{\alpha}) 
% - \frac{1}{\Gamma(\alpha)}\right).
% \end{align*}

Therefore, we have
\begin{align*}
& J^{\alpha}(L_nf)(t) = -\frac{1}{\la_n}(L_nf)(t)
+ \frac{1}{\la_n}\int^t_0 (t-\xi)^{\alpha-1}\frac{1}{\Gamma(\alpha)} 
f(\xi) d\xi\\
=& -\frac{1}{\la_n}(L_nf)(t) + \frac{1}{\la_n}(J^{\alpha}f)(t),
\quad n\in \N,
\end{align*}
that is,
$$
(L_nf)(t) = -\la_n J^{\alpha}(L_nf)(t) + (J^{\alpha}f)(t), \quad 0<t<T.
$$
Hence, $L_nf \in \HH(0,T) = J^{\alpha}L^2(0,T)$.  In view of   
$\pppa = (J^{\alpha})^{-1}$ (\cite{KRY}), we have 
$$
\pppa (L_nf) = -\la_n L_nf + f \quad \mbox{in }(0,T).
$$

Using the inequality \eqref{(5.4)}, we obtain
$$
\la_n\Vert L_nf\Vert_{L^2(0,T)} \le \la_n\Vert s^{\alpha-1}
\MLTWO(-\la_ns^{\alpha})\Vert_{L^1(0,T)}\Vert f\Vert_{L^2(0,T)}
\le \Vert f\Vert_{L^2(0,T)}.
$$
Therefore,
\begin{align*}
& \Vert L_nf\Vert_{\HH(0,T)} \le C\Vert \pppa(L_nf)\Vert_{L^2(0,T)}
\le C(\Vert -\la_nL_nf\Vert_{L^2(0,T)} + \Vert f\Vert_{L^2(0,T)})\\
\le& C\Vert f\Vert_{L^2(0,T)}, \quad n\in \N, \quad f\in L^2(0,T).
\end{align*}
Thus, the estimate \eqref{(5.5)} is proved.

Now we set $f_n(s) := (F(s), \, \va_n)$ for $0<s<T$ and $n\in \N$.
Since 
$$
\pppa \int^t_0 K(t-s)F(s) ds = \sumn \pppa(L_nf_n)(t)\va_n,
$$
we have
$$
\left\Vert \pppa \int^t_0 K(t-s)F(s) ds\right\Vert^2_{L^2(\OOO)}
= \sumn \vert \pppa(L_nf_n)(t)\vert^2.
$$ 
Then by applying \eqref{(5.5)} we obtain 
\begin{align*}
& \left\Vert \pppa \int^t_0 K(t-s)F(s) ds\right\Vert^2
_{\HH(0,T;L^2(\OOO))}
\le C\left\Vert \pppa\int^t_0 K(t-s)F(s) ds\right\Vert^2
_{L^2(0,T;L^2(\OOO))}\\
=& C\sumn \Vert \pppa(L_nf_n)\Vert^2_{L^2(0,T)}
\le C\sumn \Vert L_nf_n\Vert^2_{\HH(0,T)}\\
\le& C\sumn \Vert f_n\Vert^2_{L^2(0,T)}
= C\int^T_0 \sumn \vert (F(s),\va_n) \vert^2 ds\\
=& C\int^T_0 \Vert F(s)\Vert_{L^2(\OOO)}^2 ds 
= C\Vert F\Vert^2_{L^2(0,T;L^2(\OOO))}.
\end{align*}
Thus the proof of the estimate (ii) is completed.

The estimate (iii) in Lemma \ref{l5.1} follows from the standard estimates
of the operator $S(t)$ or can be derived by the same arguments as those that 
were employed in 
Section 6 of Chapter 4 in \cite{KRY} for the case of the 
homogeneous Dirichlet boundary condition and we omit the technical details. 
\end{proof}

Now we proceed with the
{\bf Proof of Theorem \ref{t2.1}}.

In the equation \eqref{(2.2)}, we regard the terms 
$\sum_{j=1}^d b_j(x,t)\ppp_ju$ and 
$(c_0+c(x,t))u$ as non-homogeneous terms and rewrite it in terms of 
the operator $A_0$ as follows
\begin{equation}
\label{(5.6)}
\left\{ \begin{array}{rl}
& \pppa (u-a) + A_0u(x,t) = F(x,t)\\
+& \sum_{j=1}^d b_j(x,t)\ppp_ju + (c_0+c(x,t))u, \quad 
x\in \OOO,\, 0<t<T,\\
& \NUNU u + \sigma(x) u = 0 \quad \mbox{on }\ppp\OOO \times (0,T),\\
& u(x,\cdot) - a(x) \in \HH(0,T) \quad \mbox{for almost all }x\in \OOO.
\end{array}\right.
\end{equation}
Then we represent the equation from \eqref{(5.6)} in the form  (\cite{GLY}, 
\cite{KRY})
$$
u(t) = S(t)a + \int^t_0 K(t-s)F(s) ds 
$$
\begin{equation}
\label{(5.7)}
+ \int^t_0 K(t-s) \left(\sum_{j=1}^d b_j(s)\ppp_ju(s) 
+ (c_0+c(s))u(s) \right) ds, \quad 0<t<T.        
\end{equation}
Moreover, it is known that if 
$u\in L^2(0,T;H^2(\OOO))$ satisfies the initial condition 
$u-a \in \HH(0,T;L^2(\OOO))$ and 
the equation \eqref{(5.7)}, then $u$ is a solution to the problem 
\eqref{(5.6)}.
With the notations
\begin{equation}
\label{(5.8)}
\left\{ \begin{array}{rl}
& G(t):=  \int^t_0 K(t-s)F(s) ds + S(t)a, \\
& Qu(t) = Q(t)u(t) := \sum_{j=1}^d b_j(\cdot,t)\ppp_ju(t) 
+ (c_0+c(\cdot,t))u(t), \\
& Ru(t):= \int^t_0 K(t-s) \left(\sum_{j=1}^d b_j(\cdot,s)\ppp_ju(s) 
+ (c_0+c(\cdot,s))u(s) \right) ds,\\
&\qquad \qquad \qquad \mbox{for $0<t<T$}, 
\end{array}\right.
\end{equation}
the equation \eqref{(5.7)} can be represented in form of a fixed point equation
$u = Ru + G$ in $L^2(0,T;H^2(\OOO))$.

Lemma \ref{l5.1} yields that $G \in L^2(0,T;H^2(\OOO))$.
Moreover, since $\Vert A_0^{\hhalf}a\Vert \le C\Vert a\Vert_{H^1(\OOO)}$
and $\DDD(A_0^{\hhalf}) = H^1(\OOO)$ (e.g., \cite{Fu}),
the estimate \eqref{(5.3)} implies
$$
\Vert S(t)a\Vert_{H^2(\OOO)} \le C\Vert A_0S(t)a\Vert
= C\Vert A_0^{\hhalf}S(t)A_0^{\hhalf}a\Vert
\le Ct^{-\hhalf\alpha}\Vert a\Vert_{H^1(\OOO)}
$$
and thus
$$
\Vert S(t)a\Vert^2_{L^2(0,T;H^2(\OOO))} 
\le C\left(\int^T_0 t^{-\alpha}dt \right)\Vert a\Vert^2_{H^1(\OOO)}
\le \frac{CT^{1-\alpha}}{1-\alpha}\Vert a\Vert^2_{H^1(\OOO)}.
$$
Consequently, the inclusion $S(t)a \in L^2(0,T;H^2(\OOO))$ holds valid.

For $0<t<T$, we next estimate $\Vert Rv(\cdot,t)\Vert_{H^2(\OOO)}$ for 
$v(\cdot,t) \in \mathcal{D}(A_0)$ as follows:
\begin{align*}
& \Vert Rv(\cdot,t)\Vert_{H^2(\OOO)}
\le C\Vert A_0Rv(\cdot,t)\Vert_{L^2(\OOO)}\\
\le & \int^t_0 \left\Vert A_0^{\hhalf}K(t-s)A_0^{\hhalf}
\left(\sum_{j=1}^d b_j(s)\ppp_ju(s) + (c_0+c(s))u(s) \right)
\right\Vert ds\\
\le & C\int^t_0 \Vert A_0^{\hhalf}K(t-s)\Vert 
\left\Vert A_0^{\hhalf}
\left(\sum_{j=1}^d b_j(s)\ppp_ju(s) + (c_0+c(s))u(s) \right)
\right\Vert ds\\
\le& C\int^t_0 (t-s)^{\hhalf\alpha -1}\Vert v(s)\Vert_{H^2(\OOO)}ds
= C\left( \Gamma\left(\hhalf\alpha\right)J^{\hhalf\alpha}\Vert v\Vert
_{H^2(\OOO)}\right)(t).
\end{align*}
For derivation of this estimate, we employed the inequalities
$$
\Vert A_0^{\hhalf}b_j(s)\ppp_jv(t)\Vert 
\le C\Vert b_j(s)\ppp_jv(s)\Vert_{H^1(\OOO)}
\le C\Vert v(s)\Vert_{H^2(\OOO)},
$$
$$
\Vert (c(s)+c_0)v(s)\Vert_{H^1(\OOO)} \le C\Vert v(s)\Vert_{H^2(\OOO)}
$$
that hold true because $b_j \in C^1(\ooo{\OOO}\times [0,T])$) and 
$c+c_0\in C([0,T];C^1(\ooo{\OOO}))$.

Since $(J^{\hhalf\alpha}w_1)(t) \ge (J^{\hhalf\alpha}w_2)(t)$ 
if $w_1(t) \ge w_2(t)$ for $0\le t\le T$, and
$J^{\hhalf\alpha}J^{\hhalf\alpha}w = J^{\alpha}w$ for 
$w_1, w_2, w \in L^2(0,T)$, we have
\begin{align*}
&\Vert R^2v(t)\Vert_{H^2(\OOO)} = \Vert R(Rv)(t)\Vert_{H^2(\OOO)}\\
\le & C\left( \Gamma\left(\hhalf\alpha\right)J^{\hhalf\alpha}
\left(C\Gamma\left(\hhalf\alpha\right)J^{\hhalf\alpha}
\Vert v\Vert_{H^2(\OOO)}\right) \right)(t)\\
= & \left( C\Gamma\left(\hhalf\alpha\right)\right)^2
(J^{\alpha}\Vert v\Vert_{H^2(\OOO)})(t).
\end{align*}
Repeating this argumentation $m$-times, we obtain
\begin{align*}
& \Vert R^mv(t)\Vert_{H^2(\OOO)}
\le \left( C\Gamma\left(\hhalf\alpha\right)\right)^m
\left( J^{\hhalf\alpha m}\Vert v\Vert_{H^2(\OOO)}\right)(t)\\
\le & \frac{\left( C\Gamma\left(\hhalf\alpha\right)\right)^m}
{\Gamma\left( \hhalf\alpha m\right)}
\int^t_0 (t-s)^{\frac{m}{2}\alpha -1} \Vert v(\xi)\Vert_{H^2(\OOO)}ds, \quad 
0<t<T.
\end{align*}
Applying the Young inequality to the integral in the last estimate, 
we reach the inequality
\begin{align*}
& \Vert R^mv(t)\Vert_{L^2(0,T; H^2(\OOO))}^2
\le \left(  \frac{\left( C\Gamma\left(\hhalf\alpha\right) \right)^m}
{\Gamma\left( \hhalf\alpha m\right)}\right)^2
\Vert t^{\frac{\alpha m}{2}-1}\Vert_{L^1(0,T)}^2 
\Vert v\Vert_{L^2(0,T;H^2(\OOO))}^2\\
=& \frac{\left( CT^{\frac{\alpha}{2}}
\Gamma\left(\hhalf\alpha \right)\right)^{2m}}
{\Gamma\left( \hhalf\alpha m +1\right)^2}
\Vert v\Vert_{L^2(0,T;H^2(\OOO))}^2.
\end{align*}
Employing the known asymptotic behavior of the gamma function, we obtain the 
relation
$$
\lim_{m\to\infty} \frac{\left( CT^{\frac{\alpha}{2}}
\Gamma\left(\hhalf\alpha \right)\right)^m}
{\Gamma\left( \hhalf\alpha m +1\right)} = 0
$$
that means that for sufficiently large $m\in \N$, the mapping
$$
R^m: L^2(0,T;H^2(\OOO))\, \longrightarrow \, L^2(0,T;H^2(\OOO))
$$ 
is a contraction.  Hence, by the Banach fixed point theorem, the equation 
\eqref{(5.7)} possesses a unique fixed point.
Therefore, by the first equation in \eqref{(2.3)}, we reach the inclusion
$\pppa (u-a) \in L^2(0,T;L^2(\OOO))$.
Since $\Vert \eta\Vert_{\HH(0,T)} \sim \Vert \pppa \eta\Vert_{L^2(0,T)}$
for $\eta \in \HH(0,T)$ (\cite{KRY}), we finally obtain the estimate
$$
\Vert u-a\Vert_{\HH(0,T;L^2(\OOO))} + \Vert u\Vert_{L^2(0,T;H^2(\OOO))}
\le C(\Vert a\Vert_{H^1(\OOO)} + \Vert F\Vert_{L^2(0,T;L^2(\OOO))}).
$$
The proof of Theorem \ref{t2.1} is completed.
%$\square$
\end{proof}

%%%%%%%%%%% Section 3 %%%%%%%%%%%%%%%
\section{Key lemma} 
\label{sec3}

\setcounter{section}{3}

For derivation of the comparison principles for solutions to the linear and 
semilinear time-fractional diffusion equations, we need some auxiliary results 
that are formulated and proved in this section.

%The proof is based on the representation of solution based on the 
%principal part of the linear equation in (1.1) and as 

In addition to the operator $-A_0$ defined by \eqref{(3.2)}, 
we define an elliptic operator $-A_1$ with a positive zeroth-order
coefficient:
\begin{equation}\label{(3.1a)}
(-A_1(t)v)(x):= (-A_1v)(x) 
\end{equation}
$$
:= \sumij \ppp_i(a_{ij}(x)\ppp_jv(x)) 
+ \sum_{j=1}^d b_j(x,t)\ppp_jv - b_0(x,t)v,
$$
where $b_0(x,t) > 0$, $\in C^1([0,T];C^1(\ooo{\OOO})) \cap 
C([0,T];C^2(\ooo{\OOO}))$, and \\
$\min_{(x,t)\in \ooo{\OOO}\times 
[0,T]} b_0(x,t)$ is sufficiently large.
We recall that the pointwise Caputo derivative $\ddda$ is defined by 
\begin{equation}
\label{(4.1)}
\ddda y(t) = \frac{1}{\Gamma(1-\alpha)}
\int^t_0 (t-s)^{-\alpha}\frac{dy}{ds}(s) ds \quad
\mbox{for $y\in W^{1,1}(0,T)$}.                      
\end{equation}

Henceforth we write $y'(t) = \frac{dy}{dt}(t)$ if there is no fear of 
confusion.

In what follows, we employ an extremum principle for the Caputo fractional 
derivative formulated below.

\begin{lemma}[\cite{Lu1}]
\label{l4.1}
We assume that $y\in C[0,T]$ and $t^{1-\alpha}y' \in C[0,T]$.  
If $y=y(t)$ attains its minimum over the interval 
$[0,T]$ at the point $t_0 \in (0, \,T]$, then 
$$
\ddda y(t_0) \le 0.
$$
\end{lemma}

In the lemma, the assumption $t_0>0$ is essential.

In \cite{Lu1}, Lemma \ref{l4.1} was formulated and proved under a 
weaker 
regularity condition posed on the function $y$, but for our arguments  
we can assume that $y\in C[0,T]$ and 
$t^{1-\alpha}y' \in C[0,T]$.

Employing Lemma \ref{l4.1}, we now prove our key lemma that is a basis for 
further arguments in this article.

\begin{lemma}[Positivity of a smooth solution]
\label{l4.2}
We assume that \\
$\min_{(x,t)\in \ooo{\OOO}\times 
[0,T]} b_0(x,t)$ is a sufficiently large positive constant.
For $a\in H^1_0(\OOO)$ and $F\in L^2(0,T;L^2(\OOO))$,
we assume that there exists a solution $u\in C([0,T];C^2(\ooo{\OOO}))$
satisfying $t^{1-\alpha}\ppp_tu \in C([0,T];C(\ooo{\OOO}))$ 
to the initial-boundary value problem
\begin{equation}
\label{(4.2)}
\left\{ \begin{array}{rl}
& \pppa (u-a) + A_1u = F(x,t), \quad x\in \OOO,\, 0<t<T, \\
& \NUNU u + \sigma(x)u = 0 \quad \mbox{on } \ppp\OOO\times (0,T),\\
& u(x,\cdot) - a(x) \in \HH(0,T) \quad \mbox{for almost all } x\in \OOO.
\end{array}\right.
\end{equation}
If $F \ge 0$ in $\OOO\times (0,T)$ and $a\ge 0$ in $\OOO$, then 
$$
u \ge 0 \quad \mbox{in } \OOO\times (0,T).
$$
\end{lemma}

We note that the regularity of the solution to the problem \eqref{(4.2)} 
at the point $t=0$ is a more delicate question  compared to the case 
$\alpha=1$. In particular, we can not expect that $u(x,\cdot) \in 
C^1[0,T]$.  This can be illustrated by a simple example of the equation 
$\pppa y(t) = y(t)$ with $y(t)-1 \in \HH(0,T)$ whose unique solution 
$y(t) = \MLONE(t^{\alpha})$ does not belong to the space $C^1[0,T]$.

For the case of $\alpha=1$, a similar positivity property for the 
time-fractional diffusion equation under the Robin boundary condition is 
well-known.

\begin{proof}
First we introduce an auxiliary function $\psi \in C^1([0,T];C^2(\ooo{\OOO}))$ 
that satisfies the conditions
\begin{equation}
\label{(4.3)}
\left\{ \begin{array}{rl}
& A_1\psi(x,t) = 1, \quad (x,t) \in \OOO\times [0,T], \\
& \NUNU \psi + \sigma \psi = 1 \quad \mbox{on } \ppp\OOO\times [0,T].
\end{array}\right.
\end{equation}
The proof of the existence of such $\psi$ is postponed to Appendix.

Now, choosing $M>0$ sufficiently large and $\ep>0$ sufficiently small,
we set 
\begin{equation}
\label{(w_u)}
w(x,t) := u(x,t) + \ep(M + \psi(x,t) + t^{\alpha}), \quad x\in \OOO,\,
0<t<T.
\end{equation}

For a fixed $x\in \OOO$, by the assumption on the regularity 
of $u$, we have
\begin{equation}\label{(3.6)}
t^{1-\alpha}\ppp_tu(x,\cdot) \in C[0,T].
\end{equation}
Then, $\ppp_tu(x,\cdot) \in L^1(0,T)$, that is, $u(x,\cdot)
\in W^{1,1}(0,T)$.  Moreover,
\begin{equation}\label{(3.7)}
u(x,0) - a(x) = 0, \quad x\in \OOO.
\end{equation}
           
On the other hand, we have
$$
\pppa w = \ddda w = \ddda (w+c)
$$
with any constant $c$, if $w\in H_{\alpha}(0,T) \cap W^{1,1}(0,T)$ 
and $w(0) = 0$
(e.g., Theorem 2.4 of Chapter 2 in \cite{KRY})

Since $u(x,\cdot) - a \in H_{\alpha}(0,T)$ and $u(x,\cdot) \in W^{1,1}(0,T)$ 
for almost all $x\in \OOO$, by \eqref{(3.7)} we see
$\pppa (u-a) = \ddda (u-a) = \ddda u$.

Furthermore, since $\ep(M+\psi(\cdot,t)+t^{\alpha}) \in W^{1,1}(0,T)$, 
we obtain
\begin{align*}
& \ddda w = \ddda (u + \ep(M+\psi(x,t)+t^{\alpha})))
= \ddda u + \ep\ddda (M+\psi(x,t)+t^{\alpha})\\
=& \pppa (u-a) + \ep(\ddda (\psi + t^{\alpha}))
= \pppa (u-a) + \ep(\ddda\psi + \Gamma(\alpha+1))
\end{align*}
and
\begin{align*}
& A_1w = A_1u + \ep A_1\psi + \ep A_1t^{\alpha} + \ep A_1M\\
=& A_1u + \ep + \ep b_0(x,t)t^{\alpha} + b_0(x,t)\ep M.
\end{align*}
 
Now we choose a constant $M>0$ such that $M + \psi(x,t) \ge 0$ and
$\ddda \psi(x,t) + b_0(x,t)M > 0$ for $(x,t) \in \ooo{\OOO} \times [0,T]$, so 
that 
\begin{equation}
\label{(4.4)}
\ddda w + A_1w 
\end{equation}
$$
= F + \ep(\Gamma(\alpha+1) + \ddda \psi + 1 
+ b_0(x,t)t^{\alpha} + b_0(x,t) M) > 0 \quad \mbox{in } \OOO\times (0,T).
$$
Moreover, because of the relation $\NUNU w = \NUNU u + \ep \NUNU \psi$, 
we obtain the following estimate:
$$
\NUNU w + \sigma w = \NUNU u + \sigma u + \ep + \sigma \ep t^{\alpha}
+ \ep\sigma M
$$
\begin{equation}
\label{(4.5)}
\ge \ep + \sigma\ep t^{\alpha} + \ep\sigma M \ge \ep \quad 
\mbox{on } \ppp\OOO\times (0,T).
\end{equation}
Evaluating the representation \eqref{(w_u)} at the point $t=0$ immediately leads 
to the formula
$$
%\begin{equation}
%\label{(4.6)}
w(x,0) = u(x,0) + \ep(\psi(x,0) + M), \quad x\in \OOO.
$$
%\end{equation}

Let us assume that the inequality 
$$
\min_{(x,t)\in \ooo{\OOO}\times [0,T]} w(x,t) \ge 0
$$
does not hold, that is, there exists
$(x_0,t_0) \in \ooo{\OOO}\times [0,T]$ such that 
\begin{equation}
\label{(4.7)}
w(x_0,t_0):= \min_{(x,t)\in \ooo{\OOO}\times [0,T]} w(x,t) < 0.
\end{equation}
Since $M>0$ is sufficiently large, by $u(x,0) \ge 0$, we have 
$$
w(x,0) = u(x,0) + \ep(\psi(x,0) + M) \ge u(x,0) \ge 0, \quad 
x\in \ooo{\OOO},
$$
and thus $t_0\ne 0$ has to be valid.

Next we show that $x_0 \not\in \ppp\OOO$.  Indeed, let us assume that 
$x_0 \in \ppp\OOO$.  Then the estimate \eqref{(4.5)} yields that 
$\NUNU w(x_0,t_0) + \sigma(x_0)w(x_0,t_0) \ge \ep$.
By \eqref{(4.7)} and $\sigma(x_0)\ge 0$, we obtain 
$$
\NUNU w(x_0,t_0) \ge -\sigma(x_0)w(x_0,t_0) + \ep 
\ge \ep > 0,
$$
which implies
$$
\ppp_{\nu_A}w(x_0,t_0) = 
\sumij a_{ij}(x_0)\nu_j(x_0)\ppp_iw(x_0,t_0)
= \nabla w(x_0,t_0) \cdot \mathcal{A}(x_0)\nu(x_0)
$$
\begin{equation}
\label{(4.8)}
= \sum_{i=1}^d (\ppp_iw)(x_0,t_0)[\mathcal{A}(x_0)\nu(x_0)]_i > 0.
\end{equation}
Here $\mathcal{A}(x) = (a_{ij}(x))_{1\le i,j\le d}$ and 
$[b]_i$ means the $i$-th element of a vector $b$.

For sufficiently small $\ep_0>0$ and $x_0\in \ppp\OOO$, we now verify 
\begin{equation}
\label{(4.9)}
x_0 - \ep_0\mathcal{A}(x_0)\nu(x_0) \in \OOO.             
\end{equation}
Indeed, since the matrix 
$\mathcal{A}(x_0)$ is positive-definite, the inequality  
$$
(\nu(x_0)\, \cdot \, -\ep_0\mathcal{A}(x_0)\nu(x_0))
= -\ep_0(\mathcal{A}(x_0)\nu(x_0)\, \cdot \, \nu(x_0)) < 0
$$
holds true. In other words, we have 
$$
\angle (\nu(x_0), \, (x_0 - \ep_0\mathcal{A}(x_0)\nu(x_0)) - x_0)
> \frac{\pi}{2}.
$$
Because the boundary $\ppp\OOO$ is smooth, the domain $\OOO$ is locally 
located on 
one side of $\ppp\OOO$. In a small neighborhood  of the point 
$x_0\in \ppp\OOO$, the boundary $\ppp\OOO$ can 
be described in the local coordinates composed of 
its tangential component in $\R^{d-1}$
and the normal component along $\nu(x_0)$.  
Consequently, if $y \in \R^d$ satisfies 
$\angle (\nu(x_0), y-x_0) > \frac{\pi}{2}$, then $y\in \OOO$.
Therefore, for a sufficiently small $\ep_0>0$, 
we have $x_0-\ep_0\mathcal{A}(x_0)\nu(x_0) \in \OOO$ and 
we proved the formula \eqref{(4.9)}. 
% $\square$

Moreover, for sufficiently small $\ep_0>0$, we can prove that
\begin{equation}\label{(3.11a)}
w(x_0 - \ep_0\mathcal{A}(x_0)\nu(x_0),\,t_0) < w(x_0,t_0).
\end{equation}
{\bf Verification of \eqref{(3.11a)}.}
\\
Indeed, the inequality \eqref{(4.7)} yields
$$
\sum_{i=1}^d (\ppp_iw)(x_0-\eta\mathcal{A}(x_0)\nu(x_0),\,t_0)
[\mathcal{A}(x_0)\nu(x_0)]_i > 0 \quad\mbox{if }\vert \eta\vert < \ep_0.
$$
Then, by the mean value theorem, we obtain the inequality 
\begin{align*}
&w(x_0 - \xi\mathcal{A}(x_0)\nu(x_0),\,t_0) - w(x_0,t_0)\\
= & \xi\sum_{i=1}^d \ppp_iw(x_0 - \theta\mathcal{A}(x_0)\nu(x_0),\,t_0) 
(-[\mathcal{A}(x_0)\nu(x_0)]_i) < 0,
\end{align*}
where $\theta$ is a number between $0$ and $\xi\in (0,\ep_0)$.
Thus \eqref{(3.11a)} is verified. $\square$
\\

By combining \eqref{(3.11a)} with \eqref{(4.9)}, we 
conclude that there exists $\www{x_0} \in \OOO$ such that 
$w(\www{x_0},t_0) < w(x_0,t_0)$, which contradicts the assumption 
\eqref{(4.7)} and 
we have proved that $x_0 \not\in \ppp\OOO$.

By the definition, the function $w$ attains its minimum 
at the point $(x_0,t_0)$. Because $0 < t_0 \le T$, 
Lemma \ref{l4.1} yields the inequality 
\begin{equation}
\label{(4.10)}
\ddda w(x_0,t_0) \le 0.                 
\end{equation}
Since $x_0 \in \OOO$, the necessary condition for an extremum point leads to the equality
\begin{equation}
\label{(4.11)}
\nabla w(x_0,t_0) = 0.                      
\end{equation}
Moreover, because $w$ attains its minimum at the point $x_0 \in \OOO$, 
in view of the sign of the Hessian, the inequality
\begin{equation}
\label{(4.12)}
\sumij a_{ij}(x_0)\ppp_i\ppp_j w(x_0,t_0) \ge 0   
\end{equation}
holds true (see e.g. the proof of Lemma 1 in Section 1 of Chapter 2 in 
\cite{Fr}).

The inequalities $b(x_0,t_0)>0$, $w(x_0,t_0) < 0$, and
\eqref{(4.10)}-\eqref{(4.12)} let to derive the inequality 
\begin{align*}
& \ddda w(x_0,t_0) + A_1w(x_0,t_0)\\
= & \ddda w (x_0,t_0) - \sumij a_{ij}(x_0)\ppp_i\ppp_jw(x_0,t_0)
- \sum_{i=1}^d (\ppp_ia_{ij})(x_0)\ppp_jw(x_0,t_0)\\
-& \sum_{i=1}^d b_i(x_0,t_0)\ppp_iw(x_0,t_0) + b(x_0,t_0)w(x_0,t_0) < 0,
\end{align*}
which contradicts the inequality \eqref{(4.4)}.

We thus have proved that 
$$
u(x,t) + \ep(M+\psi(x,t)+t^{\alpha}) = w(x,t) \ge 0, \quad 
(x,t) \in \OOO\times (0,T).
$$
Since $\ep>0$ is arbitrary, we let $\ep \downarrow 0$ to obtain the inequality 
$u(x,t) \ge 0$ for $(x,t) \in \OOO\times (0,T)$ and 
 the proof of Lemma \ref{l4.2} is completed. 
\end{proof}

Let us finally mention that the positivity of the function $b_0$ 
in the operator $-A_1$ 
is an essential condition for validity of our proof of Lemma \ref{l4.2}. 
However, in the next 
section, we remove this condition while deriving the comparison principle for 
the solutions to the linear time-fractional diffusion equations.
%
%
%
%%%%%%%%%%% Section 4 %%%%%%%%%%%%%%%
\section{Comparison principles for the linear time-fractional diffusion equations} 
\label{sec4}

\setcounter{section}{4}

According to the results formulated in Theorem \ref{t2.1}, in this section, we 
consider the solutions to the initial-boundary value problem \eqref{(2.3)} 
within the class 
\begin{equation}\label{(4.1a)}
\{ u; \, u-a\in \HH(0,T;L^2(\OOO)), \, u\in L^2(0,T;H^2(\OOO))\}.
\end{equation}
By $a\ge 0$, we always mean that the inequality $a\ge 0$ holds almost 
everywhere in a set under consideration.

Our first main result concerning the comparison
principle for the linear time-fractional diffusion equations is as follows:

\begin{theorem}
\label{t2.2}
Let $a \in H^1_0(\OOO)$ and $F \in L^2(\OOO\times (0,T))$.
If $F(x,t) \ge 0$ in $\OOO\times (0,T)$ and
$a(x) \ge 0$ in $\OOO$, then
$u(F,a)(x,t) \ge 0$ in $\OOO\times (0,T)$,
where by $u(F,a)$ we denote the solution $u$ in the class \eqref{(4.1a)}
to the initial-boundary value problem \eqref{(2.3)} in the class \eqref{(4.1a)}
with the initial condition $a$ and the right-hand side $F$. 
\end{theorem}
%%%%%%%%%%%%%%%%
We emphasize that the non-negativity of solution $u$ for \eqref{(2.3)} holds for 
general class \eqref{(4.1a)}, and $u$ does not necessarily satisfy 
$u \in C([0,T];C^2(\ooo{\OOO}))$ and $t^{1-\alpha}\ppp_tu \in 
C([0,T]; C(\ooo{\OOO}))$.  Thus Theorem \ref{t2.2} is widely applicable.  

Theorem \ref{t2.2} immediately yields the following comparison property:
\begin{corollary}
\label{c2.1}
Let $a_1, a_2  \in H^1_0(\OOO)$ and $F_1, F_2 \in L^2(\OOO\times (0,T))$
satisfy $a_1(x) \ge a_2(x)$ in $\OOO$ and 
$F_1(x,t) \ge F_2(x,t)$ in $\OOO\times (0,T)$, respectively.  
Then $u(F_1, a_1)(x,t) \ge u(F_2,a_2)(x,t)$ in $\OOO\times (0,T)$.
\end{corollary}

\begin{proof}
We first prove Corollary \ref{c2.1}.
Setting $a:= a_1-a_2$, $F:= F_1 - F_2$ and $u:= u(F_1,a_1) - u(F_2,a_2)$, 
we obtain that $a\ge 0$ in $\OOO$ and $F\ge 0$ in $\QQQQ$ and 
$$
\left\{ \begin{array}{rl}
& \pppa (u-a) + Au = F \ge 0 \quad \mbox{in } \QQQQ, \\
& \NUNU u + \sigma u = 0 \quad \mbox{on } \ppp\OOO.
\end{array}\right.
$$
Therefore, Theorem \ref{t2.2} implies that $u\ge 0$, that is,
$u(F_1,a_1) \ge u(F_2,a_2)$ in $\QQQQ$.
\end{proof}

We apply Corollary \ref{c2.1} to derive lower and upper bounds for the solution 
$u$ by choosing initial values and non-homogeneous terms suitably.
Here we are restricted to demonstrate only one example.
\begin{example}
\label{ex1}
Let 
$$
-Av(x) = \sum_{i,j=1}^d \ppp_i(a_{ij}(x)\ppp_jv(x)) + \sum_{j=1}^d
b_j(x,t)\ppp_jv(x),
$$
where the coefficients $a_{ij}, b_j$, $1\le i,j\le d$ satisfy \eqref{(1.2)}.
For constants $\beta \ge 0$ and $\delta>0$, we assume that
$a=0$ in $\OOO$ and $F \in L^2(0,T;L^2(\OOO))$ satisfies
$$
F(x,t) \ge \delta t^{\beta}, \quad x\in \OOO,\, 0<t<T.
$$
Then 
\begin{equation}\label{(4.2a)}
u(F,0)(x,t) \ge \frac{\delta\Gamma(\beta+1)}{\Gamma(\alpha+\beta+1)}
t^{\alpha+\beta}, \quad x\in \OOO, \, 0\le t \le T.
\end{equation}
Indeed, setting 
$$
\underline{u}(x,t):= \frac{\delta\Gamma(\beta+1)}{\Gamma(\alpha+\beta+1)}
t^{\alpha+\beta}, \quad x\in \OOO, \, t>0,
$$
we see that 
$$
\left\{\begin{array}{rl}
& \pppa \underline{u} + A\underline{u} = \delta t^{\beta} \quad \mbox{in $\OOO 
\times (0,T)$}, \\
& \ppp_{\nu_A}\underline{u} = 0 \quad \mbox{on $\ppp\OOO \times 
(0,T)$}, \\
& \underline{u}(x,\cdot) \in H_{\alpha}(0,T).
\end{array}\right.
$$
By $F \ge \delta t^{\beta}$ in $\OOO \times (0,T)$, we apply Corollary \ref{c2.1}
to $u$ and $\underline{u}$, and we can conclude \eqref{(4.2a)}.
%$\square$

In particular, noting that $u \in L^2(0,T;H^2(\OOO)) \subset 
L^2(0,T;C(\ooo{\OOO}))$ by the Sobolev embedding and the spatial 
dimensions $d \le 3$, we see that $u(F,0)(x,t) > 0$ for almost all
$t>0$ and all $x\in \ooo{\OOO}$.
\end{example}

Now we proceed to the proof of Theorem \ref{t2.2}.
\begin{proof}
We use the operators $Qu(t)$ and $G(t)$ defined by \eqref{(5.8)}.
In terms of these operators, the solution $u(t):= u(F,a)(t)$ satisfies 
the equation
\begin{equation}
\label{(6.1)}
u(F,a)(t) = G(t) + \int^t_0 K(t-s)Qu(s) ds, \quad 0<t<T.   
\end{equation}

We divide the proof of Theorem \ref{t2.2} into three parts.

{\bf (I) First part of the proof of Theorem \ref{t2.2}: existence of 
smoother solution.}

In the formulation of Lemma \ref{l4.2}, 
we assumed the existence of a solution $u$  to
the initial-boundary value problem \eqref{(4.2)} that satisfies the inclusions $u\in C([0,T];C^2(\ooo{\OOO}))$
and  $t^{1-\alpha}\ppp_tu \in C([0,T];C(\ooo{\OOO}))$.
On the other hand, Theorem \ref{t2.1} asserts the unique existence  
of solution $u$ to the initial-boundary value problem \eqref{(2.3)}
in a class $u(F,a) \in L^2(0,T;H^2(\OOO))$ and $u(F,a) - a 
\in H_{\alpha}(0,T;L^2(\OOO))$.  

In this part of the proof, we lift up the regularity of the solution
$u\in L^2(0,T;H^2(\OOO)) \cap (\HH(0,T;L^2(\OOO)) + \{ a\})$ of the problem \eqref{(4.2)} 
with $a \in C^{\infty}_0(\OOO)$ and $F\in C^{\infty}_0(\OOO\times 
(0,T))$  to the inclusion 
$u \in C([0,T];C^2(\ooo{\OOO}))$ satisfying the condition
$t^{1-\alpha}\ppp_tu \in C([0,T];C(\ooo{\OOO}))$.

More precisely, we can state

\begin{lemma}
\label{l6.1}
Let $a_{ij}$, $b_j$, $c$ satisfy the condition \eqref{(1.2)} and 
the inclusions 
$a\in C^{\infty}_0(\OOO)$, $F \in C^{\infty}_0(\OOO\times (0,T))$ hold true.
Then the solution $u=u(F,a)$ to the problem \eqref{(2.3)} satisfies 
$$
u \in C([0,T];C^2(\ooo{\OOO})), \quad 
t^{1-\alpha}\ppp_tu \in C([0,T];C(\ooo{\OOO})).
$$
\end{lemma}

\begin{proof}
We recall that $c_0>0$ is a positive fixed constant and 
$$
-A_0v = \sumij \ppp_i(a_{ij}(x)\ppp_jv) - c_0v,\ \DDD(A_0) = \{ v \in H^2(\OOO);\, \NUNU v + \sigma v = 0 \,\,
\mbox{on } \ppp\OOO\}.$$
Then $\DDD(A_0^{\hhalf}) = H^1(\OOO)$ and
$\Vert A_0^{\hhalf}v\Vert \sim \Vert v\Vert_{H^1(\OOO)}$
(\cite{Fu}).  Moreover, the estimates \eqref{(5.3)} hold for the operators 
$S(t)$ and 
$K(t)$ defined by \eqref{(5.1)} and \eqref{(5.2)}.

Here we write $u'(t) = \frac{du}{dt}(t) = \frac{\ppp u}{\ppp t}(\cdot,t)$
if there is no fear of confusion.

The relation \eqref{(6.1)} shows  that its solution $u(t)$ can be constructed 
as a fixed point of the equation
\begin{equation}
\label{(6.2)}
A_0u(t) = A_0G(t) + \int^t_0 A_0^{\hhalf}K(t-s)A_0^{\hhalf}Qu(s) ds,
\quad 0<t<T.                                 
\end{equation}
As already proved, the fixed point $u$ satisfies  
$u\in L^2(0,T;H^2(\OOO)) \cap (\HH(0,T;L^2(\OOO)) + \{ a\})$.

% Henceforth we write
% $$
% u'(t):= \frac{du}{dt}(t), \,\, \mbox{etc.}
% $$
We start with the derivation of some estimates for 
$\Vert A_0^{\kappa}u(t)\Vert$,
$\kappa=0,1$ and $\Vert A_0u'(t)\Vert$ for $0<t<T$.
We set 
$$
%\begin{equation}
%\label{(6.3)}
D:= \sup_{0<t<T} (\Vert A_0F(t)\Vert + \Vert A_0F'(t)\Vert 
+ \Vert A_0^2F(t)\Vert) + \Vert a\Vert_{H^4(\OOO)}.       
$$
%\end{equation}
Since $F\in C^{\infty}_0(\OOO\times (0,T))$, we see  
$F\in L^{\infty}(0,T;\DDD(A_0^2))$ and $D < \infty$.
Moreover, in view of \eqref{(5.3)}, for $\kappa=1,2$, we can 
estimate 
\begin{align*}
& \left\Vert A_0^{\kappa}\int^t_0 K(t-s)F(s) ds \right\Vert
\le C\int^t_0 \Vert K(t-s)\Vert \Vert A_0^{\kappa}F(s)\Vert ds\\
\le& C\left( \int^t_0 (t-s)^{\alpha-1}ds \right) 
\sup_{0<s<T} \Vert A_0^{\kappa}F(s)\Vert \le CD,
\end{align*}
\begin{align*}
& \left\Vert A_0\frac{d}{dt}\int^t_0 K(t-s)F(s) ds \right\Vert
= \left\Vert A_0\frac{d}{dt}\int^t_0 K(s)F(t-s) ds \right\Vert\\
=& \left\Vert A_0K(t)F(0) + A_0\int^t_0 K(s)F'(t-s) ds \right\Vert\\
\le& C\left\Vert A_0\int^t_0 K(s)F'(t-s) ds \right\Vert
\le C \int^t_0 s^{\alpha-1}
\Vert A_0F'(t-s) \Vert ds < CD, \quad \mbox{etc.}
\end{align*}
The regularity conditions in \eqref{(1.2)} lead to the estimates
$$
\Vert A_0^{\hhalf}Q(s)u(s)\Vert \le C\Vert Q(s)u(s)\Vert_{H^1(\OOO)}
$$
$$
= C\left\Vert \sum_{j=1}^d b_j(s)\ppp_ju(s) + (c_0+c(s))u(s)
\right\Vert_{H^1(\OOO)}
$$
\begin{equation}
\label{(6.4)}
\le C\Vert u(s)\Vert_{H^2(\OOO)} \le C\Vert A_0u(s)\Vert,\quad 
0<s<T.      
\end{equation}
Moreover, 
$$
\Vert A_0S(t)a\Vert = \Vert S(t)A_0a\Vert \le C\Vert a\Vert_{H^2(\OOO)}
\le CD
$$
by \eqref{(5.3)}.  Then 
\begin{align*}
& \Vert A_0u(t)\Vert \le CD 
+ \int^t_0 \Vert A_0^{\hhalf}K(t-s)\Vert \Vert A_0^{\hhalf}Q(s)u(s)\Vert ds\\
\le& CD + C\int^t_0 (t-s)^{\hhalf\alpha -1}\Vert A_0u(s)\Vert ds, \quad
0<s<T.
\end{align*}
The generalized Gronwall inequality yields 
$$
\Vert A_0u(t)\Vert \le CD + C\int^t_0 (t-s)^{\hhalf\alpha -1}D ds
\le CD, \quad 0<t<T,
$$
which implies
$$
\Vert A_0u\Vert_{L^{\infty}(0,T;H^2(\OOO))} \le CD. 
$$
Moreover, we can repeat the same arguments in the space
$C([0,T]; L^2(\OOO))$ as the ones employed for 
the iterations $R^n$ of the operator in the proof 
of Theorem \ref{t2.1} and apply the fixed point theorem to the  
equation \eqref{(6.1)}, so that $A_0u \in C([0,T];L^2(\OOO))$.
Therefore we reach  
\begin{equation}
\label{(6.5)}
u\in C([0,T];H^2(\OOO)), \quad 
\Vert u\Vert_{C([0,T];H^2(\OOO))}
\le CD.                               
\end{equation}
Choosing $\ep_0 > 0$ sufficiently small, we have
\begin{equation}
\label{(6.6)}
A_0^{\frac{3}{2}}u(t) = A_0^{\frac{3}{2}}G(t) 
+ \int^t_0 A_0^{\frac{3}{4}+\ep_0} K(t-s)A_0^{\frac{3}{4}-\ep_0}Q(s)u(s)ds,
\quad 0<t<T.                                                
\end{equation}

Next, according to \cite{Fu}, the inclusion
$$
\DDD(A_0^{\frac{3}{4}-\ep_0}) \subset H^{\frac{3}{2}-2\ep_0}(\OOO)
$$
holds true.
Now we proceed to the proof of the inclusion $Q(s)u(s) \in 
\DDD(A_0^{\frac{3}{4}-\ep_0})$.  
By \eqref{(5.3)}, we obtain the inequality
$$
\Vert A_0^{\frac{3}{2}}u(t)\Vert \le CD
+ \int^t_0 (t-s)^{(\frac{1}{4}-\ep_0)\alpha-1} 
\Vert A_0^{\frac{3}{4}-\ep_0}Q(s)u(s)\Vert ds,
$$
which leads to the estimate
$$
\Vert u(t)\Vert_{H^3(\OOO)} \le CD
+ \int^t_0 (t-s)^{(\frac{1}{4}-\ep_0)\alpha-1} 
\Vert u(s)\Vert_{H^3(\OOO)} ds, \quad 0<t<T.
$$
In the last estimate, we employed the inequality
$$
\Vert A_0^{\frac{3}{4}-\ep_0}Q(s)u(s)\Vert 
\le C\Vert Q(s)u(s)\Vert_{H^{\frac{3}{2}}(\OOO)}
\le C\Vert Q(s)u(s)\Vert_{H^2(\OOO)}
\le C\Vert u(s)\Vert_{H^3(\OOO)},
$$
which follows from by the regularity conditions \eqref{(1.2)} posed on the 
coefficients $b_j, c$. 

The generalized Gronwall inequality yields the estimate 
$$
\Vert u(t)\Vert_{H^3(\OOO)} \le C\left(1 
+ t^{\alpha\left(\frac{1}{4}-\ep_0\right)} 
\right)D
$$
for $0<t<T$.         

For the relation \eqref{(6.6)}, we repeat the same arguments in 
$C([0,T];L^2(\OOO))$ as the ones 
employed in the proof of 
Theorem \ref{t2.1} to estimate 
$A_0^{\frac{3}{2}}u(t)$ in the norm $C([0,T];L^2(\OOO))$ by the fixed point 
theorem and reach the inclusion
 $A_0^{\frac{3}{2}}u \in C([0,T];L^2(\OOO))$.

Collecting the last estimates, we obtain
\begin{equation}
\label{(6.7)}
\left\{ \begin{array}{rl}
& u \in C([0,T];\mathcal{D}(A_0^{\frac{3}{2}})) \subset C([0,T];H^3(\OOO)), \\
& \Vert u(t)\Vert_{H^3(\OOO)} \le C\left( 
1 + t^{\alpha\left(\frac{1}{4}-\ep_0\right)} \right)D,
\quad 0<t<T.                          
\end{array}\right.
\end{equation}

Next we estimate $\Vert Au'(t)\Vert$.  First we represent $u'(t)$ in the form
\begin{align*}
& u'(t) = G'(t) + \frac{d}{dt}\int^t_0 K(t-s)Q(s)u(s) ds\\
= & G'(t) + \frac{d}{dt}\int^t_0 K(s)Q(t-s)u(t-s) ds\\
= & G'(t) + K(t)Q(0)u(0) \\
+ & \int^t_0 K(s) (Q(t-s)u'(t-s)
+ Q'(t-s)u(t-s)) ds, \quad 0<t<T,
\end{align*}
so that  
\begin{equation}
\label{(6.8)}
A_0u'(t) = A_0G'(t) + A_0K(t)Q(0)u(0) 
\end{equation}
$$
+ \int^t_0 A_0^{\hhalf}K(s) A_0^{\hhalf}(Q(t-s)u'(t-s)
+ Q'(t-s)u(t-s)) ds, \quad 0<t<T.       
$$
Similarly to the arguments applied for derivation of \eqref{(6.4)}, we obtain 
the inequality
$$
\Vert A_0^{\hhalf}(Q(t-s)u'(t-s) + Q'(t-s)u(t-s))\Vert 
\le C\Vert A_0u'(t-s)\Vert, \quad 0<t<T.
$$
The inclusion $Q(0)u(0) = Q(0)a \in C^2_0(\OOO) \subset \DDD(A_0)$ follows  
from the regularity conditions \eqref{(1.2)}
and the inclusion $a \in C^{\infty}_0(\OOO)$.  
Furthermore, by \eqref{(5.2)} and \eqref{(5.3)}, we see
$$
\Vert A_0S'(t)a\Vert 
= \Vert A_0^2K(t)a\Vert = \Vert K(t)A_0^2a\Vert 
\le Ct^{\alpha-1}\Vert A_0^2a\Vert \le Ct^{\alpha-1}\Vert a\Vert_{H^4(\OOO)}
$$
and 
$$
\Vert K(t)A_0(Q(0)a)\Vert \le Ct^{1-\alpha}\Vert A_0(Q(0)a)\Vert 
\le Ct^{\alpha-1}\Vert a\Vert_{H^3(\OOO)}.
$$
Hence, the representation \eqref{(6.8)} leads to the estimate
$$
\Vert A_0u'(t)\Vert \le Ct^{\alpha-1}D  
+ C\int^t_0 s^{\hhalf\alpha -1}\Vert A_0u'(t-s)\Vert ds,
\quad 0<t<T.
$$
We define a normed space
$$
\www{X}:= \{v\in C([0,T];L^2(\OOO)) \cap C^1((0,T];L^2(\OOO));\, 
t^{1-\alpha}\ppp_tv \in C([0,T];L^2(\OOO))\}
$$
and
$$
\Vert v\Vert_{\www{X}}:= \max_{0\le t\le T} \Vert t^{1-\alpha}\ppp_tv
(\cdot,t)\Vert_{L^2(\OOO)}
+ \max_{0\le t\le T} \Vert v(\cdot,t)\Vert_{L^2(\OOO)}.
$$
We readily verify that $\www{X}$ is a Banach space.

Arguing similarly to the proof of Theorem \ref{t2.1} and applying the fixed 
point theorem in $\www{X}$, we conclude that 
$A_0u \in \www{X}$, that is, $t^{1-\alpha}A_0u' \in C([0,T];L^2(\OOO))$. 
Using the inclusion  $\DDD(A_0) \subset 
C(\ooo{\OOO})$ by the spatial dimensions $d=1,2,3$, the Sobolev embedding 
theorem yields 
\begin{equation}
\label{(6.9)}
u' \in C(\ooo{\OOO} \times (0,T]), \quad \Vert A_0u'(t)\Vert 
\le CDt^{\alpha-1}, \quad 0\le t\le T.
\end{equation}

Now we proceed to the estimation of $A_0^2u(t)$.
Since $\frac{d}{ds}(-A_0^{-1}S(s)) = K(s)$ for $0<s<T$ by \eqref{(5.2)}, 
the integration by parts yields 
\begin{align*}
& \int^t_0 K(t-s)Q(s)u(s) ds = \int^t_0 K(s)Q(t-s)u(t-s) ds \\
= & \left[ -A_0^{-1}S(s)Q(t-s)u(t-s)\right]^{s=t}_{s=0}\\
- &
\int^t_0 A_0^{-1}S(s)(Q'(t-s)u(t-s)+Q(t-s)u'(t-s)) ds\\
= & A_0^{-1}Q(t)u(t) - A_0^{-1}S(t)Q(0)u(0)
\end{align*}
\begin{equation}
\label{(6.10)}
- \int^t_0 A_0^{-1}S(s)(Q'(t-s)u(t-s)+Q(t-s)u'(t-s)) ds,
\quad 0<t<T.                        
\end{equation}
The Lebesgue convergence theorem and the estimate
$\vert \MLONE(\eta)\vert \le \frac{C}{1+\eta}$ for all 
$\eta>0$ (Theorem 1.6 in \cite{Po}), we deduce 
$$
\Vert (S(t) - 1)a\Vert^2
= \sumn \vert (a,\va_n)\vert^2 (\MLONE(-\la_nt^{\alpha}) - 1)^2
\, \longrightarrow\, 0
$$
as $t \to \infty$ for $a \in L^2(\OOO)$.  
\\
Hence,
$u \in C([0,T];L^2(\OOO))$ and $\lim_{t\downarrow 0}
\Vert (S(t)-1)a\Vert = 0$ and so 
$$
\lim_{s\downarrow 0} S(s)Q(t-s)u(t-s) = S(0)Q(t)u(t) \quad
\mbox{in } L^2(\OOO)
$$
and
$$
\lim_{s\uparrow t} S(s)Q(t-s)u(t-s) = S(t)Q(0)u(0) \quad
\mbox{in }L^2(\OOO),
$$
which justify the last equality in the formula \eqref{(6.10)}.

Thus, in terms of \eqref{(6.10)}, the representation \eqref{(5.7)} 
can be rewritten in the form
$$
A_0^2(u(t) - A_0^{-1}Q(t)u(t)) = A_0^2G(t) -A_0S(t)Q(0)u(0) 
$$
\begin{equation}
\label{(6.11)}
- \int^t_0 A_0^{\hhalf}S(s)A_0^{\hhalf}
(Q'(t-s)u(t-s) + Q(t-s)u'(t-s)) ds, \quad
0<t<T.                                       
\end{equation}
Since $u(0) = a \in C^{\infty}_0(\OOO)$ and $F \in C^{\infty}_0(\OOO
\times (0,T))$, in view of \eqref{(1.2)} we have the inclusion
$$
A_0^2G(\cdot) \in C([0,T];L^2(\OOO)), \quad 
A_0S(t)Q(0)u(0) = S(t)(A_0Q(0)a) \in C([0,T];L^2(\OOO)).
$$
Now we use the conditions \eqref{(1.2)} and \eqref{(5.3)} and repeat the 
arguments employed for derivation 
of \eqref{(6.4)} by means of \eqref{(6.5)} and \eqref{(6.9)} to obtain the 
estimate
\begin{align*}
& \left\Vert \int^t_0 A_0^{\hhalf}S(s)A_0^{\hhalf} 
(Q'(t-s)u(t-s) + Q(t-s)u'(t-s)) ds \right\Vert\\
\le &C\int^t_0 s^{-\hhalf\alpha}\Vert Q'(t-s)u(t-s) + Q(t-s)u'(t-s)
\Vert_{H^1(\OOO)} ds\\
\le& C\int^t_0 s^{-\hhalf\alpha}(\Vert A_0u'(t-s)\Vert
+ \Vert A_0u(t-s)\Vert) ds
\le Ct^{\hhalf\alpha}D
\end{align*}
and the inclusion
$$
-\int^t_0 A_0^{\hhalf}S(s)A_0^{\hhalf}(Q'(t-s)u(t-s) + Q(t-s)u'(t-s)) ds 
\in C([0,T];L^2(\OOO)).
$$
Therefore, 
$$
A_0^2(u(t) - A_0^{-1}Q(t)u(t)) = A_0(A_0u(t) - Q(t)u(t))
\in C([0,T];L^2(\OOO)),
$$
that is,
$$
A_0u(t) - Q(t)u(t) \in C([0,T]; \DDD(A_0)) \subset C([0,T];H^2(\OOO)).
$$
On the other hand, the estimate \eqref{(6.7)} implies  
$$
Q(t)u(t) \in C([0,T];H^2(\OOO))
$$
and we obtain 
\begin{equation}
\label{(6.12a)}
A_0u(t) \in C([0,T];H^2(\OOO)).   
\end{equation}

For further arguments, we define the Schauder spaces $C^{\theta}(\ooo{\OOO})$ 
and 
$C^{2+\theta}(\ooo{\OOO})$ with $0<\theta<1$ (see e.g., \cite{GT},
\cite{LU}) as follows:  A function $w$ is defined to belong to the space 
$C^{\theta}(\ooo{\OOO})$ if 
$$
\sup_{x \ne x', x, x'\in \OOO} 
\frac{\vert w(x) - w(x')\vert}{\vert x-x'\vert^{\theta}}
< \infty.
$$  
For  $w \in C^{\theta}(\ooo{\OOO})$, we define the norm
$$
\Vert w\Vert_{C^{\theta}(\ooo{\OOO})}
:= \Vert w\Vert_{C(\ooo{\OOO})}
+ \sup_{x\ne x', x, x'\in \OOO} 
\frac{\vert w(x) - w(x')\vert}{\vert x-x'\vert^{\theta}}.
$$
and for $w\in C^{2+\theta}(\ooo{\OOO})$, the norm is given by
$$
\Vert w\Vert_{C^{2+\theta}(\ooo{\OOO})}
:= \Vert w\Vert_{C^2(\ooo{\OOO})}
+ \sum_{\vert \tau\vert=2}
\sup_{x\ne x', x, x'\in \OOO} 
\frac{\vert \ppp_x^{\tau}w(x) - \ppp_x^{\tau}w(x')\vert}
{\vert x-x'\vert^{\theta}} < \infty.
$$

In the last formula, the notations 
 $\tau := (\tau_1, ..., \tau_d) \in (\N \cup \{0\})^d$,
$\ppp_x^{\tau}:= \ppp_1^{\tau_1}\cdots \ppp_d^{\tau_d}$, and
$\vert \tau\vert:= \tau_1 + \cdots + \tau_d$ are employed.

For $d=1,2,3$, the Sobolev embedding theorem says that $H^2(\OOO) \subset 
C^{\theta}(\ooo{\OOO})$ with some $\theta \in (0,1)$
(\cite{Ad}).  

Therefore, in view of \eqref{(6.12a)}, we see $h:= A_0u(\cdot,t) 
\in C^{\theta}(\ooo{\OOO})$ for each $t \in [0,T]$.
We apply the Schauder 
estimate (see e.g., \cite{GT} or \cite{LU})
for solutions to the elliptic boundary value problem
$$
A_0u(\cdot,t) = h\in C^{\theta}(\ooo{\OOO}) \quad \mbox{in } \OOO
$$
with the boundary condition $\NUNU u(\cdot,t) + \sigma(\cdot)u(\cdot,t) = 0$ 
on $\ppp\OOO$ to arrive at the inclusion
\begin{equation}
\label{(6.13)}
u \in C([0,T]; C^{2+\theta}(\ooo{\OOO})).     
\end{equation}
The inclusions \eqref{(6.9)} and \eqref{(6.13)} complete the proof of 
Lemma \ref{l6.1}.
\end{proof}

{\bf (II) Second part of the proof of Theorem \ref{t2.2}}.

In this part, we weaken the regularity conditions posed on $u$ in 
Lemma \ref{l4.2} and prove the same
results provided that $u\in L^2(0,T;H^2(\OOO))$ and 
$u-a \in \HH(0,T;L^2(\OOO))$.

Let $F \in L^2(0,T;L^2(\OOO))$ and $a\in H^1_0(\OOO)$ be arbitrarily 
chosen such that $F\ge 0$ in $\OOO\times (0,T)$ and $a\ge 0$ in $\OOO$.

Now we apply the standard mollification procedure (see, e.g., \cite{Ad}) 
and construct the sequences
$F_n \in C^{\infty}_0(\OOO\times (0,T))$ and $a_n \in C^{\infty}_0(\OOO)$,
$n\in \N$ such that $F_n\ge 0$ in $\OOO\times (0,T)$ and $a_n\ge 0$ in 
$\OOO$, $n\in \N$ and $\lim_{n\to\infty}
\Vert F_n-F\Vert_{L^2(0,T;L^2(\OOO))} = 0$ and 
$\lim_{n\to\infty}\Vert a_n-a\Vert_{H^1_0(\OOO)} = 0$.
Then Lemma \ref{l6.1} yields the inclusion
$$
u(F_n,a_n) \in C([0,T];C^2(\ooo{\OOO})), \quad
t^{1-\alpha}\ppp_tu(F_n,a_n) \in C([0,T];C(\ooo{\OOO})),
\quad n\in \N
$$
and thus Lemma \ref{l4.2} ensures the inequalities
\begin{equation}
\label{(6.14)}
u(F_n,a_n) \ge 0 \quad \mbox{in } \OOO\times (0,T), \, n\in \N.
\end{equation}
Since Theorem \ref{t2.1} hold true  for the the 
initial-boundary value problem \eqref{(4.2)} with $F$ and $a$ 
replaced by $F-F_n$ and $a-a_n$, respectively, we have 
$$
\Vert u(F,a) - u(F_n,a_n)\Vert_{L^2(0,T; H^2(\OOO))}
$$
$$ 
\le C(\Vert a-a_n\Vert_{H^1(\OOO)} 
+ \Vert F-F_n\Vert_{L^2(0,T;L^2(\OOO))}) \, \longrightarrow \, 0
$$
as $n\to \infty$.  Therefore, we can choose a subsequence $m(n)\in \N$
such that $u(F,a)(x,t) = \lim_{m(n)\to \infty} u(F_{m(n)},a_{m(n)})(x,t)$
for almost all $(x,t) \in \OOO\times (0,T)$.
Then the inequality \eqref{(6.14)} leads to the desired result that  
$u(F,a)(x,t) \ge 0$ for almost all 
$(x,t) \in \OOO\times (0,T)$.

{\bf (III) Third part: Completion of the proof of Theorem \ref{t2.2}.}

Let $a\ge 0$, $\in H^1_0(\OOO)$, $F \ge 0$,
$\in L^2(0,T;L^2(\OOO))$, and let $u=u(F,a) \in 
L^2(0,T;H^2(\OOO))$ satisfy the initial condition $u-a \in \HH(0,T;L^2(\OOO))$ 
and the boundary value problem.
Now we will prove Theorem \ref{t2.2} for $u$ in this class without assumptions
on the sign of the zeroth-order coefficient.

In \eqref{(3.1a)}, we choose a constant function $b_0>0$ as $b_0(x,t)$, and
assume that $b_0>0$ is sufficiently large.
We define the operator $A_1$ by \eqref{(3.1a)}.
Then we verify that the initial-boundary value problem \eqref{(2.3)}
is equivalent to 
\begin{equation}
\label{(6.15)}
\left\{ \begin{array}{rl}
& \pppa (u-a) + A_1u = (b_0+c(x,t))u + F, \quad 
(x,t) \in \OOO\times (0,T), \\
& \NUNU u + \sigma u = 0 \quad \mbox{on $\ppp\OOO\times (0,T)$}.
\end{array}\right.
\end{equation}
In what follows, we choose sufficiently large $b_0>0$ such that 
$b_0 \ge \Vert c\Vert_{C(\ooo{\OOO} \times [0,T])}$.

In the previous parts of the proof, we already interpreted the solution $u$ as 
a unique fixed point for the equation \eqref{(6.1)}. 
Now let us construct an approximating sequence 
$u_n$, $n\in \N$ for $u$ as follows. Setting 
$$
u_1(x,t) = a(x) \ge 0, \quad (x,t) \in \OOO\times (0,T),
$$
we inductively define a sequence $\{u_n\}_{n\in \N}$ 
by solving $u_{n+1}$ with given $u_n$:
\begin{equation}
\label{(6.16)}
\left\{ \begin{array}{rl}
&\pppa (u_{n+1}-a) + A_1u_{n+1} = (b_0+c(x,t))u_n + F
\quad \mbox{in } \OOO\times (0,T),\\
& \NUNU u_{n+1} + \sigma u_{n+1} = 0 \quad \mbox{on } \ppp\OOO\times (0,T),\\
& u_{n+1} - a \in \HH(0,T;L^2(\OOO)), \quad n\in \N.
\end{array}\right.
\end{equation}
We set $u_0(x,t) := 0$ for $(x,t) \in \OOO\times (0,T)$.
First we show that
\begin{equation}
\label{(6.17)}
u_n(x,t) \ge 0, \quad (x,t) \in \OOO\times (0,T), \quad n\in \N.
\end{equation}
Indeed, the inequality \eqref{(6.17)} holds for $n=1$.  Now we assume that 
$u_n \ge 0$ in $\OOO\times (0,T)$.
Then $(b_0+c(x,t))u_n + F \ge 0$ in $\QQQQ$, and thus by the results 
established in the second part of the proof of Theorem \ref{t2.2}, 
we obtain $u_{n+1} \ge 0$ in $\QQQQ$.  
Thus, by the principle of mathematical induction, 
the inequality \eqref{(6.17)} holds true for all $n\in \N$.

We rewrite \eqref{(6.16)} as
$$
\pppa (u_{n+1}(t) - a) + A_0u_{n+1}(t) 
= (Q(t)u_{n+1}-(c(t)+b_0)u_{n+1})
+ (b_0+c(t))u_n + F,
$$
where we recall that $A_0$ and $Q(t)$ are defined by \eqref{(3.2)} and \eqref{(5.8)}.
Next we estimate $w_{n+1}:= u_{n+1} - u_n$.  By the relation \eqref{(6.16)}, 
we obtain
$$
\left\{ \begin{array}{rl}
&\pppa w_{n+1} + A_0w_{n+1} = (Q(t)w_{n+1} - (c(t)+b_0)w_{n+1})
+ (b_0+c(x,t))w_n \\ 
& \qquad \qquad \quad \mbox{in } \OOO\times (0,T),\\
& \NUNU w_{n+1} + \sigma w_{n+1} = 0 \quad \mbox{on } \ppp\OOO\times (0,T),\\
& w_{n+1} \in \HH(0,T;L^2(\OOO)), \quad n\in \N.
\end{array}\right.
$$
In terms of the operator $K(t)$ defined by \eqref{(5.2)}, acting  
similarly to our analysis of the fixed point equation \eqref{(6.1)}, 
we have 
\begin{align*}
& w_{n+1}(t) = \int^t_0 K(t-s)(Qw_{n+1})(s) ds
- \int^t_0 K(t-s)(c(s)+b_0)w_{n+1}(s) ds\\
+ & \int^t_0 K(t-s)(b_0+c(s))w_n(s) ds, \quad 0<t<T, 
\end{align*}
which leads to the estimation
\begin{align*}
& \Vert A_0^{\hhalf}w_{n+1}(t)\Vert 
\le \int^t_0 \Vert A_0^{\hhalf}K(t-s)\Vert 
\Vert Q(s)w_{n+1}(s)\Vert ds\\
+& \int^t_0 \Vert A_0^{\frac{1}{2}}K(t-s)\Vert
\Vert (c(s)+b_0)w_{n+1}(s)\Vert ds
+ \int^t_0 \Vert A_0^{\hhalf}K(t-s)\Vert \Vert (b_0+c(s))w_n(s)\Vert ds\\
\le& C\int^t_0 (t-s)^{\hhalf\alpha-1}
\Vert A_0^{\hhalf}w_{n+1}(s)\Vert ds 
+ C\int^t_0 (t-s)^{\hhalf\alpha-1}\Vert A_0^{\hhalf}w_n(s)\Vert ds
\quad \mbox{for $0<t<T$.}
\end{align*}
Here, by \eqref{(1.2)} we used the norm estimates
$$
\Vert Q(s)w_{n+1}(s)\Vert 
\le C\Vert w_{n+1}(s)\Vert_{H^1(\OOO)} 
\le C\Vert A_0^{\hhalf}w_{n+1}(s)\Vert
$$
and
$$
\Vert (c(s)+b_0)w_\ell(s)\Vert \le C\Vert w_\ell(s)\Vert_{H^1(\OOO)}
\le C\Vert A_0^{\hhalf}w_\ell(s)\Vert, \quad \ell=n, n+1.
$$
Thus  
\begin{align*}
& \Vert A_0^{\hhalf}w_{n+1}(t)\Vert 
\le C\int^t_0 (t-s)^{\frac{1}{2}\alpha-1} 
\Vert A_0^{\hhalf}w_{n+1}(s)\Vert ds\\
+ & C\int^t_0 (t-s)^{\frac{1}{2}\alpha-1} 
\Vert A_0^{\hhalf}w_n(s)\Vert ds, \quad 0<t<T.
\end{align*}
The generalized Gronwall inequality yields the inequality
\begin{align*}
&\Vert \Ahalf w_{n+1}(t)\Vert
\le C\int^t_0 (t-s)^{\hhalf\alpha-1}\Vert \Ahalf w_n(s)\Vert ds\\
+& C\int^t_0 (t-s)^{\hhalf\alpha-1}
\left( \int^s_0 (s-\xi)^{\hhalf\alpha-1}
\Vert \Ahalf w_n(\xi)\Vert d\xi\right)ds.
\end{align*}
For the second term of the right-hand side of the last inequality, 
the following calculations hold true:
\begin{align*}
& \int^t_0 (t-s)^{\hhalf\alpha-1}
\left( \int^s_0 (s-\xi)^{\hhalf\alpha-1}\Vert \Ahalf w_n(\xi)\Vert d\xi\right)
ds\\
=& \int^t_0 \Vert \Ahalf w_n(\xi)\Vert \left( \int^t_{\xi}
(t-s)^{\hhalf\alpha-1} (s-\xi)^{\hhalf\alpha-1} ds \right) d\xi\\
=& \frac{\Gamma\left( \hhalf\alpha\right)\Gamma\left( \hhalf\alpha\right)}
{\Gamma(\alpha)}\int^t_0 (t-\xi)^{\alpha-1}\Vert \Ahalf w_n(\xi)\Vert d\xi\\
=& \frac{\Gamma\left( \hhalf\alpha\right)^2}{\Gamma(\alpha)}T^{\hhalf\alpha}
\int^t_0 (t-s)^{\hhalf\alpha-1}\Vert \Ahalf w_n(s)\Vert ds.
\end{align*}
Thus, we can choose a constant $C>0$ dependent on $\alpha$ and $T$, 
such that 
\begin{equation}
\label{(ineq1)}
\Vert \Ahalf w_{n+1}(t)\Vert 
\le C\int^t_0 (t-\xi)^{\hhalf\alpha-1}\Vert \Ahalf w_n(s)\Vert ds,
\quad 0<t<T,\, n\in \N.
\end{equation}
Recalling that 
$$
\int^t_0 (t-s)^{\hhalf\alpha -1}\eta(s)ds = \Gamma\left( \hhalf\alpha\right)
(J^{\hhalf\alpha}\eta)(t), \quad t>0,
$$
and setting $\eta_n(t):= \Vert \Ahalf w_n(t)\Vert$, we can rewrite 
\eqref{(ineq1)} as follows:
\begin{equation}
\label{(6.18)}
\eta_{n+1}(t) \le C\Gamma\left( \hhalf\alpha\right)(J^{\hhalf\alpha}\eta_n)(t),
\quad 0<t<T, \, n\in \N.                         
\end{equation}
Since $J^{\hhalf\alpha}$ preserves the signs of $\eta$ and 
$J^{\beta_1}(J^{\beta_2}\eta)(t) = J^{\beta_1+\beta_2}\eta(t)$ for 
each $\beta_1, \beta_2 > 0$, applying the inequality \eqref{(6.18)} 
repeatedly, we reach the inequality
\begin{align*}
& \eta_n(t) \le \left(C\Gamma\left( \hhalf\alpha\right)\right)^{n-1}
(J^{(n-1)\frac{\alpha}{2}}\eta_1)(t) \\
= & \frac{\left(C\Gamma\left( \hhalf\alpha\right)\right)^{n-1}}
{\Gamma\left( \frac{\alpha}{2}(n-1)\right)}
\left( \int^t_0 (t-s)^{(n-1)\hhalf\alpha-1} ds\right) \Vert \Ahalf a\Vert\\
= & \frac{\left(C\Gamma\left( \hhalf\alpha\right)\right)^{n-1}}
{\Gamma\left( \frac{\alpha}{2}(n-1)\right)}
\frac{t^{(n-1)\hhalf\alpha}}{(n-1)\hhalf\alpha} \Vert \Ahalf a\Vert
\le C_1\frac{\left(C\Gamma\left( \hhalf\alpha\right) T^{\frac{\alpha}{2}}
\right)^{n-1}}
{\Gamma\left( \frac{\alpha}{2}(n-1)\right)}.
\end{align*}
The known asymptotic behavior of the gamma function justifies the relation
$$
\lim_{n\to \infty} 
\frac{\left(C\Gamma\left( \hhalf\alpha\right) T^{\frac{\alpha}{2}}
\right)^{n-1}}
{\Gamma\left( \frac{\alpha}{2}(n-1)\right)} = 0.
$$
Thus we have proved that the sequence $u_N = w_0 + \cdots + w_N$ converges to 
the solution 
$u$ in $L^{\infty}(0,T;H^1(\OOO))$ as $N \to \infty$.
Therefore, we can choose a subsequence $m(n)\in \N$ for $n \in \N$
such that $\lim_{m(n)\to\infty} u_{m(n)}(x,t) = u(x,t)$ for 
almost all $(x,t) \in \QQQQ$.
This statement in combination with the inequality \eqref{(6.17)} means that 
$u(x,t) \ge 0$ for almost all
$(x,t) \in \QQQQ$.  Thus the proof of Theorem \ref{t2.2} is completed.
\end{proof}

\vspace{0.2cm}

Now let us fix a source function $F = F(x,t) \ge 0$ and an initial 
value $a \in H^1_0(\OOO)$ in the initial-boundary value problem \eqref{(2.3)}  
and denote 
% by $u_c = u_c(x,t)$,
by $u(c,\sigma) = u(c,\sigma)(x,t)$  the solution 
to the  problem \eqref{(2.3)} with
the functions $c=c(x,t)$ and $\sigma = \sigma(x)$.
Then the following comparison property regarding the coefficients 
$c$ and $\sigma$ is valid:

\begin{theorem}
\label{t2.3}
Let $a\in H^1_0(\OOO)$ and $a\ge 0$ in $\OOO$ and $F \in L^2(\OOO\times (0,T))$
satisfy $F\ge 0$ in $\OOO \times (0,T)$.
\\
(i) Let $c_1, c_2 \in C^1([0,T]; C^1(\ooo{\OOO})) \cap
C([0,T];C^2(\ooo{\OOO}))$ satisfy 
$c_1(x,t) \ge c_2(x,t)$ in $\OOO$.
Then $u(c_1,\sigma)(x,t) \ge u(c_2,\sigma)(x,t)$ in $\OOO \times (0,T)$.
\\
(ii) Let $c(x,t) < 0$ for $(x,t) \in \QQQQ$.
We arbitrarily fix a constant $\sigma_0>0$.
If the smooth functions $\sigma_1, \sigma_2$ on $\ppp\OOO$ satisfy
$$
\sigma_2(x) \ge \sigma_1(x) \ge \sigma_0 \quad \mbox{for $x\in \ppp\OOO$},
$$
then $u(c,\sigma_1) \ge u(c, \sigma_2) \quad \mbox{in } \QQQQ$.
\end{theorem}

\begin{proof}
(i) Because $a\ge 0$ in $\OOO$ and $F\ge 0$ in $\QQQQ$, Theorem \ref{t2.2} 
yields the inequality
$u(c_2,\sigma)\ge 0$ in $\QQQQ$.
Setting $u(x,t):= u(c_1,\sigma)(x,t) - u(c_2,\sigma)(x,t)$ for $(x,t) 
\in \QQQQ$, we have
$$
\left\{ \begin{array}{rl}
& \pppa u - \sumij \ppp_i(a_{ij}\ppp_ju) 
- \sum_{j=1}^d b_j \ppp_ju\\
-& c_1(x,t)u = (c_1-c_2)u(c_2,\sigma)(x,t) \quad \mbox{in }\QQQQ, \\
& \NUNU u + \sigma u = 0 \quad \mbox{on } \ppp\OOO,\\
& u \in \HH(0,T;L^2(\OOO)).
\end{array}\right.
$$
Since $u(c_2,\sigma)(x,t) \ge 0$ and $(c_1-c_2)(x,t) \ge 0$ for $(x,t) 
\in \QQQQ$,
 Theorem \ref{t2.2} leads to the estimate $u(x,t)\ge 0$ for $(x,t) \in \QQQQ$, 
which is equivalent to the inequality $u(c_1,\sigma)(x,t) 
\ge u(c_2,\sigma)(x,t)$ for $(x,t) \in \QQQQ$.
\\
(ii) Similarly to the procedure applied for the second part of the proof of 
Theorem \ref{t2.2}, we choose the 
sequences $F_n \ge 0$, $\in C^{\infty}_0(\QQQQ)$ and 
$a_n \ge 0$, $\in C^{\infty}_0(\OOO)$, $n\in \N$ such that 
$F_n \longrightarrow F$ in $L^2(\QQQQ)$ and
$a_n \longrightarrow a$ in $H^1_0(\OOO)$.  Let $u_n$, $v_n$ be the 
solutions to initial-boundary value problem \eqref{(2.3)} with 
$F=F_n$, $a=a_n$ and with the coefficients 
$\sigma_1$ and $\sigma_2$ in the boundary condition,
respectively.  We note that $v_n, u_n \in C(\ooo{\OOO} \times [0,T])$
and $t^{1-\alpha}\ppp_tv_n, \, t^{1-\alpha}\ppp_tu_n 
\in C([0,T]:C(\ooo{\OOO}))$, $n\in \N$ by Lemma \ref{l6.1}, 
and Theorem \ref{t2.2} yields that 
\begin{equation}\label{(4.22a)}
v_n(x,t) \ge 0, \quad (x,t) \in \ppp\OOO\times (0,T).
\end{equation} 

Then Theorem \ref{t2.1} yields the relation
\begin{equation}
\label{(6.19)}
\lim_{n\to\infty}\Vert u_n - u(c,\sigma_1)\Vert_{L^2(0,T;L^2(\OOO))}
= \lim_{n\to\infty}\Vert v_n - u(c,\sigma_2)\Vert_{L^2(0,T;L^2(\OOO))}
= 0.                                         
\end{equation}
Let us now define an auxiliary function $w_n:= u_n - v_n$.
Then, 
\begin{equation}\label{(4.23)}
t^{1-\alpha}\ppp_tw_n \in C([0,T]:C(\ooo{\OOO})), \quad 
w_n \in C([0,T];C^2(\ooo{\OOO})), \quad n\in \N
\end{equation}
and 
\begin{equation}\label{(4.24)}
\left\{ \begin{array}{rl}
& \pppa w_n + Aw_n = 0 \quad \mbox{in } \QQQQ,\\
& \NUNU w_n + \sigma_1w_n = (\sigma_2-\sigma_1)v_n
  \quad \mbox{on } \ppp\OOO\times (0,T),\\
& w_n(x,\cdot) \in \HH(0,T) \quad \mbox{for almost all } x\in \OOO.
\end{array}\right.
\end{equation}
Hence, by \eqref{(4.22a)} and $\sigma_2 \ge \sigma_1$ on $\ppp\OOO$, 
we have
\begin{equation}\label{(4.23a)}
\NUNU w_n + \sigma_1w_n \ge 0 \quad \mbox{on $\ppp\OOO\times (0,T)$}.
\end{equation}

We show a variant of Lemma \ref{l4.2}.
\begin{lemma}
\label{l4.2a}
Let the elliptic operator $-A$ be defined by \eqref{(2.1)}.
We assume \eqref{(1.2)} and $c(x,t) < 0$ for $x \in \ooo{\OOO}$ and
$0\le t \le T$ and there exists a constant $\sigma_0>0$ such that 
$$
\sigma(x) \ge \sigma_0 \quad \mbox{for all $x\in \ppp\OOO$}.
$$
For $a \in H^1_0(\OOO)$ and $F\in L^2(\OOO\times (0,T))$, we further assume
that there exists a solution $u\in C([0,T];C^2(\ooo{\OOO}))$ 
satisfying $t^{1-\alpha}\ppp_tu \in C([0,T];C(\ooo{\OOO}))$ to
%\begin{equation}\label{(4.24a)}
$$
\left\{ \begin{array}{rl}
& \pppa (u-a) + Au = F \quad \mbox{in $\OOO\times (0,T)$}, \\
& \ppp_{\nu_A}u + \sigma(x)u \ge 0, \quad \mbox{on 
$\ppp\OOO \times (0,T)$}, \\
& u(x,\cdot) - a\in H_{\alpha}(0,T) \quad \mbox{for almost all
$x\in \OOO$}.
\end{array}\right.
$$
%\end{equation}
If $F \ge 0$ in $\OOO \times (0,T)$ and $a\ge 0$ in $\OOO$, then 
$u \ge 0$ in $\OOO\times (0,T)$.
\end{lemma}

In this lemma, at the expense of extra condition $\sigma > 0$
on $\ppp\OOO$, we do not need assume that $\min_{(x,t)\in 
\ooo{\OOO}\times [0,T]} (-c(x,t))$ is sufficiently large, which is the main 
difference from Lemma \ref{l4.2}.
The proof is much simpler than Lemma \ref{l4.2} and
postponed to the end of this section.

Now we complete the proof of Theorem \ref{t2.3}.
Since $c(x,t) < 0$ for $(x,t)\in \QQQQ$ and $\sigma_1 \ge \sigma_0 > 0$
on $\ppp\OOO$, in terms of \eqref{(4.23)} and \eqref{(4.23a)} we can apply 
Lemma \ref{l4.2a} to \eqref{(4.24)} and deduce the inequality $w_n \ge 0$ in $\QQQQ$, 
that is,
$u_n \ge v_n$ in $\QQQQ$ for $n\in \N$.  In view of \eqref{(6.19)}, 
we can choose
a suitable subsequence of $w_n,\ n\in \N$ and pass to the limit, 
thus arriving at the inequality  
$u(c,\sigma_1) \ge u(c,\sigma_2)$ in $\QQQQ$.
The proof of Theorem \ref{t2.3} is completed.
\end{proof}

Finally, let us mention one direction for further research in connection with 
the results formulated and proved in this sections. In order to remove the 
negativity condition posed on the coefficient $c=c(x,t)$ in Theorem \ref{t2.3} 
(ii),
one needs a unique existence result for  solutions to the initial-boundary 
value problems  of type \eqref{(2.3)}
with non-zero Robin boundary condition similar to the one formulated in Theorem 
\ref{t2.1}.
There are several works that treat the case of the initial-boundary value 
problems  with non-zero Dirichlet boundary conditions
(see, e.g., \cite{Ya18} and the references therein).
However, to the best of the authors' knowledge, analogous  results are not 
available for the case of the initial-boundary value problems with 
the Neumann and Robin boundary conditions.  Thus, in Theorem \ref{t2.3} (ii), 
we assumed the condition 
$c<0$ in $\QQQQ$, although our conjecture is that this result is valid for 
an arbitrary coefficient $c=c(x,t)$.

We conclude this section with 
\\
{\bf Proof of Lemma \ref{l4.2a}.}
The proof is simple because we do not need the function 
$\psi$ defined by \eqref{(4.3)}.
We set
$$
\www{w}(x,t):= u(x,t) + \ep(1+t^{\alpha}), \quad x\in \OOO,\,
0<t<T.
$$
Using $c<0$ on $\ooo{\OOO} \times [0,T]$ and $\sigma \ge \sigma_0>0$ on 
$\ppp\OOO$ and calculating similarly to the proof of Lemma \ref{l4.2}, 
we have
$$
\ddda \www{w} + A\www{w} = F + \ep\Gamma(\alpha+1)
- c(x,t)\ep(1+t^{\alpha}) > 0 \quad \mbox{in $\OOO\times (0,T)$},
$$
$$
\ppp_{\nu_A}\www{w} + \sigma \www{w} = \ppp_{\nu_A}u + \sigma u 
+ \sigma\ep(1+t^{\alpha}) \ge \sigma_0\ep
\quad \mbox{on $\ppp\OOO \times (0,T)$}
$$
and
$$
\www{w}(x,0) = a(x) + \ep \ge \ep \quad \mbox{in $\OOO$}.
$$
Therefore, we can follow the same arguments after \eqref{(4.7)} in the proof of
Lemma \ref{l4.2}.  Thus the proof of Lemma \ref{l4.2a} is complete.
$\square$
%
%
%%%%%%%%%%%%%%%% Section 5 %%%%%%%%%%%%

\section{Well-posedness for an initial-boundary value problems for
semilinear time-fractional diffusion equations} 
\label{sec5}

\setcounter{section}{5}

In order to extend the comparison principle to semilinear 
time-fractional diffusion equations, we have to establish the unique
existence and a priori estimates for initial-boundary 
value problems, because the arguments rely on the construction of 
the solutions, as we see from the third part of the proof of 
Theorem \ref{t2.2} in Section \ref{sec4}.

Thus, in this section, we first formulate and prove some existence results 
for the following initial-boundary value problem for the semilinear 
time-fractional diffusion equation:
\begin{equation}
\label{(3.1)}
\left\{ \begin{array}{rl}
& \pppa (u(x,t)-a(x)) + Au(x,t) = f(u)(x,t), \quad x \in \Omega, 
\thinspace 0<t<T, \\
& \NUNU u + \sigma(x)u(x,t) = 0, \quad x\in \ppp\OOO, \, 0<t<T,\\
& u(x,\cdot) - a(x) \in \HH(0,T) \quad \mbox{for almost all } 
x\in \OOO,
\end{array}\right.
\end{equation}
where the second order spatial differential operator $-A$ is defined by 
\eqref{(2.1)} and the source function $f(u)(\cdot,t)$ may depend not only on 
$u$ but also on $x$ and its spatial derivatives, 
provided that $f$ satisfies the conditions \eqref{(3.3)} stated below.

We recall that for a fixed but sufficiently large constant $c_0>0$,
the operator $A_0$ in 
$L^2(\OOO)$ is defined by the relation
$$
\left\{ \begin{array}{rl}
& (-A_0v)(x) = \sumij \ppp_i(a_{ij}(x)\ppp_jv(x)) - c_0v(x), \quad
x\in \OOO, \\
& \mathcal{D}(A_0) = \left\{ v\in H^2(\OOO);\,
\NUNU v + \sigma v = 0 \quad \mbox{on } \ppp\OOO \right\},
\end{array}\right.
$$
where $\sigma$ is smooth, $\sigma\ge 0$ on $\ppp\OOO$,  
and the coefficients $a_{ij}$, $b_j$ and $c$ satisfy the conditions 
\eqref{(1.2)}.

Furthermore we recall that we number all the eigenvalues of $A$ as
$0 < \la_1 \le \la_2 \le \cdots$ with the multiplicities,
and we can choose linearly independent $\va_n$, $n\in \N$ such that 
$A\va_n = \la_n\va_n$.
Then $A_0^{\gamma}v = \sum_{n=1}^{\infty} \la_n^{\gamma} (v,\va_n)\va_n$
and 
$$
\mathcal{D}(A_0^{\gamma})
= \left\{ v\in L^2(\OOO): \thinspace
\sum_{n=1}^{\infty} \la_n^{2\gamma} (v,\va_n)^2 < \infty\right\}
$$
with
$$
\Vert A_0^{\gamma}v\Vert = \left( \sum_{n=1}^{\infty}
\la_n^{2\gamma} (v,\va_n)^2 \right)^{\frac{1}{2}}.
$$

Henceforth we choose a constant $\gamma$ such that 
$$
\frac{3}{4} < \gamma \le 1.
$$
For a semilinear term 
$f: \DDDA \longrightarrow L^2(\OOO)$, we assume that  
for some constant $m>0$, there exists a constant $C_f=C_f(m)>0$ such 
that
\begin{equation}
\label{(3.3)}
\left\{ \begin{array}{rl}
& \mbox{(i)} \Vert f(v)\Vert \le C_f, \quad 
\Vert f(v_1) - f(v_2)\Vert \le C_f\Vert v_1-v_2\Vert_{\DDDA} \\ 
& \mbox{if } \Vert v\Vert_{\DDD(A_0^{\gamma})}, 
\Vert v_1\Vert_{\DDD(A_0^{\gamma})}, 
\Vert v_2\Vert_{\DDD(A_0^{\gamma})}\le m \\
& \mbox{and}\\
& \mbox{(ii) there exists a constant} \ep \in \left(0, \frac{3}{4}\right)
\mbox{ such that }\\
& \Vert f(v)\Vert_{H^{2\ep}(\OOO)} \le C_f(m) \quad
\mbox{if } \Vert v\Vert_{\mathcal{D}(A_0^{\gamma})} \le m.
\end{array}\right.
\end{equation}
In what follows, $C>0$, $C_0, C_1 > 0$, etc., denote generic constants which 
are independent of the functions $u, v$, etc. under consideration, and 
we write 
$C_f$, $C(m)$ when we need to specify the dependence on related quantities.

Before we state and prove the main results of this section, let us discuss 
some examples of the source functions which satisfy the condition 
\eqref{(3.3)}. 

\begin{example}
\label{e1}
For $f\in C^1(\R)$, by setting $f(u):= f(u(x,t))$ for $(x,t) \in \QQQQ$,
we define 
$f: \DDD(A_0^{\gamma})\, \longrightarrow L^2(\OOO)$, 
$\frac{3}{4} < \gamma < 1$.
\end{example} 

Let us verify that the function from Example \ref{e1} satisfies the condition 
\eqref{(3.3)}. By the Sobolev embedding and the spatial dimensions 
$d \le 3$, we see that 
$\DDDA \subset H^{2\gamma}(\OOO) \subset C(\ooo{\OOO})$.  
Therefore $\Vert v\Vert_{\DDD(A_0^{\gamma})}\le m$ yields
$\vert v(x)\vert \le C_0m$ for $a \in \OOO$.  Hence, by using the mean
value theorem, the first condition in \eqref{(3.3)} is satisfied if
we choose $C_f(m) = \Vert f\Vert_{C^1[-C_0m,\,C_0m]}$.

Next we have to verify the second condition in \eqref{(3.3)}.
To this end, we prove the following lemma:

\begin{lemma}
\label{l11.1}
Let $f \in C^1[-C_0m,C_0m]$ and $0 < \ep < \frac{3}{4}$.
Then
$$
\Vert f(w)\Vert_{\DDD(A_0^{\ep})} \le C(1+m)\Vert f\Vert_{C^1[-C_0m,C_0m]}
\quad \mbox{if } \Vert w\Vert_{\DDDA} \le m.
$$
\end{lemma}

\begin{proof}
Because of the restrictions $0 < \ep < \frac{3}{4}$, we have the equality 
$\DDD(A_0^{\ep}) = H^{2\ep}(\OOO)$
(\cite{Fu}, \cite{Gri}).  
Therefore, since $f \in C^1[-C_0m,C_0m]$, the following norm estimates
in the Sobolev-Slobodecki space $H^{2\ep}(\OOO)$ are valid:
\begin{align*}
& \Vert f(w)\Vert^2_{H^{2\ep}(\OOO)}
= \Vert f(w)\Vert_{L^2(\OOO)}^2
+ \int_{\OOO}\int_{\OOO} \frac{\vert f(w(x)) - f(w(y))\vert^2}
{\vert x-y\vert^{d+4\ep}} dxdy\\
\le& \Vert f\Vert^2_{C^1[-C_0m,C_0m]}\vert \OOO\vert^2
+ \int_{\OOO}\int_{\OOO} \frac{\Vert f\Vert^2_{C^1[-C_0m,C_0m]}
\vert w(x) - w(y)\vert^2}{\vert x-y\vert^{d+4\ep}} dxdy\\
\le &C(1 + \Vert w\Vert^2_{H^{2\ep}(\OOO)})\Vert f\Vert^2_{C^1[-C_0m,C_0m]}
\le C(1+m^2)\Vert f\Vert^2_{C^1[-C_0m,C_0m]}.
\end{align*}
\end{proof}

%$\square$

The result formulated in Lemma \ref{l11.1} ensures the validity of the second  
condition from \eqref{(3.3)} for the function defined in Example \ref{e1}.

\begin{example}
\label{e2}
Let 
$$
f(v)(x):= \sum_{k=1}^d \mu_k(x)v(x)\ppp_kv(x),\quad x\in \OOO,
$$
where $\mu_k \in C^1(\ooo{\OOO})$.  In particular, this semilinear term
appears in the time-fractional Burgers equation $\pppa u = \ppp_x^2u - 
u\ppp_xu$ among its particular cases.
\end{example}

Let us verify that  the function defined in Example \ref{e2} satisfies the 
condition \eqref{(3.3)}. 
Under the conditions $\frac{3}{4} < \gamma < 1$,
the first condition from \eqref{(3.3)} is verified as follows.
For $v \in \mathcal{D}(A_0^{\gamma})$, it follows from 
$\frac{3}{4} < \gamma < 1$ that 
$\Vert \nabla v\Vert \le C\Vert v\Vert_{H^1(\OOO)}
\le C\Vert A_0^{\gamma}v\Vert$.
Moreover, $d=1,2,3$ and the Sobolev embedding imply  
$\Vert v\Vert_{C(\ooo{\OOO})} \le C\Vert v\Vert_{\DDDA}$.

Therefore, we obtain the desired inequalities:
\begin{align*}
& \Vert f(v)\Vert \le C\sum_{k=1}^d \Vert v\ppp_kv\Vert 
\le C\Vert v\Vert_{C(\ooo{\OOO})}\sum_{k=1}^d \Vert \ppp_kv\Vert\\
\le& C\Vert v\Vert_{C(\ooo{\OOO})}\Vert v\Vert_{\DDDA}
\le C\Vert v\Vert^2_{\DDDA} \le Cm^2 =: C_f(m)
\end{align*}
and
\begin{align*}
& \Vert f(v_1)-f(v_2)\Vert 
= \left\Vert \sum_{k=1}^d \mu_k(v_1-v_2)\ppp_kv_1
+ \mu_kv_2\ppp_k(v_1-v_2)\right\Vert\\
\le& C( \Vert v_1-v_2\Vert_{C(\ooo{\OOO})}\Vert \nabla v_1\Vert
+ \Vert v_2\Vert_{C(\ooo{\OOO})}\Vert \nabla (v_1-v_2)\Vert)\\
\le &C\max\{ \Vert v_1\Vert_{\DDDA},\, \Vert v_2\Vert_{\DDDA}\}
\Vert v_1 - v_2\Vert_{\DDDA}.
\end{align*}

Now we verify the second condition from \eqref{(3.3)}.  Since $\mu_k\in 
C^1(\ooo{\OOO})$, it suffices to prove that we can choose some constant
$C(m) > 0$ such that  
\begin{equation}
\label{(11.1)}
\Vert v\ppp_kv\Vert_{H^{2\ep}(\OOO)} \le C(m), \quad 
k=1,..., d \quad \mbox{if } \Vert v\Vert_{\DDD(A_0^{\gamma})} \le m.
\end{equation}
Since $d=1,2,3$ and $\gamma > \frac{3}{4}$, the Sobolev embedding
\begin{equation}
\label{(11.2)}
H^{2\gamma}(\OOO) \subset C^{\theta}(\ooo{\OOO})   
\end{equation}
holds true, where $\theta < 2\gamma - \frac{3}{2}$ (e.g., 
Theorem 1.4.4.1 (p.27) in \cite{Gri}).
In \eqref{(11.2)}, $C^{\theta}(\ooo{\OOO})$ denotes the Schauder space of 
uniform H\"older continuous functions on $\ooo{\OOO}$, which was defined 
at the end of the proof of Lemma \ref{l6.1}.

Now we choose small $\ep \in (0,1)$ such that
\begin{equation}
\label{(11.3)}
\ep < \min\left\{ \hhalf\theta, \, \frac{1}{4} \right\}.
\end{equation}
Now, using the inequalities
$$
\Vert v\Vert_{H^1(\OOO)} \le C\Vert v\Vert_{\DDD(A_0^{\gamma})}
\le Cm,
$$ 
we estimate 
\begin{align*}
& \Vert v\ppp_kv\Vert^2_{H^{2\ep}(\OOO)}
= \Vert v\ppp_k v\Vert^2_{L^2(\OOO)}
+ \int_{\OOO}\int_{\OOO}
\frac{\vert (v\ppp_k v)(x) - (v\ppp_k v)(y)\vert^2}
{\vert x-y\vert^{d+4\ep}} dxdy\\
=:& I_1+I_2.
\end{align*}
The inclusion \eqref{(11.2)} leads to the inequalities
\begin{align*}
&I_1 \le \Vert v\Vert^2_{C(\ooo{\OOO})}\Vert \ppp_kv\Vert^2_{L^2(\OOO)}
\le C\Vert v\Vert_{\DDD(A_0^{\gamma})}^2\Vert v\Vert_{H^1(\OOO)}^2\\
\le & C\Vert v\Vert^4_{\DDD(A_0^{\gamma})} \le Cm^4.
\end{align*}
As for the term $I_2$, since 
\begin{align*}
& \vert (v\ppp_kv)(x) - (v\ppp_kv)(y)\vert^2
\le (\vert v(x)(\ppp_kv(x) - \ppp_kv(y))\vert 
+ \vert v(x) - v(y)\vert \vert \ppp_kv(y)\vert)^2\\
\le & 2\Vert v\Vert^2_{C(\ooo{\OOO})}
\vert \ppp_kv(x) - \ppp_kv(y)\vert^2
+ 2\vert \ppp_kv(y)\vert^2 \vert v(x) - v(y)\vert^2,
\end{align*}
we obtain the estimate
\begin{align*}
& I_2 \le C\Vert v\Vert_{C(\ooo{\OOO})}^2
\int_{\OOO}\int_{\OOO}
\frac{\vert (\ppp_k v)(x) - (\ppp_k v)(y)\vert^2}
{\vert x-y\vert^{d+4\ep}} dxdy\\
+ & C\int_{\OOO} \int_{\OOO} \vert \ppp_kv(y)\vert^2 
 \frac{\vert v(x) - v(y)\vert^2}{\vert x-y\vert^{d+4\ep}} dxdy\\
=: & I_{21} + I_{22}.
\end{align*}
The Sobolev embedding yields the inequality $\Vert v\Vert_{C(\ooo{\OOO})}
\le C\Vert v\Vert_{\DDD(A_0^{\gamma})}$.
Moreover, since $0<\ep<\frac{1}{4}$ by the condition \eqref{(11.3)}, 
the inequality 
$\gamma > \frac{3}{4}$ leads to the estimates $2\ep + 1 < 2\frac{1}{4} + 1 
= \frac{3}{2} < 2\gamma$, so that we obtain the inequalities
$$
\Vert \ppp_kv\Vert_{H^{2\ep}(\OOO)} \le C\Vert v\Vert_{H^{2\ep+1}(\OOO)}
\le C\Vert v\Vert_{H^{2\gamma}(\OOO)} \le C\Vert v\Vert_{\DDD(A_0^{\gamma})}.
$$
Hence
$$
I_{21} \le C\Vert v\Vert^2_{C(\ooo{\OOO})}
\Vert \ppp_kv\Vert^2_{H^{2\ep}(\OOO)}
\le C\Vert v\Vert^2_{C(\ooo{\OOO})}
\Vert v\Vert^2_{\DDD(A_0^{\gamma})} \le Cm^4.
$$
Since \eqref{(11.2)} yields
$$
\vert v(x) - v(y)\vert \le C\Vert v\Vert_{\DDD(A_0^{\gamma})}
\vert x-y\vert^{\theta}, \quad x,y \in \ooo{\OOO},
$$
we obtain 
$$
I_{22} \le C\Vert v\Vert_{\DDD(A_0^{\gamma})}^2 \int_{\OOO}
\vert \ppp_kv(y)\vert^2\left(\int_{\OOO} \vert x-y\vert^{2\theta-d-4\ep}
dx \right) dy.
$$
Introducing the polar coordinate $r:= \vert x-y\vert$ for the integral in $x$ 
with a fixed $y\in \OOO$ and choosing $R_0>0$ sufficiently large such that 
$\max_{x,y\in \ooo{\OOO}} \vert x-y\vert < R_0$, we proceed as follows: 
\begin{align*}
& \int_{\OOO} \vert x-y\vert^{2\theta-d-4\ep}dx 
\le \int_{\vert x-y\vert < R_0} \vert x-y\vert^{2\theta-d-4\ep} dx\\
\le& C\int^{R_0}_0 r^{2\theta-d-4\ep}r^{d-1} dr
\le C\int^{R_0}_0 r^{(2\theta-4\ep)-1} dr =: C_1 < \infty.
\end{align*}
To derive the last inequality, we used the condition $2\theta-4\ep > 0$ that 
holds valid because of \eqref{(11.3)}.
Hence, we reach the final estimate
$$
I_{22} \le CC_1\Vert v\Vert^2_{\DDD(A_0^{\gamma})}\int_{\OOO}
\vert \ppp_kv(y)\vert^2 dy
\le C\Vert v\Vert_{\DDD(A_0^{\gamma})}^2\Vert v\Vert^2_{H^1(\OOO)}
\le Cm^4.
$$
Thus, the condition \eqref{(3.3)} is satisfied for the semilinear term 
$f$ defined in Example \ref{e2}.

Now we proceed to the derivation of some results regarding the local unique 
existence of a solution to the problem \eqref{(3.1)}.

\begin{theorem}
\label{t3.1}
Let a semilinear term $f$ satisfy the condition \eqref{(3.3)} with $m>0$ and 
$\Vert a\Vert_{\DDDA} \le m$.  

Then there exists a constant $T=T(m)>0$
such that the semilinear initial-boundary value problem \eqref{(3.1)} 
possesses a unique solution $u=u_a$ satisfying the inclusions 
\begin{equation}
\label{(3.4)}
u_a \in L^2(0,T;H^2(\OOO)) \cap C([0,T];\DDDA), \quad
u_a-a \in \HH(0,T;L^2(\OOO)).                      
\end{equation}

Moreover there exists a constant $C(m)>0$, such that
\begin{equation}
\label{(3.5)}
\Vert u_a - u_b\Vert_{L^2(0,T;H^2(\OOO))}
\le C\Vert a-b\Vert_{\DDDA}                  
\end{equation}
provided that $\Vert a\Vert_{\DDDA}, \Vert b\Vert_{\DDDA} \le m$.
\end{theorem}

In the theorem, the unique existence local in time of solution to an 
initial-boundary value problem for the semilinear 
time-fractional diffusion equation is established in the class 
\eqref{(3.4)}.
In the well-studied case of the 
conventional diffusion equation, i.e., $\alpha=1$ in \eqref{(3.1)},
one of the available methodologies for the proof of this result 
is the theory of the analytic semigroups.  For the case $0<\alpha<1$, 
our proof relies on 
a similar idea, that is, on employing the operators $S(t)$ and 
$K(t)$ defined by \eqref{(5.1)} and \eqref{(5.2)},
even if the operator  $S(t)$ does not possess any semigroup properties.
In several published works, 
a similar approach can be found for self-adjoint $A$ in the case where
the constant $C_f$ can be chosen independently of $m>0$.
We can refer for example to Section 6.4.1 of 
Chapter 6 in \cite{J}.  However, in many semilinear equations, 
the constant $C_f$ depends on $m$, in other words,
$f$ is not globally Lipschitz continuous, and 
we need more cares for the arguments.

The results formulated in this theorem are similar to the ones well-known for 
the parabolic equations, that is, $\alpha=1$ in \eqref{(3.1)}, and 
see, e.g.,  \cite{He} or \cite{Pa}.

\begin{proof}

First we recall the notation
\begin{equation}
\label{(7.1)}
Q(t)u(t) := \sum_{j=1}^d b_j(\cdot,t)\ppp_ju(\cdot,t) 
+ (c(\cdot,t)+c_0)u(\cdot,t) \quad \mbox{in $\OOO$}, \, 0<t<T.   
\end{equation}
Moreover, 
\begin{equation}
\label{(7.2)}
\Vert \Ahalf Q(t)u(t)\Vert \le C\Vert Q(t)u(t)\Vert_{H^1(\OOO)}
\le C\Vert Q(t)u(t)\Vert_{H^2(\OOO)} \le C\Vert A_0u(t)\Vert            
\end{equation}
for $u\in \mathcal{D}(A_0)$.

Similarly to \eqref{(5.7)}, we can formally rewrite the problem 
\eqref{(3.1)} as follows:
\begin{equation}
\label{(7.3)}
u(t) = S(t)a + \int^t_0 K(t-s)Q(s)u(s) ds 
\end{equation}
$$
+ \int^t_0 K(t-s)f(u(s)) ds,
\quad 0<t<T.                                        
$$

For a fixed $\gamma \in \left( \frac{3}{4}, \, 1\right)$ 
in the condition \eqref{(3.3)} and a fixed initial value  
$a \in \DDD(A_0^{\gamma})$, 
we define an operator $L: L^2(0,T;L^2(\OOO))
\, \longrightarrow \,  L^2(0,T;L^2(\OOO))$ by 
$$
(Lu)(t) := S(t)a + \int^t_0 K(t-s)Q(s)u(s) ds + \int^t_0 K(t-s)f(u(s)) ds,
\quad 0<t<T.
$$
Choosing a constant $m>0$ arbitrarily, we set 
\begin{equation}
\label{(7.4)}
V:= \{ v\in C([0,T];\DDD(A_0^{\gamma}));\, 
\Vert u - S(\cdot)a\Vert_{C([0,T];\DDD(A_0^{\gamma}))} \le m\}.  
\end{equation}

Then we prove the following lemma:

\begin{lemma}
\label{l7.1}
Let $H \in C([0,T]; L^2(\OOO))$.  Then
$$
\int^t_0 A_0^{\gamma}K(t-s)H(s) ds \in C([0,T];L^2(\OOO)),
$$
that is,
$$
\int^t_0 K(t-s)H(s) ds \in C([0,T];\mathcal{D}(A_0^{\gamma})).
$$
\end{lemma}

\begin{proof}
Let $0< \eta < t \le T$.  We have
\begin{align*}
& \int^t_0 A_0^{\gamma}K(t-s)H(s) ds - \int^{\eta}_0 A_0^{\gamma}K(\eta-s)
H(s) ds\\
=& \int^t_0 A_0^{\gamma}K(s)H(t-s) ds - \int^{\eta}_0 A_0^{\gamma}K(s)
H(\eta-s) ds\\
=& \int^t_{\eta} A_0^{\gamma}K(s)H(t-s) ds 
+ \int^{\eta}_0 A_0^{\gamma}K(s) (H(t-s) - H(\eta-s)) ds\\
=: &I_1 + I_2.
\end{align*}
For the first integral, by \eqref{(5.3)} and $\gamma < 1$, we have
\begin{align*}
&\Vert I_1\Vert
\le C\int^t_{\eta} s^{\alpha(1-\gamma)-1} \max_{0\le s\le t}
\Vert H(t-s)\Vert ds\\
\le& C\Vert H\Vert_{C([0,T];\DDD(A_0^{\gamma}))}
\frac{t^{\alpha(1-\gamma)} - \eta^{\alpha(1-\gamma)}}
{\alpha(1-\gamma)} \, \longrightarrow \, 0
\end{align*}
as $\eta \uparrow t$.  

%Next let $\chi_{(0,\eta)}(s) = 
%\left\{\begin{array}{rl}
%1, \,\, &0<s<\eta,\\
%0, \,\, & \mbox{otherwise}.
%\end{array}\right.$
Next we have the estimates
\begin{align*}
& \Vert I_2\Vert
=\left\Vert \int^{\eta}_0 A_0^{\gamma}K(s) (H(t-s) - H(\eta-s)) ds\right\Vert\\
\le & C\int^{\eta}_0 s^{(1-\gamma)\alpha-1} 
\max_{0\le \eta \le t \le T} \Vert H(t-s) - H(\eta-s)\Vert ds.
\end{align*}
Hence, by $H\in C([0,T];L^2(\OOO))$, the function
$$
\vert s^{(1-\gamma)\alpha-1} \vert
\max_{0\le \eta \le t \le T} \Vert H(t-s) - H(\eta-s)\Vert 
$$
is an integrable function with respect to $s \in (0,\eta)$ and 
$$
\lim_{\eta \uparrow t} s^{(1-\gamma)\alpha-1} 
\max_{0\le \eta \le t \le T} \Vert H(t-s) - H(\eta-s)\Vert = 0
$$
for almost all $s\in (0,\eta)$.
Therefore, the Lebesgue convergence theorem implies
$\lim_{\eta \uparrow t} \Vert I_2\Vert = 0$.  Thus the proof of 
Lemma \ref{l7.1} is completed.
\end{proof}

Now we proceed to the proof of Theorem \ref{t3.1}.
In view of \eqref{(5.1)}, the inclusion $a \in \DDDA$ 
implies
\begin{equation}
\label{(7.5)}
S(t)a \in C([0,T];\DDDA).                       
\end{equation}
Indeed, 
\begin{align*}
& \Vert A_0^{\gamma}(S(t)a - S(s)a)\Vert^2
= \Vert S(t)(A_0^{\gamma}a) - S(s)(A_0^{\gamma}a)\Vert^2\\
=& \sum_{n=1}^{\infty} \vert E_{\alpha,1}(-\la_nt^{\alpha})
- E_{\alpha,1}(-\la_ns^{\alpha})\vert^2 \vert (A_0^{\gamma}a,\va_n)\vert^2.
\end{align*}
Applying the Lebesgue convergence theorem, in view of 
$$
\vert E_{\alpha,1}(-\la_nt^{\alpha})\vert \le \frac{C}{1+\la_nt^{\alpha}}
\quad \mbox{for all $n\in \N$ and $t>0$}
$$
(e.g., Theorem 1.6 (p.35) in \cite{Po}), we can verify \eqref{(7.5)}.

Because of the condition \eqref{(3.3)} and $\mathcal{D}(A_0^{\gamma})
\subset H^1(\OOO)$, 
for $v \in C([0,T];\DDD(A_0^{\gamma}))$, we have 
$f(v) \in C([0,T];L^2(\OOO))$ and $Qv \in C([0,T];L^2(\OOO))$.
Applying Lemma \ref{l7.1}, in view of \eqref{(7.5)}, we see
\begin{equation}
\label{(7.6)}
Lv \in C([0,T]; \DDDA) \quad \mbox{for } v \in C([0,T]; \DDDA).
\end{equation}

For the further proof, we have to prove 
the following properties (i) and (ii) which 
are valid for sufficiently small $T>0$:
\\
(i) $LV \subset V$, $V$ being the set defined by \eqref{(7.4)}.
\\
(ii) There exists a constant $\rho \in (0,1)$ such that 
$$
\Vert Lu_1 - Lu_2\Vert_{C[0,T];\DDDA)} 
\le \rho\Vert u_1 - u_2\Vert_{C[0,T];\DDDA)}, \quad 0<t<T
$$
for $u_1, u_2 \in V$.
\\
{\bf Proof of (i).}
\\
Let $u \in V$.  Then, the inclusion \eqref{(7.6)} says that 
$Lu \in C([0,T];\DDDA)$.
Now we consider 
\begin{equation}
\label{(7.7)}
A_0^{\gamma}(Lu(t) - S(t)a)
= \int^t_0 A_0^{\gamma}K(t-s)Q(s)u(s) ds + \int^t_0 A_0^{\gamma}K(t-s)
f(u(s)) ds.                                     
\end{equation}
For any $u \in V$, using 
$$
\Vert a\Vert_{\mathcal{D}(A_0^{\gamma})}
= \Vert A_0^{\gamma}a\Vert \le m, \quad
\Vert u-S(\cdot)a\Vert_{C([0,T];\mathcal{D}(A_0^{\gamma}))} \le m,
$$
we obtain
\begin{equation}
\label{(7.8)}
\Vert u(t)\Vert_{\DDDA} \le m + \Vert A_0^{\gamma}S(t)a\Vert 
= m + \Vert S(t)A_0^{\gamma}a\Vert 
\le m + C_1m =: C_2m.
\end{equation}
The first condition in \eqref{(3.3)} implies that
\begin{equation}
\label{(7.9)}
\Vert f(u(t)) \Vert \le C_f(C_2m) \quad \mbox{for all $u\in V$
and $0<t<T$.}    
\end{equation}
Furthermore, \eqref{(7.8)} implies 
\begin{equation}
\label{(7.10)}
\Vert Q(s)u(s)\Vert \le C\Vert u(s)\Vert_{H^1(\OOO)}
\le C_3\Vert A_0^{\hhalf}u(s)\Vert \le C_4\Vert u(t)\Vert_{\DDD(A_0^{\gamma})}
\le C_4C_2m                         
\end{equation}
because of the inequalities $\gamma \ge \hhalf$.
Applying \eqref{(7.9)} and \eqref{(7.10)} in \eqref{(7.7)}, by means of \eqref{(5.3)},
we obtain the estimates
\begin{align*}
& \Vert Lu(t) - S(t)a \Vert_{\DDDA} \\
= & \left\Vert \int^t_0 A_0^{\gamma}K(t-s)Q(s)u(s) ds
+ \int^t_0 A_0^{\gamma}K(t-s)f(u(s)) ds\right\Vert\\
\le & C\int^t_0 (t-s)^{(1-\gamma)\alpha-1} (C_2C_4m + C_f(C_2m))ds
\le C_5\frac{t^{(1-\gamma)\alpha}}{(1-\gamma)\alpha}
\le C_5\frac{T^{(1-\gamma)\alpha}}{(1-\gamma)\alpha}.
\end{align*}
The constant $C_5>0$ depends on $m>0$ but is independent 
on $T>0$.
Therefore, choosing $T>0$ sufficiently small, we complete the proof of 
the property (i).  $\square$
\\
{\bf Proof of (ii).}
\\
Estimate \eqref{(7.8)} yields that $\Vert u_1(t)\Vert_{\DDDA}$,
$\Vert u_2(t)\Vert_{\DDDA} \le C_2m$ for each $u_1, u_2 \in V$.
Therefore, the condition \eqref{(3.3)} yields
$$
\Vert f(u_1(s)) - f(u_2(s))\Vert \le C_f(C_2m)\Vert u_1(s)-u_2(s)\Vert
_{\DDDA}, \quad 0<s<T.
$$
Hence, by \eqref{(7.10)} we have the following chain of estimates:
\begin{align*}
& \Vert Lu_1(t) - Lu_2(t)\Vert_{\DDDA} \\
= & \biggl\Vert \int^t_0 A_0^{\gamma}K(t-s)Q(s)(u_1(s)-u_2(s)) ds \\
+  & \int^t_0 A_0^{\gamma}K(t-s)(f(u_1(s)) - f(u_2(s))) ds\biggr\Vert\\
\le & C\int^t_0 (t-s)^{\alpha(1-\gamma)-1}\Vert (u_1-u_2)(s)\Vert
_{\DDDA} ds\\
+& C_f(C_2m)\int^t_0 (t-s)^{\alpha(1-\gamma)-1}\Vert (u_1-u_2)(s)\Vert
_{\DDDA} ds\\
\le & C_6T^{\alpha(1-\gamma)}
\sup_{0<s<T} \Vert u_1(s) - u_2(s)\Vert_{\DDDA}.
\end{align*}
In the last inequality, the constant $C_6>0$ is independent of $T>0$, and thus 
one can further choose a sufficiently small constant $T>0$ satisfying 
the inequality
$$
\rho:= C_6T^{\alpha(1-\gamma)} < 1.
$$
The proof of the property (ii) is completed.
$\square$

The properties (i) and (ii) allow to apply  the contraction theorem to the 
equation $u=Lu$, which says that this equation  has a unique 
solution $u\in V$  
for $0<t<T$.  This solution $u\in C([0,T];\mathcal{D}(A_0^{\gamma}))$ 
satisfies the estimates 
\eqref{(7.8)} and the equation
\begin{equation}
\label{(7.11)}
u(t) = S(t)a + \int^t_0 K(t-s)Q(s)u(s) ds 
\end{equation}
$$
+ \int^t_0 K(t-s)f(u(s)) ds, \quad 0<t<T.          
$$

Next we have to prove the rest inclusions of \eqref{(3.4)}.
In the condition \eqref{(3.3)}, we can choose $\ep>0$ such that $0<\ep<\hhalf$.
By the equation \eqref{(7.11)}, we obtain
\begin{align*}
& A_0u(t) = A_0^{1-\gamma}S(t)A_0^{\gamma}a
+ \int^t_0 A_0^{1-\ep}K(t-s)A_0^{\ep}Q(s)u(s) ds\\
+& \int^t_0 A_0^{1-\ep}K(t-s)A_0^{\ep}f(u(s)) ds, \quad 0<t<T.
\end{align*}
In the last equation, $a\in \DDD(A_0^{\gamma})$ and the estimates 
$$
\Vert A_0^{\ep}Q(s)u(s)\Vert \le C\Vert A_0^{\hhalf}(Q(s)u(s))\Vert 
\le C\Vert u(s)\Vert_{H^2(\OOO)} \le C\Vert A_0u(s)\Vert
$$
hold true because of \eqref{(7.2)}.  Furthermore, by the second condition 
in \eqref{(3.3)}, the inequality \eqref{(7.8)} 
yields the estimate $\Vert A_0^{\ep}f(u(s))\Vert 
\le C_f(C_2m)$ by the second condition from \eqref{(3.3)}.  
Thus we reach the chain of the inequalities
\begin{align*}
& \Vert A_0u(t)\Vert
\le Ct^{-\alpha(1-\gamma)}\Vert A_0^{\gamma}a\Vert 
+ C\int^t_0 (t-s)^{\alpha\ep-1}\Vert A_0u(s)\Vert ds\\
+ & C\int^t_0 (t-s)^{\alpha\ep-1}C_f(C_2m) ds\\
\le & Ct^{-\alpha(1-\gamma)}\Vert A_0^{\gamma}a\Vert + C_f(C_2m)
+ C\int^t_0 (t-s)^{\alpha\ep-1}\Vert A_0u(s)\Vert ds, \quad
0<t<T.
\end{align*}
By $0<\alpha<1$, we have $-\alpha(1-\gamma) > -1$ and thus the inclusion
$t^{-\alpha(1-\gamma)} \in L^1(0,T)$. 
Application of the generalized Gronwall inequality yields the estimates
\begin{align*}
& \Vert A_0u(t)\Vert
\le (Ct^{-\alpha(1-\gamma)}\Vert A_0^{\gamma}a\Vert + C_f(C_2m)) \\ 
+ & C\int^t_0 (t-s)^{\alpha\ep-1}(s^{-\alpha(1-\gamma)}
\Vert A_0^{\gamma}a\Vert + C_f(C_2m)) ds\\
\le&  Ct^{-\alpha(1-\gamma)}\Vert A_0^{\gamma}a\Vert 
+ C_f(C_2m)
+ (\Vert A_0^{\gamma}a\Vert + C_f(C_2m))t^{\alpha(\ep-(1-\gamma))},
\quad 0<t<T.
\end{align*}
Therefore, noting that $-\alpha(1-\gamma) < \alpha(\ep-(1-\gamma))$, we have
$$
\Vert A_0u(t)\Vert \le C_7(1+T^{\alpha\ep})(t^{-\alpha(1-\gamma)} + 1),
\quad 0<t<T,
$$
where $C_7>0$ depends on $\Vert A_0^{\gamma}a\Vert$ and
$C_f(C_2m)$, $\alpha$, $\ep$.
For $\hhalf < \gamma \le 1$, we can directly verify that 
$-2\alpha(1-\gamma) > -1$, 
% by $-\alpha(1-\gamma) < \alpha(\ep-(1-\gamma))$, 
so that 
$\int^T_0 \Vert A_0u(t)\Vert^2 dt < \infty$, that is, the inclusion
\begin{equation}
\label{(7.12)}
u \in L^2(0,T;H^2(\OOO))                         
\end{equation}
holds true. 

It remains to prove $u-a \in H_{\alpha}(0,T;L^2(\OOO))$.  
The inequality \eqref{(7.9)} implies the inclusion 
$f(u) \in L^2(0,T;L^2(\OOO))$.

The estimate \eqref{(7.10)} and the inclusion \eqref{(7.12)} result in the 
inclusion $Qu \in L^2(0,T;L^2(\OOO))$.  
Therefore, Lemma \ref{l5.1} (ii) yields the inclusion
$$
\int^t_0 K(t-s)Q(s)u(s) ds + \int^t_0 K(t-s)f(u(s)) ds
\in \HH(0,T;L^2(\OOO)).
$$
Applying Lemma \ref{l5.1} (iii) to the equation \eqref{(7.3)}, we reach  
the inclusion $u-a 
\in \HH(0,T;L^2(\OOO))$, which completes the proof of the relation 
\eqref{(3.4)} from the theorem.

Finally, we have to prove the estimate \eqref{(3.5)}. By the construction of 
the solutions $u_a, u_b$ as the fixed points, we have the inequalities
\begin{equation}\label{(5.20)}
\Vert u_a\Vert_{L^2(0,T;H^2(\OOO))}, \,
\Vert u_b\Vert_{L^2(0,T;H^2(\OOO))} \le C(m).
\end{equation}
On the other hand, 
\begin{align*}
&u_a(t) - u_b(t)
= S(t)a - S(t)b + \int^t_0 K(t-s)Q(s)(u_a-u_b)(s) ds\\
+& \int^t_0 K(t-s)(f(u_a(s)) - f(u_b(s))) ds, \quad 0<t<T. 
\end{align*}
In view of \eqref{(5.20)}, we can use the condition \eqref{(3.3)} and apply
the generalized Gronwall inequality.
Further details of the derivations are similar to the ones employed
in the proof of Theorem \ref{t2.1}
and we omit them here. 
Thus, the proof of Theorem \ref{t3.1} is completed.
\end{proof}

%%%%%%%%%%%%%%%% Section 6 %%%%%%%%%%%%

\section{Comparison principles for the semilinear time-fractional diffusion 
equations} 
\label{sec6}

\setcounter{section}{6}

In this section, we derive comparison principles for semilinear
time-fractional diffusion equations.  To this end, we have to 
restrict semilinear terms $f$ which depend only on $x$ and $u(x,t)$,
but should not depend on its derivatives.
 
We introduce the class of semilinear terms in terms of 
$\mathfrak{f} \in C^1(\ooo{\OOO}\times [-m,\, m])$.  
For such a function, we can naturally define a mapping
$f: \{v\in \DDD(A_0^{\gamma}):\, \Vert v\Vert_{\DDD(A_0^{\gamma})} \le m\} 
\, \longrightarrow \, L^2(\OOO)$ by 
\begin{equation}
\label{(3.6a)}
f(v):= \mathfrak{f}(x,\, v(x)), \quad x\in \OOO,\, 0<t<T.      
\end{equation}
In what follows, we identify $f$ with $\mathfrak{f}$ through the relation 
\eqref{(3.6a)}.
For a fixed constant $M>0$, we set 
\begin{equation}
\label{(3.7a)}
\mathcal{F}_M:= \{ f\in C^1(\ooo{\OOO}\times [-m,\,m]);\,
\Vert f\Vert_{C^1(\ooo{\OOO}\times [-m,\,m])} \le M\}.      
\end{equation}

Now we are ready to formulate and to prove the first comparison principle for 
the initial-boundary value problems for semilinear time-fractional 
diffusion equation \eqref{(3.1)}.

\begin{theorem}
\label{t3.2}
For $f_1, f_2 \in \mathcal{F}_M$ and $a_1, a_2 \in H^1_0(\OOO)$,  
we assume that there exist solutions $u(f_k,a_k)$, $k=1,2$ to the 
initial-boundary value problem \eqref{(3.1)}
with 
the semilinear terms $f_k$, $k=1,2$ and the initial values $a_k$, $k=1,2$, 
respectively, which 
satisfy the inclusions \eqref{(3.4)} and the estimates
\begin{equation}
\label{(3.8)}
\vert u(f_k,a_k)(x,t)\vert \le m, \quad x\in \OOO,\, 
0<t<T, \, k=1,2.                              
\end{equation}
If $f_1 \ge f_2$ on $\OOO \times (-m, m)$ and $a_1\ge a_2$ in $\OOO$, then 
\begin{equation}
\label{(3.9)}
u(f_1,a_1) \ge u(f_2,a_2) \quad \mbox{in } \OOO \times (0,T).
\end{equation}
\end{theorem}

\begin{proof}

Henceforth, for simplicity, we denote $f_k(x,u_j(x,t))$ by 
$f_k(u_j)$ and
$u(f_k,a_k)$ by $u_k$ for $j,k=1,2$.  

For $d=1,2,3$, application of the Sobolev embedding yields  
$\DDD(A_0^{\gamma}) \subset H^{2\gamma}(\OOO) \subset C(\ooo{\OOO})$
by $\gamma > \frac{3}{4}$.
Hence, \eqref{(3.4)} implies
\begin{equation}
\label{(8.1)}
u_1, u_2 \in C(\ooo{\OOO} \times [0,T]).     
\end{equation}
On the other hand, we have the representations
\begin{align*}
& f_1(u_1(x,t)) - f_2(u_2(x,t))
= f_1(u_1(x,t)) - f_1(u_2(x,t)) + (f_1-f_2)(u_2(x,t))\\
=& g(x,t)(u_1(x,t) - u_2(x,t)) + H(x,t), \quad
(x,t) \in \QQQQ,
\end{align*}
where we set
$$
g(x,t) := 
\left\{ \begin{array}{rl}
& \frac{f_1(u_1(x,t)) - f_1(u_2(x,t))}
{u_1(x,t) - u_2(x,t)} \quad \mbox{if } u_1(x,t) \ne u_2(x,t),\\
& f_1'(u_1(x,t)) \quad \mbox{if } u_1(x,t) = u_2(x,t),
\end{array}\right.
$$
and
$$
H(x,t) := (f_1-f_2)(u_2(x,t)), \quad (x,t) \in \QQQQ.
$$
By the inclusions \eqref{(8.1)} and $f \in C^1(\ooo{\OOO} \times [-m,m])$, we 
can verify that $g, H \in C(\ooo{\OOO}\times [0,T])$.  
Furthermore $f_1\ge f_2$ implies $H \ge 0$ in $\QQQQ$.

Setting now $y:= u_1 - u_2$ and $a:= a_1-a_2$, the function $y$ is a solution 
to the following initial-boundary value problem:
\begin{equation}
\label{(8.2)}
\left\{ \begin{array}{rl}
& \pppa (y-a) + Ay - g(x,t)y = H \ge 0 \quad 
\mbox{in } \QQQQ, \\
& \NUNU y + \sigma y = 0 \quad \mbox{on } \ppp\OOO \times (0,T), \\
& y(x,\cdot) - a(x) \in \HH(0,T) \quad \mbox{for almost all } x\in \OOO.
\end{array}\right.
\end{equation}
Since $g \in C(\ooo{\OOO}\times [0,T])$ does not in general satisfy 
the regularity condition \eqref{(1.2)}, we cannot directly apply Theorem 
\ref{t2.2}.
Thus,  we first approximate $g$ by $g_n \in C^{\infty}(\ooo{\OOO}\times
[0,T])$, $n\in \N$ such that $g_n \longrightarrow g$ in 
$C(\ooo{\OOO}\times [0,T])$ as $n\to \infty$ and consider a sequence of 
the following problems:
\begin{equation}
\label{(8.3)}
\left\{ \begin{array}{rl}
& \pppa (y_n-a) + Ay_n - g_n(x,t)y_n = H \quad 
\mbox{in } \QQQQ, \\
& \NUNU y_n + \sigma y_n = 0 \quad \mbox{on } \ppp\OOO \times (0,T), \\
& y_n(x,\cdot) - a(x) \in \HH(0,T) \quad \mbox{for almost all } x\in \OOO.
\end{array}\right.
\end{equation}
Then Theorem \ref{t2.1} yields that for any $n\in \N$ there exists a unique 
solution 
$y_n \in L^2(0,T;H^2(\OOO))$ to the problem \eqref{(8.3)} such that 
$y_n-a \in \HH(0,T;L^2(\OOO))$.  Moreover, since 
$g_n$ satisfies the regularity condition \eqref{(1.2)}, Theorem \ref{t2.2} 
yields that 
\begin{equation}
\label{(8.4)}
y_n\ge 0 \quad \mbox{in } \QQQQ \mbox{ for each } n\in \N.
\end{equation}
On the other hand, setting $z_n:= y_n-y$ in $\QQQQ$, the equations 
\eqref{(8.2)} and 
\eqref{(8.3)} allow to characterize $z_n$ as solution to the problem
$$
\left\{ \begin{array}{rl}
& \pppa z_n + Az_n = (g_n-g)y + g_nz_n \quad 
\mbox{in } \QQQQ, \\
& \NUNU z_n + \sigma z_n = 0 \quad \mbox{on } \ppp\OOO \times (0,T), \\
& z_n(x,\cdot)  \in \HH(0,T) \quad \mbox{for almost all } x\in \OOO.
\end{array}\right.
$$
Similarly to our treatment of the problem \eqref{(5.6)}, we rewrite the first 
equation in the form
\begin{align*}
&\pppa z_n + A_0z_n\\
= &(g_n-g)y + g_nz_n + \sum_{j=1}^d b_j(t)\ppp_jz_n
+ (c_0+c(t))z_n \quad \mbox{in $\QQQQ$},
\end{align*}
and repeat the same arguments as the ones employed for \eqref{(5.7)} to obtain 
the  following chain of the estimates:
\begin{align*}
& \Vert A_0^{\hhalf} z_n(t)\Vert
= \biggl\Vert \int^t_0 A_0^{\hhalf}K(t-s)(g_n-g)(s)y(s) ds\\
+& \int^t_0 A_0^{\hhalf}K(t-s)
\left( (g_n(s)+c_0+c(s))z_n(s) 
+ \sum_{j=1}^d b_j(s)\ppp_jz_n(s) \right) ds \biggr\Vert\\
\le& C\int^t_0 (t-s)^{\hhalf\alpha-1} \Vert g_n-g\Vert
_{C(\ooo{\OOO}\times [0,T])}\Vert y(s)\Vert ds\\
+& \int^t_0 (t-s)^{\hhalf\alpha-1} (\Vert g_n\Vert_{C(\ooo{\OOO}\times [0,T])}
+ c_0 + \Vert c\Vert_{C(\ooo{\OOO}\times [0,T])} \\
+ & \max_{1\le j \le d} \Vert b_j\Vert_{C(\ooo{\OOO}\times [0,T])})
\Vert z_n(s)\Vert_{H^1(\OOO)} ds \\
\le & Ct^{\hhalf\alpha} \Vert g_n-g\Vert_{C(\ooo{\OOO}\times [0,T])}
+ C\int^t_0 (t-s)^{\hhalf\alpha-1} \Vert A_0^{\hhalf}z_n(s)\Vert ds\\
\le & C_1\Vert g_n-g\Vert_{C(\ooo{\OOO}\times [0,T])}
+ C_1\int^t_0 (t-s)^{\hhalf\alpha-1} \Vert A_0^{\hhalf}z_n(s)\Vert ds,
\quad 0<t<T.
\end{align*}
Here the constant $C_1>0$ depends also on $\Vert y\Vert_{C([0,T];L^2(\OOO))}$.
In the derivations, we used the relation 
$\sup_{n\in \N}\Vert g_n\Vert_{C(\ooo{\OOO}\times [0,T])} < \infty$,
which is justified by  
$g_n \longrightarrow g$ in $C(\ooo{\OOO}\times [0,T])$ and 
$y:= u_1-u_2 \in C(\ooo{\OOO}\times [0,T])$ that follows from 
\eqref{(8.1)}.

Therefore, the generalized Gronwall inequality yields
\begin{align*}
& \Vert A_0^{\hhalf} z_n(t)\Vert
\le C_2\Vert g_n - g\Vert_{C(\ooo{\OOO}\times [0,T])}
+ C_2\int^t_0 (t-s)^{\hhalf\alpha-1} 
\Vert g_n - g\Vert_{C(\ooo{\OOO}\times [0,T])} ds\\
\le &C_3\Vert g_n - g\Vert_{C(\ooo{\OOO}\times [0,T])},
\quad 0<t<T,\, n\in \N.
\end{align*}
Hence, $z_n \longrightarrow 0$ in $L^{\infty}(0,T;\DDD(A_0^{\hhalf}))$
as $n\to \infty$, and thus  
$y_n \longrightarrow y$ in $L^{\infty}(0,T;\DDD(A_0^{\hhalf}))$
as $n\to \infty$.  Because of the relation
\eqref{(8.4)}, we reach the inequality $y\ge 0$ in $\QQQQ$.
Thus the proof of Theorem \ref{t3.2} is completed.
\end{proof}

We note that in the formulation of Theorem \ref{t3.2}, the boundedness 
condition 
\eqref{(3.8)} has been assumed. However, this condition is not automatically 
guaranteed for solutions to the initial-boundary value problem \eqref{(3.1)}. 
In the rest of this section, we present another comparison principle in terms 
of the upper and lower solutions, which does not require this 
boundedness condition. 

\begin{definition}
\label{duls}
Let $f \in \mathcal{F}_M$.
The functions $\UPP$ and $\LOP$ satisfying \eqref{(3.4)} are 
called an upper solution and a lower solution for the solution 
$u$ to the problem \eqref{(3.1)}, respectively, 
if
\begin{equation}
\label{(3.10)}
\left\{ \begin{array}{rl}
& \pppa (\UPP - \ooo{a}) + A\UPP \ge f(\UPP) \quad 
\mbox{in } \OOO\times (0,T),               \\
& \ooo{a}(x) \ge u(x,0), \quad x\in \OOO,\\
& \NUNU \UPP + \sigma\UPP = 0 \quad \mbox{on } \ppp\OOO\times (0,T),\\
& \UPP - \ooo{a} \in \HH(0,T;L^2(\OOO)),
\end{array}\right.
\end{equation}
and
\begin{equation}
\label{(3.11)}
\left\{ \begin{array}{rl}
& \pppa (\LOP - \underline{a}) + A\LOP \le f(\LOP) \quad \mbox{in } \OOO\times 
(0,T),\\
& \underline{a}(x) \le u(x,0), \quad x\in \OOO,\\
& \NUNU \LOP + \sigma\LOP = 0 \quad \mbox{on } \ppp\OOO\times (0,T),\\
& \LOP - \underline{a} \in \HH(0,T;L^2(\OOO)),
\end{array}\right.
\end{equation}
with $\ooo{a}, \underline{a} \in L^2(\OOO)$.
\end{definition}

We recall that $\mathcal{F}_M$ is defined by \eqref{(3.7a)}.

Then the following result holds true:
\begin{theorem}
\label{t3.3}
We arbitrarily choose $T>0$.  We assume that there exist 
an upper solution $\ooo{u}$ and a lower solution $\underline{u}$
 to the problem \eqref{(3.1)} and that  
$\vert \UPP(x,t)\vert$, $\vert \LOP(x,t)\vert \le m$ for $x\in \OOO$ and 
$0<t<T$.

Then there exists a unique solution $u_a$ to the problem \eqref{(3.1)} in 
the class \eqref{(3.4)} and
$$
\LOP(x,t) \le u_a(x,t) \le \UPP(x,t), \quad x\in \OOO,\, 0<t<T.
$$
\end{theorem}

\begin{proof}

For the functions $\UPPY, \LOWY \in L^2(0,T;L^2(\OOO))$ satisfying 
the condition 
$\LOWY(x,t) \le \UPPY(x,t)$ for $(x,t) \in \QQQQ$, we set 
$$
[\LOWY,\, \UPPY]:= \{ y\in L^2(0,T;L^2(\OOO));\, 
\LOWY(x,t) \le y(x,t) \le \UPPY(x,t)\},
$$
where the inequalities hold true for almost all  $(x,t) \in \QQQQ$. 

The key element of our proof of Theorem \ref{t3.3} is the following fixed 
point theorem in an ordered Banach space:

\begin{lemma}
\label{l9.1}
Let an operator $L: [\LOWY,\, \UPPY] \, \longrightarrow [\LOWY,\, \UPPY]
\subset L^2(0,T;L^2(\OOO))$ be compact and increasing, that is, 
$Lv \ge Lw$ in $\QQQQ$ if $v, w\in [\LOWY,\, \UPPY]$ and $v\ge w$ in 
$\QQQQ$.

Then the operator $L$ possesses fixed points $y^*$, $y_* \in 
[\LOWY,\, \UPPY]$ such that
$$
y_* = \lim_{k\to\infty} L^k\LOWY, \quad
y^* = \lim_{k\to\infty} L^k\UPPY \quad \mbox{in } L^2(0,T;L^2(\OOO)).
$$
\end{lemma}

The proof of Lemma \ref{l9.1} can be found in \cite{Am}.
Lemma \ref{l9.1} implies the following result:

\begin{lemma}
\label{l9.2}
In Lemma \ref{l9.1}, we further assume that the operator $L$ possesses 
a unique fixed point.
Then
$$
\LOWY \le y \le \UPPY \quad \mbox{in } \QQQQ.
$$
\end{lemma}

\begin{proof}
{\bf Proof of Lemma \ref{l9.2}.}
Indeed, since $L[\LOWY,\, \UPPY] \subset [\LOWY,\,\UPPY]$, we see
that $\LOWY \le L\LOWY$ and $L\UPPY \le \UPPY$.
Moreover, since the operator $L$ is increasing, the inequality $L\LOWY \le
L\UPPY$ holds true.  Hence,
$$
\LOWY \le L\LOWY \le L\UPPY \le \UPPY \quad \mbox{in } \QQQQ.
$$
Applying the operator $L$ once again, we obtain
$$
\LOWY \le L\LOWY \le L^2\LOWY \le L^2\UPPY \le L\UPPY \le \UPPY.
$$
Repeating the same arguments, we reach the inequalities 
$$
\LOWY \le L^k\LOWY \le L^k\UPPY \le \UPPY \quad \mbox{in } \QQQQ \mbox{ for }
k\in \N.
$$
By the uniqueness of the fixed point, letting $k \to \infty$, we see
that $\LOWY \le y \le \UPPY$.  
\end{proof}
%$\square$
% Thus the proof of Lemma \ref{l9.2} follows 
% directly from Lemma \ref{l9.1}.

Now we proceed to the proof of Theorem \ref{t3.3} and 
 define an operator $L$ from $[\UPP,\, \LOP] \subset  
L^2(0,T;L^2(\OOO))$ to itself as follows: 
For $u\in [\underline{u}, \, \ooo{u}]$, the function $v:=Lu$ satisfies
\begin{equation}
\label{(9.1)}
\left\{ \begin{array}{rl}
& \pppa (v-a) + Av + (M+1)v = (M+1)u + f(u) \quad \mbox{in } \QQQQ,\\
& v - a \in \HH(0,T;L^2(\OOO)), \\
& \NUNU v + \sigma v = 0 \quad \mbox{on } \ppp\OOO \times (0,T),
\end{array}\right.
\end{equation}
where $M>0$ is the constant in the definition \eqref{(3.7a)} of the set 
$\mathcal{F}_M$. 

We verify that $L: [\underline{u},\, \ooo{u}] \subset 
L^2(0,T;L^2(\OOO)) \, \longrightarrow \, L^2(0,T;L^2(\OOO))$ is well-defined.
Indeed let $u \in [\underline{u},\, \ooo{u}]$.
By $\vert \underline{u}(x,t)\vert$, $\vert \ooo{u}(x,t)\vert \le m$ for
all $x\in \ooo{\OOO}$ and $0\le t \le T$, we see that 
$\vert u(x,t)\vert \le m$ for almost all $x\in \OOO$ and 
$t\in (0,T)$.  Therefore $f(u(x,t))$ can be defined for almost all
$(x,t) \in \OOO \times (0,T)$, and $f(u(x,t)) \in L^{\infty}(\OOO\times 
(0,T))$.  Hence, since $(M+1)u + f(u) 
\in L^{\infty}(\OOO\times (0,T)) \subset L^2(0,T;L^2(\OOO))$, we apply 
Theorem \ref{t2.1} to \eqref{(9.1)}, and we verify that there exists a unique solution 
$v \in L^2(0,T;H^2(\OOO))$ to \eqref{(9.1)} satisfying 
$v-a \in H_{\alpha}(0,T;L^2(\OOO))$, and 
\begin{equation}\label{(6.12)}
\Vert v-a\Vert_{H_{\alpha}(0,T;L^2(\OOO))} 
+ \Vert v\Vert_{L^2(0,T;H^2(\OOO))}
\end{equation}
\begin{align*}
\le& C(\Vert a\Vert_{H^1(\OOO)} + \Vert (M+1)u+f(u)\Vert
_{L^2(0,T;L^2(\OOO))}\\
\le& C(\Vert a\Vert_{H^1(\OOO)} 
+ ((M+1)m+M)T^{\hhalf}\vert \OOO\vert^{\hhalf})
=: C_4
\end{align*}
by using 
\begin{align*}
& \Vert (M+1)u+f(u)\Vert_{L^{\infty}(\OOO\times (0,T))} \\
\le & (M+1)\Vert u\Vert_{L^{\infty}(\OOO\times (0,T))}
+ \Vert f(u)\Vert_{L^{\infty}(\OOO\times (0,T))}
\le m(M+1)+M.
\end{align*}
Thus $L$ is well-defined and the estimate \eqref{(6.12)} holds.
$\square$

For the operator $L$, we will verify the following properties (i) - (iii):

(i) $L$ is a compact operator from $[\LOP, \, \UPP] \subset 
L^2(0,T;L^2(\OOO))$ into itself.  In other words, the set 
$L\, [\LOP, \, \UPP]$ is relatively compact in $L^2(0,T;L^2(\OOO))$.
\\
{\bf Verification of (i).}
Let $u \in [\LOP, \, \UPP]$.  Then the estimate \eqref{(6.12)} means 
$$
\Vert Lu\Vert_{L^2(0,T;H^2(\OOO))} 
+ \Vert Lu - a\Vert_{\HH(0,T;L^2(\OOO))} 
\le C_4.
$$
For $a\in H^1(\OOO)$, we deduce that
$$
\Vert Lu-a\Vert_{L^2(0,T;H^1(\OOO))}
+ \Vert Lu - a \Vert_{\HH(0,T;L^2(\OOO))}
\le C_4 + \sqrt{T}\Vert a\Vert_{H^1(\OOO)}.
$$
Since  
$$
L^2(0,T;H^1(\OOO)) \cap \HH(0,T;L^2(\OOO))
\subset L^2(0,T;L^2(\OOO))
$$ 
is a compact embedding (see, e.g., \cite{Te}), we obtain that 
$v \longrightarrow Lv - a$ is compact from $[\UPP,\, \LOP]
\subset L^2(0,T;L^2(\OOO))$ to $L^2(0,T;L^2(\OOO))$.
Thus the compactness of the operator $L$ is proved.
$\square$

(ii) $Lv \ge Lw$ in $\QQQQ$ if $v\ge w$ in $\QQQQ$ and
$v,w \in [\LOP, \,\UPP]$.
\\
{\bf Verification of (ii).}
Setting $y:= v-w$ and $z:= Lv - Lw$  
% similarly to the arguments for derivation of \eqref{(8.2)} 
and 
applying the mean value theorem for  $f(v) - f(w)$, we obtain the 
representation
$$
\pppa (z-a) + Az + (M+1)z=(M+1)y 
+ f'(\mu(x,t))y \quad \mbox{in } \QQQQ,
$$
where $\mu(x,t)$ is a number between $v(x,t)$ and $u(x,t)$.

The inclusion  $f \in \mathcal{F}_M$ implies $\vert f'(\mu(x,t))\vert \le M$.
Therefore, in view of $y \ge 0$ in $\QQQQ$, we estimate
$$
(M+1)y + f'(\mu(x,t))y \ge (M+1-M)y(x,t) \ge 0 \quad \mbox{in } \QQQQ.
$$
Hence, the function $z$ satisfies
$$
\left\{ \begin{array}{rl}
& \pppa (z-a) + Az + (M+1)z \ge 0 \quad \mbox{in } \QQQQ,\\
& \NUNU z + \sigma z = 0 \quad \mbox{on } \ppp\OOO \times (0,T).
\end{array}\right.
$$
Theorem \ref{t2.1} implies that $z\in L^2(0,T;H^2(\OOO))$ and 
$z-a \in \HH(0,T;L^2(\OOO))$ and Theorem \ref{t2.2} yields the inequality
$z\ge 0$ in $\QQQQ$.  Therefore $Lv \ge Lw$ in $\QQQQ$.
$\square$

(iii) $L[\LOP,\,\UPP] \subset [\LOP,\, \UPP]$.
\\
{\bf Verification of (iii).}
In order to prove the above inclusion,  we show that $\LOP \le L\LOP$ and
$L\UPP \le \UPP$ in $\QQQQ$.  Indeed, let $\LOP \le u \le \UPP$ in 
$\QQQQ$.  
Since $L$ is increasing, we have the inequalities $L\LOP \le Lu \le L\UPP$.
Hence, $\LOP \le Lu \le \UPP$ in $\QQQQ$, which proves (iii).
$\square$
\\
{\bf Proof of $\LOP \le L\LOP$.}
\\
Setting $v:= L\LOP$, we have
$$
\left\{ \begin{array}{rl}
& \pppa (v - a) 
+ Av + (M+1)v = (M+1)\LOP + f(\LOP) \quad \mbox{in } \QQQQ,\\
& v - a \in \HH(0,T;L^2(\OOO)), \\
& \NUNU v + \sigma u = 0 \quad \mbox{on } \ppp\OOO \times (0,T).
\end{array}\right.
$$
On the other hand, since $\LOP$ is a lower solution, we see that 
$$
\left\{ \begin{array}{rl}
& \pppa (\LOP - \underline{a}) + A\LOP + (M+1)\LOP 
\le (M+1)\LOP + f(\LOP) \quad \mbox{in } \QQQQ,\\
& \LOP - \underline{a} \in \HH(0,T;L^2(\OOO)), \\
& \NUNU \LOP + \sigma \LOP = 0 \quad \mbox{on } \ppp\OOO \times (0,T).
\end{array}\right.
$$
Therefore, $z:= v - \LOP = L\LOP - \LOP$ satisfies
$$
\left\{ \begin{array}{rl}
& \pppa (z - (a -\underline{a})) + Az + (M+1)z \ge 0 
\quad \mbox{in } \QQQQ,\\
& z - (a-\underline{a}) \in \HH(0,T;L^2(\OOO)), \\
& \NUNU z + \sigma z = 0 \quad \mbox{on } \ppp\OOO \times (0,T),\\
& a-\underline{a} \ge 0 \quad \mbox{in } \OOO.
\end{array}\right.
$$
According to Theorem \ref{t2.2}, the inequality 
$z\ge 0$ holds true in $\QQQQ$, that is, $L\LOP \ge \LOP$ in $\QQQQ$.
By similar arguments, we prove that $L\UPP \le \UPP$ in $\QQQQ$.
Thus the property (iii) is proved.
$\square$

% Finally, in Theorem \ref{t2.1}, we already proved that the uniqueness of the 
% fixed point or we can repeat the same argument based on (5.7) and the 
% generalized Gronwall inequality.
The properties (i)-(iii) allow to apply Lemmas \ref{l9.1} and \ref{l9.2},  which
completes the proof of 
Theorem \ref{t3.3}.
\end{proof}

Theorem \ref{t3.3} ensures that, for a given $T>0$, it is sufficient 
to determine an upper and a lower solutions to the problem \eqref{(3.1)} 
in order to guarantee the existence of the solution.
This technique is called a monotonicity method which is also related to the 
Perron method or the concept of viscosity solutions.
For applications of the monotonicity method to parabolic equations, 
that is, $\alpha=1$ in \eqref{(3.1)}, we refer, e.g., to
\cite{Am}, \cite{Ke},  \cite{Pao1}, and \cite{Pao2}.  Theorem \ref{t3.3} 
asserts a corresponding result for the case 
$0<\alpha<1$ and readily yields the following statement:

\begin{proposition}
\label{p3.1}
Let $u$ be a solution to the problem \eqref{(3.1)} satisfying  
\eqref{(3.4)} and $\vert u(x,t)\vert \le m$ for $(x,t)\in \OOO\times (0,T)$ with some
constant $m>0$. 
Then:

(i) For any  upper solution $\UPP$ to the problem \eqref{(3.1)} satisfying
\eqref{(3.4)} and $\vert \UPP(x,t)\vert \le m$ for $(x,t) \in \OOO\times (0,T)$, 
we have
$$
u(x,t) \le \UPP(x,t) \quad \mbox{for all $(x,t)\in \OOO\times (0,T)$}.
$$ 

(ii) For any lower solution $\LOP$ to the problem \eqref{(3.1)} 
satisfying \eqref{(3.4)} and $\vert \LOP(x,t)\vert \le m$ for 
$(x,t) \in \OOO\times (0,T)$, we have
$$
\LOP(x,t) \le u(x,t) \quad \mbox{for all $(x,t)\in \OOO\times (0,T)$}.
$$ 
\end{proposition}

In the rest of this section, we focus on a special case of the operator $-A$ 
in the form
\begin{equation}
\label{div}
\left\{\begin{array}{rl}
& -Av := \mbox{div}\, (p(x)\nabla v(x)) + c(x)v, \quad x\in \OOO,\\
& \DDD(A) = \{ v\in H^2(\OOO); \, \NUNU v + \sigma v = 0 \quad 
\mbox{on } \ppp\OOO \},
\end{array}\right.
\end{equation}
where $p\in C^2(\ooo{\OOO})$, $p>0$ on $\ooo{\OOO}$ and 
$c\in C^2(\ooo{\OOO})$, $c\le 0$, $c \not\equiv 0$ in $\OOO$.

% Then by $c\le 0, \not\equiv 0$ we see that the minimum eigenvalue $\la_1$ is 
% positive and we can prove that $\la_1$ is simple and we can choose 
% the corresponding eigenfunction $\va_1$ which is positive on $\ooo{\OOO}$.
% The properties on $\la_1$ and $\va_1$ are formulated as 
% Lemma \ref{l9.3} presented below.

For this operator, Theorem \ref{t3.3} leads to the following result:

\begin{proposition}
\label{p3.2}
We assume that $f \in C^1(\R)$ is 
a monotone decreasing function and that 
there exists a solution $u_{\infty}\in H^2(\OOO)$ to the 
boundary-value problem
\begin{equation}
\label{(3.12)}
\left\{ \begin{array}{rl}
& Au_{\infty}(x) = f(u_{\infty}(x)), \quad x \in \OOO, \\
& \NUNU u_{\infty} + \sigma u_{\infty} = 0 \quad \mbox{on } \ppp\OOO,\\
& \vert u_{\infty}(x)\vert \le m, \quad x\in \OOO,
\end{array}\right.
\end{equation}
with the operator $-A$ defined by \eqref{div}.

Then for any $T>0$ and any $a\in \DDD(A_0^{\gamma})$, 
the initial-boundary value problem \eqref{(3.1)} possesses
a unique solution $u_a$ in the class \eqref{(3.4)} such that the inequality 
\begin{equation}
\label{(3.14)}
\vert u_a(x,t) - u_{\infty}(x)\vert \le CE_{\alpha,1}(-\la_1 t^{\alpha})
\vert \va_1(x)\vert,
\quad x\in \OOO, \, t\in (0,T)                               
\end{equation}
holds true with  a certain constant $C>0$. 
In particular, we have
\begin{equation}
\label{(3.15)}
\vert u_a(x,t) - u_{\infty}(x)\vert \le Ct^{-\alpha}\vert \va_1(x)\vert
\quad \mbox{for all } x\in \OOO \mbox{ as } t\to \infty.   
\end{equation}
\end{proposition}

\begin{proof}

The proof is similar to the one of 
Proposition 1.4 (pp.25-26) in \cite{Ke}, see also 
\cite{Pao2}.

Since $c(x)\le 0$, $\not\equiv 0$ on $\ooo{\OOO}$ and belongs to the space 
$C^2(\ooo{\OOO})$, all the eigenvalues of the operator 
$A$ are known to be positive:
$$
0 < \la_1 \le \la_2 \le \cdots \longrightarrow \infty,
$$
where $\la_k$, $k\in \N$ are numbered according to their multiplicities.
Now we choose an eigenfuction $\va_1$ for the minimum eigenvalue $\la_1$ 
with $\Vert \va_1\Vert = 1$.  Then we apply the following well-known result.
\begin{lemma}\label{l9.3}
The multiplicity of the eigenvalue $\la_1$ is one and
$\va_1(x) > 0$ for $x\in \ooo{\OOO}$ or
$\va_1(x) < 0$ for $x\in \ooo{\OOO}$.
\end{lemma}

As for the proof, we can refer, for example, to Lemma 1.4 (p.96) and
Theorem 1.2 (p.97) in \cite{Pao2}.

Thus we choose such an eigenfunction $\va_1$ such that 
$\va_1(x)>0$ for $x \in \ooo{\OOO}$.
Then we can find a sufficiently large constant $M_1>0$ such that 
\begin{equation}
\label{(9.2)}
u_{\infty}(x) - M_1\va_1(x) \le 
a(x) \le u_{\infty}(x) + M_1\va_1(x), \quad x\in \ooo{\OOO}.
\end{equation}
We set 
$$
\left\{ \begin{array}{rl}
&\ooo{a}(x):= \UPP(x,0) = u_{\infty}(x) + M_1\va_1(x), \\
&\underline{a}(x):= u_{\infty}(x) - M_1\va_1(x), \quad x\in \OOO
\end{array}\right.
$$
and
$$
\left\{ \begin{array}{rl}
& \LOP(x,t):= u_{\infty}(x) - M_1\MLONE(-\la_1 t^{\alpha})\va_1(x),  \\
& \UPP(x,t):= u_{\infty}(x) + M_1\MLONE(-\la_1 t^{\alpha})\va_1(x),
\quad (x,t) \in \QQQQ,
\end{array}\right.,
$$
where the Mittag-Leffler function $\MLONE$ is defined by 
$$
\MLONE(z) = \sum_{k=0}^{\infty} \frac{z^k}{\Gamma(\alpha k + 1)}, \quad
z\in \C
$$
(see, e.g., \cite{Po}).
Let $T>0$ be arbitrary.  Recalling that 
$$
-Aw(x) = \mbox{div}\, (p(x)\nabla w) + c(x)w, \quad x\in \OOO
$$
with $c(x) \le 0$, $\not\equiv 0$ on $\ooo{\OOO}$ and 
$\pppa(\MLONE(-\la_1t^{\alpha})-1) = -\la_1\MLONE(-\la_1t^{\alpha})$,
we have $\LOP - \lowa \in H_{\alpha}(0,T;L^2(\OOO))$, 
$$
\pppa (\LOP - \lowa) =  \la_1M_1\MLONE(-\la_1t^{\alpha})\va_1(x)
$$
and
$$
A\LOP = Au_{\infty} - \la_1M_1\MLONE(-\la_1t^{\alpha})\va_1(x),
\quad x\in \OOO,\, 0<t<T,
$$
so that 
$$
\pppa (\LOP - \lowa) + A\LOP = Au_{\infty} = f(u_{\infty}) \quad
\mbox{in }\QQQQ.
$$
It is known that $\MLONE(-\la_1t^{\alpha}) \ge 0$ for $t\ge 0$
(see, e.g., \cite{GKMR}) and 
$$
\LOP(x,t) \le u_{\infty}(x), \quad (x,t) \in \OOO\times (0,T).
$$
Since $f$ is decreasing, we obtain
$f(u_{\infty}(x)) \le f(\LOP(x,t))$.  Therefore,  
$$
\pppa (\LOP - \lowa) + A\LOP \le f(\LOP) \quad
\mbox{in }\QQQQ.
$$
Since $\NUNU\LOP + \sigma\LOP = 0$ on 
$\ppp\OOO \times (0,T)$, the function 
$\LOP$ is a lower solution to the initial-boundary value problem \eqref{(3.1)} 
with the operator $-A$ defined by \eqref{div}. Similarly we can verify 
that $\UPP$ is an upper solution.  Hence, Theorem \ref{t3.3} yields that 
there exists a unique solution $u\in L^2(0,T;H^2(\OOO))$ such that 
$$
u-a \in \HH(0,T;L^2(\OOO)), \quad 
\LOP(x,t) \le u(x,t) \le \UPP(x,t), \quad (x,t) \in \QQQQ,
$$
that is,
$$
\vert u(x,t) - u_{\infty}(x)\vert \le M_1\MLONE(-\la_1t^{\alpha})
\va_1(x), \quad x\in \OOO,\, 0<t<T.
$$
Since $T>0$ is arbitrary and $\vert \MLONE(-\la_1t^{\alpha})\vert
\le \frac{C}{t^{\alpha}}$ as $t \to \infty$ (Theorem 1.6 (p.35) in \cite{Po}), 
the proof of Proposition \ref{p3.2} is completed.
\end{proof}

In general, the solution to the boundary-value problem \eqref{(3.12)} may  
be not unique. However,
under our assumptions,  the estimate \eqref{(3.14)} implies that $u_{\infty}$ 
is uniquely 
determined as the limit of $u_a$ as $t \to \infty$.
In the homogeneous  case ($f\equiv 0$), the asymptotic behavior 
of the solution $u_a$ is known  (see, e.g., \cite{KRY}, \cite{SY}, \cite{Ya}): 
$\Vert u_a(\cdot, t)\Vert = O(t^{-\alpha})$ as $t \to \infty$ that is 
the same as in the relation \eqref{(3.15)}.
\\

We close this section with two examples of the results presented in 
Proposition \ref{p3.1} in the case
$A = -\Delta$ with the homogeneous Neumann boundary condition, that is,
$\sigma=0$ on $\ppp\OOO$. 
In this case, the boundary condition can be represented as follows:
$$
\NUNU v = \ppp_{\nu}v := \nabla v\cdot \nu \quad \mbox{on } \ppp\OOO.
$$

We assume that there exists a solution $u_a$ to the initial-boundary value 
problem \eqref{(3.1)} in the class \eqref{(3.4)}. 
Then we present two examples, where we can estimate
$u(x,t) - a(x)$ for small $t>0$. 

\begin{example}
\label{e3}
Let $\sigma \equiv 0$ on $\ppp\OOO$ and 
$$
f(\eta) = -\frac{\eta}{1+\vert \eta\vert}, \quad \eta \in \R.
$$
We assume that 
\begin{equation}
\label{(3.16)}
a\in C^2(\ooo{\OOO}), \quad a\ge 0 \quad \mbox{in $\OOO$}, \quad
\ppp_{\nu} a = 0 \quad\mbox{on }\ppp\OOO. 
\end{equation}
Then we choose a constant $\rho>0$ so large that
\begin{equation}\label{(6.18a)}
\Delta a(x) \le \Gamma(\alpha+1)\rho, \quad x\in \ooo{\OOO}.
\end{equation}
Then
\begin{equation}
\label{(3.17)}
0 \le u_a(x,t) - a(x) \le \rho t^{\alpha}, \quad x\in \OOO,\,
0\le t < T.                         
\end{equation}
\end{example}
{\bf Proof of \eqref{(3.17)}.}
We set 
$$
\LOP(x,t) = \underline{a}(x) = 0, \quad
\UPP(x,t) = \rho t^{\alpha} + a(x), \quad 
x\in \OOO, \, 0<t<T.
$$
Then the equation 
$$
\pppa (\LOP - \underline{a}) - \Delta \LOP - f(\LOP) = 0
$$
holds trivially, and so 
$$
\left\{ \begin{array}{rl}
& \pppa (\LOP - 0) - \Delta \LOP \le f(\LOP)\quad \mbox{in $\OOO
\times (0,T)$}, \\
& 0=\underline{a}(x) \le u(x,0) = a(x), \quad x\in \OOO,\\
& \ppp_{\nu}\LOP = 0 \quad \mbox{on $\ppp\OOO \times (0,T)$},
\end{array}\right.
$$
which means that $\LOP(x,t) \equiv 0$ is a lower solution.
Moreover, the conditions \eqref{(3.16)} and \ref{(6.18a)} imply 
$$
\left\{ \begin{array}{rl}
& \pppa (\UPP(x,t) - a(x)) - \Delta \UPP - f(\UPP)
= \rho\Gamma(\alpha+1) - \Delta a 
+ \frac{a(x)+\rho t^{\alpha}}{1+\vert a(x)+\rho t^{\alpha}\vert} \\
\ge &\rho \Gamma(\alpha+1) - \Delta a(x) \ge 0, \quad x\in \OOO, \, 0<t<T, \\
& \UPP(x,0) = a(x), \quad x\in \OOO,\\
& \ppp_{\nu} \UPP = 0\quad \mbox{on } \ppp\OOO \times (0,T). 
\end{array}\right.
$$
Therefore, $\UPP$ and $\LOP$ defined above are an upper and a lower solutions,
respectively, and Proposition \ref{p3.1} leads to the estimate \eqref{(3.17)}.
$\square$

We note that for $\alpha=1$ the equation considered in Example \ref{e3} is a
semilinear parabolic equation
$$
\ppp_tu - \Delta u = -\frac{u}{1+\vert u\vert},
$$
which is a governing equation for an enzyme process. The monotonicity method 
is well-known for this type of semilinear parabolic equations 
(e.g., \cite{Ke} and \cite{Pao1}).

\begin{example}
\label{e4}
Let $f \in C^1(0,\infty)$  be an increasing function 
in $(0, \infty)$.
We assume that 
\begin{equation}
\label{(3.18)}
a \in C^2(\ooo{\OOO}), \qquad \ppp_{\nu}a = 0 \quad \mbox{on } \ppp\OOO.  
\end{equation}
Then for an arbitrary but fixed $\ep \in (0, \alpha)$, there exists 
$T_1=T_1(\ep) > 0$ such that 
\begin{equation}
\label{(3.19)}
u_a(x,t) - a(x) \le t^{\alpha-\ep}, \quad 0<t<T_1.     
\end{equation}
\end{example}

In order to prove the inequality \eqref{(3.19)}, 
we set $\UPP(x,t) = t^{\alpha-\ep} + a(x)$ for $(x,t) \in 
\OOO\times (0,T_1)$ with $T_1$ which will be selected later.  Then 
$$
\pppa (\UPP-a) = \frac{\Gamma(\alpha-\ep+1)}{\Gamma(-\ep+1)}t^{-\ep}
$$
and $\Delta \UPP = \Delta a$, and thus
$$
\pppa (\UPP - a) - \Delta \UPP 
= \frac{\Gamma(\alpha-\ep+1)}{\Gamma(-\ep+1)}t^{-\ep} - \Delta a,
\quad x\in \OOO, \, 0<t<T_1.
$$
Moreover, setting $M_1:= \max_{x\in \ooo{\OOO}} a(x)$, we deduce
\begin{equation}\label{(6.23)}
f(\UPP(x,t)) = f(t^{\alpha-\ep}+a(x))
\le f(T_1^{\alpha-\ep}+M_1), \quad x\in \OOO,\, 0<t<T_1.
\end{equation}
We note that $T_1^{-\ep}$ can be arbitrarily large if we choose 
$T_1>0$ sufficiently small.  Therefore, in view of $\alpha-\ep>0$, 
we can choose $T_1 > 0$ sufficiently small 
such that 
\begin{equation}
\label{(3.20)}
\frac{\Gamma(\alpha-\ep+1)}{\Gamma(-\ep+1)}T_1^{-\ep}
\ge f(T_1^{\alpha-\ep} + M_1) + \Delta a(x), \quad x\in \ooo{\OOO}.
\end{equation}
For this $T_1$, in terms of \eqref{(6.23)} and \eqref{(3.20)}, we can readily see 
$$
\pppa (\UPP-a) - \Delta \UPP \ge f(\UPP) \quad \mbox{in } \OOO\times (0,T_1).
$$
Since $\UPP(x,0) = a(x) = u_a(x,0)$ and $\ppp_{\nu}\UPP = 0$ on 
$\ppp\OOO \times (0,T_1)$, the function  $\UPP$ is an upper solution for
$0<t<T_1$, which yields the inequality \eqref{(3.19)}.

Next we look for a suitable  lower solution. 
In addition to the conditions \eqref{(3.18)}, we assume that there exists
a constant $\delta_1>0$ such that 
\begin{equation}\label{(6.23a)}
a(x) \ge \delta_1 > 0 \quad \mbox{for all $x\in \ooo{\OOO}$}.
\end{equation}
Then we set  
\begin{equation}
\label{(3.21)}
M_2 := \max_{x\in \ooo{\OOO}} (-\Delta a(x)).     
\end{equation}
The lower solution is introduced in the form
$$
\LOP(x,t) = a(x) - \rho t^{\alpha}, \quad x \in \OOO,\,
0<t < \left( \frac{\delta_1}{2\rho}\right)^{\frac{1}{\alpha}}=: T_2,
$$
where $\delta_1>0$ is the constant in the condition \eqref{(6.23a)} and 
the constant $\rho>0$ will be selected later.
Then, \eqref{(6.23a)} implies
$$
\LOP(x,t) \ge \frac{\delta_1}{2}, \quad x\in \OOO,\, 0<t<T_2.
$$
Moreover, the inequalities 
$$
\pppa (\LOP - a) - \Delta \LOP = -\Gamma(\alpha+1)\rho - \Delta a
\le -\Gamma(\alpha+1)\rho + M_2
$$
and
$$
f(\LOP(x,t)) = f(a(x) - \rho t^{\alpha})
\ge f(\delta_1 - \rho T_2^{\alpha}) 
= f\left( \delta_1 - \frac{\delta_1}{2}\right) 
= f\left(\frac{\delta_1}{2}\right)
$$
hold true.
Now we choose a constant $\rho$ sufficiently large such that 
$T_2:= \left( \frac{\delta_1}{2\rho}\right)^{\frac{1}{\alpha}}$
is sufficiently small and moreover
$$
\rho \ge \frac{1}{\Gamma(\alpha+1)}
\left(M_2 - f\left(\frac{\delta_1}{2}\right)\right).
$$
Then 
$$
\pppa (\LOP - a) - \Delta \LOP \le f(\LOP) \quad \mbox{in $\OOO\times 
(0,T_2)$}.
$$
Since $\ppp_{\nu}\LOP = 0$ on $\ppp\OOO \times (0,T_2)$ and 
$\LOP(x,0) = a(x)$ for $x\in \OOO$, the function
$\LOP$ is a lower solution, which leads to the estimate
\begin{equation}
\label{(3.22)}
a(x) - \frac{1}{\Gamma(\alpha+1)}
\left(M_2 - f\left(\frac{\delta_1}{2}\right)\right)t^{\alpha} \le u_a(x,t) 
\end{equation}
for $0 < t < T_2$, where $T_2>0$ is sufficiently small.
\\

Summarizing the results presented above,  under the conditions  \eqref{(3.18)} 
and \eqref{(6.23a)}, for any $\ep \in (0,\, \alpha)$, there exists 
a constant $T_{\ep,a} > 0$ such that  
\begin{equation}
\label{(3.23)}
-\frac{1}{\Gamma(\alpha+1)}\left( M_2 - f\left( \frac{\delta_1}{2}\right)
\right)t^{\alpha} \le u_a(x,t) - a(x)
\end{equation}
$$
\le t^{\alpha-\ep}, \quad x\in \OOO,\, 0<t<T_{\ep,a}.                     
$$
$\square$

Finally let us emphasize that the construction of upper and lower solutions is 
not unique.  For example,
the right-hand side of the inequality \eqref{(3.23)} has the form 
$t^{\alpha-\ep}$. However, under the conditions \eqref{(3.18)} we can obtain 
a different estimate as follows: For a sufficiently small $T_3>0$, 
we can choose 
a large constant $M_3=M_3(T_3)>0$ such that 
\begin{equation}\label{(6.27a)}
\left\{ \begin{array}{rl}
&\Vert \Delta a\Vert_{C(\ooo{\OOO})} \le \frac{1}{2}M_3\Gamma(\alpha+1),\\
& f(M_3T_3^{\alpha}+\Vert a\Vert_{C(\ooo{\OOO})})
\le \frac{1}{2}M_3\Gamma(\alpha+1).
\end{array}\right.
\end{equation}
We remark that sufficiently large $M_3(T)>0$ can satisfy \eqref{(6.27a)}
for small $T>0$, because 
$\lim_{T\downarrow 0} f(M_3T^{\alpha}+\Vert a\Vert_{C(\ooo{\OOO})})
= f(\Vert a\Vert_{C(\ooo{\OOO})})$.

Then the inequality 
\begin{equation}
\label{(3.24)}
u_a(x,t) - a(x) \le M_3(T_3)t^{\alpha},\quad x\in \OOO, \, 0<t<T_3
\end{equation}
holds true. 
\\
{\bf Verification of \eqref{(3.24)}.} 
We set $\UPP(x,t) := a(x) + M_3t^{\alpha}$.  Then the first condition 
in \eqref{(6.27a)} yields 
$$
\pppa (\UPP-a) - \Delta \UPP = M_3\Gamma(\alpha+1) - \Delta a
\ge \frac{1}{2}M_3\Gamma(\alpha+1).
$$
Moreover, since $f$ is monotone increasing, we have the estimate
$$
f(\UPP) = f(a(x)+M_3t^{\alpha}) 
\le f(\Vert a\Vert_{C(\ooo{\OOO})}+ M_3T_3^{\alpha})
$$
for $(x,t) \in\OOO \times (0,T_3)$.
Therefore, by means of \eqref{(6.27a)}, we obtain
$$
\pppa (\UPP-a) - \Delta \UPP \ge \frac{1}{2}M_3\Gamma(\alpha+1)
\ge f(M_3T_3^{\alpha}+\Vert a\Vert_{C(\ooo{\OOO})})
\ge f(\UPP) \quad \mbox{in $\OOO \times (0,T_3)$},
$$
which means that $\UPP$ is another upper solution and 
thus the inequality \eqref{(3.24)} is verified.
$\square$

\section{Concluding remarks}
\label{sec7}
\setcounter{section}{7}

In this article, for derivation of the comparison principles and 
the monotonicity method, we considered a class of 
solutions $u\in L^2(0,T;H^2(\OOO))$ under the initial condition 
$u-a \in \HH(0,T;L^2(\OOO))$.
However, our main results remain valid for the solutions from a larger space 
of functions, for instance, for mild solutions to the semilinear 
equation \eqref{(3.1)} which  
are defined as solutions to the equation \eqref{(7.3)}.

In order to prove the comparison principles for the  mild solutions, 
one needs pointwise arguments in Lemma \ref{l4.2} and thus 
 approximations of the
solutions by smooth solutions $u\in C([0,T];C^2(\ooo{\OOO}))$
satisfying $t^{1-\alpha}\ppp_tu \in C([0,T];C(\ooo{\OOO}))$ has to be 
considered. In other words, we first 
have to  prove the comparison principles
for smooth functions belonging to such a smoother class 
and then extend them to the desired class of solutions. 

In our discussions, we restricted ourselves to the case  of the homogeneous 
boundary conditions.  However, in the definitions
\eqref{(3.10)} and \eqref{(3.11)} of the upper and lower solutions, 
it is natural to  replace 
the homogeneous boundary conditions by the inequalities 
$\NUNU \UPP + \sigma \UPP \ge 0$ and $\NUNU \LOP + \sigma\LOP \le 0$
on $\ppp\OOO \times (0,T)$.  To this end, the unique existence results for the 
initial-boundary value problems with non-homogeneous boundary conditions are 
needed, which are not sufficiently studied. 
These problems will be considered elsewhere.

For the sake of simplicity of arguments, in this article we 
discussed only the case of the spatial dimensions 
$d=1, 2, 3$.  The case $d\ge 4$ can be treated similarly, but some stronger 
regularity conditions
for the coefficients of the differential operators and more involved estimates 
are needed.

\section*{Appendix: Proof of existence of a function satisfying  the conditions (3.4)}
\label{sec8}
\setcounter{section}{8}

{\bf First Step.}
\\
In this step, we prove
\begin{lemma}\label{lem1}
We assume \eqref{(1.2)} and $M:= \min_{(x,t)\in \ooo{\OOO}\times [0,T]} b_0(x,t)>0$ 
is sufficiently large.
Then there exists a constant $\kappa_1>0$ such that
\begin{equation}\label{1}
(A_1(t)v,\, v) \ge \kappa_1\Vert v\Vert_{H^1(\OOO)}^2
\end{equation} 
$\forall v \in H^2(\OOO)$ satisfying 
$\ppp_{\nu_A}v + \sigma v = 0$ on $\ppp\OOO$ and $\forall t \in [0,T]$.
\end{lemma}
In particular, Lemma \ref{lem1} implies that all the eigenvalues of $A_0$ 
defined by \eqref{(3.2)} are positive if the constant $c_0>0$ is sufficiently large.
Henceforth we write $b=(b_1,..., b_d)$.

\begin{proof}
By using \eqref{(1.2)} and $\NUNU v + \sigma v = 0$ on $\ppp\OOO$, 
integration by parts yields
\begin{align*}
& 0 = (A_1(t)v,v)\\
= &-\int_{\OOO} \sum_{i,j=1}^d \ppp_i(a_{ij}(x)\ppp_jv)v dx
 - \frac{1}{2}\int_{\OOO} \sum_{j=1}^d b_j(x,t)\ppp_j(\vert v\vert^2) dx
+ \int_{\OOO} b_0(x,t) \vert v\vert^2 dx\\
=& \int_{\OOO} \sum_{i,j=1}^d a_{ij}(x)(\ppp_iv)(\ppp_jv) dx
- \int_{\ppp\OOO} (\ppp_{\nu_A}v)v dS\\
+ & \frac{1}{2}\int_{\OOO} (\mbox{div}\, b)\vert v\vert^2 dx
- \frac{1}{2}\int_{\ppp\OOO} (b\cdot \nu)\vert v\vert^2 dS
+ \int_{\OOO} b_0(x,t) \vert v\vert^2 dx\\
\ge & \kappa \int_{\OOO} \vert \nabla v\vert^2 dx 
+ \int_{\OOO} \left(\min_{(x,t)\in\ooo{\OOO}\times [0,T]} b_0(x,t)
- \frac{1}{2}\vert \mbox{div}\, b\vert \right)\vert v\vert^2 dx\\
+& \int_{\ppp\OOO} \left( \sigma - \frac{1}{2}(b\cdot \nu)\right)
\vert v\vert^2 dS
\end{align*}
\begin{equation}\label{(2)}
\ge \kappa \int_{\OOO} \vert \nabla v\vert^2 dx 
+ \left( M - \frac{1}{2}\Vert \mbox{div}\, b\Vert_{C(\ooo{\OOO}\times [0,T])}
\right)\int_{\OOO} \vert v\vert^2 dx - C\int_{\OOO} \vert v\vert^2 dS.
\end{equation}
Here and henceforth $C>0$, $C_{\ep}, C_{\delta} > 0$, etc. denote generic 
constants which are independent of $v$.  

By the trace theorem (e.g., Theorem 9.4 (pp.41-42) in \cite{LM}), for
$\delta \in (0, \frac{1}{2}]$, there exists
a constant $C_{\delta}>0$ such that 
$$
\Vert v\Vert_{L^2(\ppp\OOO)} \le C_{\delta}\Vert v\Vert
_{H^{\delta+\frac{1}{2}}(\OOO)} \quad \mbox{for all
$v \in H^1(\OOO)$.}
$$
We fix $\delta \in \left(0, \, \frac{1}{2}\right)$.
The interpolation inequality for the Sobolev spaces yields that for 
any $\ep>0$ we can find a constant $C_{\ep,\delta} > 0$ such that 
$$
\Vert v\Vert_{H^{\delta+\frac{1}{2}}(\OOO)} 
\le \ep\Vert \nabla v\Vert_{L^2(\OOO)} + C_{\ep,\delta}\Vert v\Vert_{L^2(\OOO)}
\quad \mbox{for all $v \in H^1(\OOO)$}
$$
(e.g., Chapter IV in \cite{Ad} or Sections 2.5 and 11 of Chapter 1 
in \cite{LM}).
Therefore,
\begin{align*}
& \Vert v\Vert^2_{L^2(\ppp\OOO)}
\le (\ep C_{\delta}\Vert \nabla v\Vert_{L^2(\OOO)}
+ C_{\delta}C_{\ep,\delta}\Vert v\Vert_{L^2(\OOO)})^2\\
\le& 2(\ep C_{\delta})^2\Vert \nabla v\Vert_{L^2(\OOO)}^2
+ 2(C_{\delta}C_{\ep,\delta})^2\Vert v\Vert_{L^2(\OOO)}^2
\end{align*}
for all $v \in H^1(\OOO)$.
Substituting this inequality into \eqref{(2)}, we obtain
\begin{align*}
&(\kappa - 2C(\ep C_{\delta})^2) \Vert \nabla v\Vert_{L^2(\OOO)}^2\\
+& ( M - \frac{1}{2}\Vert \mbox{div}\, b\Vert_{C(\ooo{\OOO}\times [0,T])}
- 2(CC_{\delta}C_{\ep,\delta})^2)\Vert v\Vert_{L^2(\OOO)}^2
\le 0.
\end{align*}
Choosing $\ep > 0$ sufficiently small such that 
$\kappa - 2C(\ep C_{\delta})^2) > 0$, since $M>0$ is sufficiently large
such that 
$$
M > \frac{1}{2}\Vert \mbox{div}\, b\Vert_{C(\ooo{\OOO}\times [0,T])}
+ 2(CC_{\delta}C_{\ep,\delta})^2,
$$
we can complete the proof of Lemma \ref{lem1}.
\end{proof}

{\bf Second Step.}
\\
By \eqref{1}, we can apply Theorem 3.2 (p.137) in \cite{LU} to see that 
there exists a constant $\theta \in (0,1)$ such that for each 
$t \in [0,T]$, a solution $\psi(\cdot,t) \in C^{2+\theta}(\ooo{\OOO})$ to 
\eqref{(4.3)} exists uniquely.

We set 
\begin{equation}\label{(8.3a)}
\eta(t):= \Vert \psi(\cdot,t)\Vert_{C^{2+\theta}(\ooo{\OOO})}, \quad
0\le t \le T.
\end{equation}
We note that $\eta(t) < \infty$ for each $t \in [0,T]$.

The space $C^{2+\theta}(\ooo{\OOO})$ is the Schauder space 
defined at the end of the proof of Lemma \ref{l6.1} in Section \ref{sec4}.

Furthermore, we obtain that for arbitrary $G \in \HOLSZ$, 
there exists a unique solution $w=w(\cdot,t)$ to
\begin{equation}\label{(3)}
\left\{ \begin{array}{rl}
& A_1(t)w = G \quad \mbox{in $\OOO$}, \\
& \ppp_{\nu_A}w + \sigma w = 0 \quad \mbox{on $\ppp\OOO$}
\end{array}\right.
\end{equation}
for each $t \in [0,T]$.

Now, we can derive:
for each $t\in [0,T]$, there exists a constant $C_t > 0$ such that 
\begin{equation}\label{(4)}
\Vert w(\cdot,t)\Vert_{\HOLST} \le C_t\Vert G\Vert_{\HOLSZ}
\end{equation}
for all $w$ satisfying \eqref{(3)}.  Here the constant $C_t>0$ depends on
$\Vert a_{ij}\Vert_{C^1(\ooo{\OOO})}$, 
$\Vert b_k\Vert_{C([0,T];C^2(\ooo{\OOO}))}$ for $1\le i,j\le d$ and 
$0\le k\le d$, but not on respective choices of these coefficients.
\\
{\bf Verification of \eqref{(4)}.}
\\
For each $t \in [0,T]$, the inequality
\begin{equation}\label{(5)}
  \Vert w(\cdot,t)\Vert_{\HOLST} \le C_t(\Vert G\Vert_{\HOLSZ}
+ \Vert w(\cdot,t)\Vert_{C(\ooo{\OOO})})
\end{equation}
holds true (see, e.g., the formula (3.7) on p. 137 in \cite{LU}).
Then we have to eliminate the second term, namely the term  
$\Vert w(\cdot,t)\Vert_{C(\ooo{\OOO})}$,
on the right-hand side of this inequality.
This can be done by a classical compactness-uniqueness argument.  
More precisely,
assume that \eqref{(5)} does not hold.  Then there exist $w_n\in \HOLST$,
$G_n\in \HOLSZ$ for $n\in \N$ such that $\Vert w_n\Vert_{\HOLST} = 1$ and 
$\lim_{n\to\infty}\Vert G_n\Vert_{\HOLSZ} = 0$.
By the Ascoli-Arzel\`a theorem, we can extract a subsequence $w_{k(n)}$
from $w_n$ with $n\in \N$ such that 
$w_{k(n)} \longrightarrow \www{w}$ in $C(\ooo{\OOO})$ as $n\to \infty$.
Applying \eqref{(5)} to 
$$
A_1(t)(w_{k(n)} - w_{k(m)}) = G_{k(n)} - G_{k(m)} \quad 
\mbox{in $\OOO$}
$$
with $\ppp_{\nu_A}(w_{k(n)} - w_{k(m)}) + 
\sigma(w_{k(n)} - w_{k(m)}) = 0$ on $\ppp\OOO$, we obtain
\begin{align*}
&  \Vert w_{k(n)} - w_{k(m)}\Vert_{\HOLST}\\
\le & C_t(\Vert G_{k(n)} - G_{k(m)}\Vert_{\HOLSZ}
+ \Vert w_{k(n)} - w_{k(m)}\Vert_{C(\ooo{\OOO})})
\,\longrightarrow 0
\end{align*}
as $n,m \to \infty$.  Hence, there exists $w_0 \in \HOLST$ such that 
$w_{k(n)} \longrightarrow w_0$ in $\HOLST$, and so 
$$
\Vert w_0\Vert_{\HOLST} = \lim_{n\to\infty} \Vert w_{k(n)}\Vert_{\HOLST}
= 1
$$
and $G_{k(n)} = A_1(t)w_{k(n)} \longrightarrow A_1(t)w_0$ in 
$\HOLSZ$.

Since $\lim_{n\to\infty} \Vert G_{k(n)}\Vert_{\HOLSZ} = 0$, we reach 
$A_1(t)w_0 = 0$ in $\OOO$ with $\ppp_{\nu_A}w_0 + \sigma w_0 = 0$ on 
$\ppp\OOO$.  Then Lemma \ref{lem1} yields $w_0=0$ in $\OOO$, which 
contradicts $\Vert w_0\Vert_{\HOLST} = 1$.
Thus the verification of \eqref{(4)} is complete.
$\square$
\\

{\bf Third Step.}
\\
We have to prove that $\psi(x,t)$ constructed in Second Step 
satisfy the inclusion $\psi \in C^1([0,T];C^2(\ooo{\OOO}))$.

Fixing $t\in [0,T]$ arbitrarily, we verify that
$d(x,s):= \psi(x,t) - \psi(x,s)$ satisfies
\begin{equation}\label{(6)}
\left\{ \begin{array}{rl}
& -A_1(t)d(\cdot,s) = (b_0(t) - b_0(s))\psi(\cdot,s)\\
- & \sum_{j=1}^d (b_j(t) - b_j(s))\ppp_j\psi(\cdot,s) \quad 
\mbox{in $\OOO$ for $0\le s, t \le T$},\\
& \ppp_{\nu_A}d + \sigma d = 0 \quad \mbox{on $\ppp\OOO$, $\,\,$ 
$0 \le s,t \le T$}.
\end{array}\right.
\end{equation}

Let $\delta>0$ be arbitrarily fixed.  We set 
$I_{\delta,t}:= [0,T] \cap \{s;\, \vert t-s\vert \le \delta\}$.

Again the application of \eqref{(8.3a)} and \eqref{(4)} to \eqref{(6)} yields 
\begin{align*}
& \Vert d(\cdot,s)\Vert_{\HOLST}\\
\le & C\left(\left\Vert \sum_{j=1}^d (b_j(t)-b_j(s))\ppp_j\psi(\cdot,s)
\right\Vert_{\HOLSZ}
+ \Vert (b_0(t)-b_0(s))\psi(\cdot,s)\Vert_{\HOLSZ}\right)\\
\end{align*}
\begin{equation}\label{(8.8a)}
\le C\sum_{j=0}^d \Vert b_j(s) - b_j(t)\Vert_{\CONE}\eta(s)
\le C\max_{0\le j \le d} \Vert b_j(s) - b_j(t)\Vert_{\CONE}\,
\sup_{s \in I_{\delta,t}}\eta(s).
\end{equation}

Setting 
$$
h(\delta):= \max_{0\le j\le d} \sup_{\vert s-t\vert \le \delta}
\Vert b_j(s) - b_j(t)\Vert_{\CONE},
$$
by the condition $b_j \in C([0,T];C^1(\ooo{\OOO}))$ for $0\le j \le d$, we
deduce that $\lim_{\delta\downarrow 0} h(\delta) = 0$.

Moreover, in terms of $\eta$ defined by \eqref{(8.3a)}, we can rewrite \eqref{(6)} as
$$
\vert \eta(s) - \eta(t)\vert \le Ch(\delta)\sup_{s\in I_{\delta,t}} \eta(s)
\quad \mbox{for $s \in I_{\delta,t}$},
$$
and so 
$$
\eta(s) \le \eta(t) + Ch(\delta)\sup_{s\in I_{\delta,t}} \eta(s)
\quad \mbox{for $s \in I_{\delta,t}$}.
$$
Choosing $\delta: =\delta(t)>0$ sufficiently small 
for given $t \in [0,T]$, we see that
$\sup_{s\in I_{\delta(t),t}} \eta(s) \le C_1\eta(t)$.
Varying $t \in [0,T]$, we can choose a finite number of intervals 
$I_{\delta(t),t}$ covering $[0,T]$, we obtain
\begin{equation}\label{(8.8)}
\Vert \psi\Vert_{L^{\infty}(0,T;C^{2+\theta}(\ooo{\OOO}))} \le 
C_2                             
\end{equation}
with some constant $C_2>0$.

Substitution of \eqref{(8.8)} into \eqref{(8.8a)} yields 
$$
\Vert d(\cdot,s)\Vert_{C^{2+\theta}(\ooo{\OOO}))}
= \Vert \psi(\cdot,t) - \psi(\cdot,s)\Vert_{C^{2+\theta}(\ooo{\OOO}))}
\le Ch(\delta)C_2
$$
for $s \in I_{\delta(t),t}$.
Consequently, 
$\lim_{s\to t} \Vert \psi(\cdot,s) - \psi(\cdot,t)\Vert
_{\HOLST} = 0$, that is,
\begin{equation}\label{(8)}
\psi \in C([0,T];\HOLST).
\end{equation}

Finally we have to prove $\psi \in C^1([0,T];\HOLST)$.
Since $-A_1(\xi)\psi(x,\xi) = 1$ in $\OOO$ with $\xi = t,s$,  
differentiating in $t$ and $s$, we can obtain
\begin{equation}\label{(9)}
\sum_{j=1}^d \ppp_i(a_{ij}(x)\ppp_j\ppp_{\xi}\psi(x,\xi))
+ \sum_{j=1}^d b_j(\xi)\ppp_j\ppp_{\xi}\psi(x,\xi)
- b_0(\xi)\ppp_{\xi}\psi(x,\xi)
\end{equation}
$$
= -\sum_{j=1}^d \ppp_{\xi}b_j(\xi)\ppp_j\psi(x,\xi)
+ (\ppp_{\xi}b_0)(\xi)\psi(x,\xi) \quad \mbox{in $\OOO$}
$$
with $\xi = t,s$.  Therefore, by subtracting the equation \eqref{(9)} with 
$\xi=s$ from the one with $\xi=t$, we deduce that 
$d_1(x,s):= (\ppp_t\psi)(x,t) - (\ppp_t\psi)(x,s)$ satisfies
\begin{equation}\label{(10)}
-A_1(t)d_1(x,s)
\end{equation}
\begin{align*}
=& \biggl[ -\sum_{i,j=1}^d (b_j(t)-b_j(s))\ppp_j\ppp_s\psi(x,s)
+ (b_0(t)-b_0(s))\ppp_s\psi(x,s)\\
- & \sum_{j=1}^d (\ppp_tb_j(t)-\ppp_sb_j(s))\ppp_j\psi(x,s)
+ (\ppp_tb_0(t)-\ppp_sb_0(s))\psi(x,s)\biggr]  \\
+ &\left[ -\sum_{j=1}^d \ppp_tb_j(t)(\ppp_j\psi(x,t) - \ppp_j\psi(x,s))
 + \ppp_tb_0(t)(\psi(x,t) - \psi(x,s))\right]\\
=:& H_1(x,t,s) + H_2(x,t,s) \quad \mbox{in $\OOO$}
\end{align*}
with $\ppp_{\nu_A}d_1 + \sigma d_1 = 0$ on $\ppp\OOO$ for 
$s,t \in [0,T]$.
Thus, if we can verify the relation $\lim_{s\to t} \Vert d_1(\cdot,s)\Vert
_{\HOLST} = 0$, then we can complete the proof of
$\psi\in C^1([0,T];\HOLST)$.
To this end, by applying Theorem 3.2 (p.137) in \cite{LU} to \eqref{(10)},
it suffices to prove that 
\begin{equation}\label{(11)}
\lim_{s\to t} \Vert H_{\ell}(\cdot,t,s)\Vert_{\HOLSZ} = 0, \quad
\ell=1,2. 
\end{equation}
Indeed the final limit implies that 
$\lim_{s\to t} \Vert d_1(\cdot,t,s)\Vert_{\HOLST}$, that is,
$\ppp_t\psi \in C^1([0,T];\HOLST)$, so that the proof of the existence of
$\psi$ is complete.
\\
{\bf Verification of \eqref{(11)}.}
\\
Applying Theorem 3.2 in \cite{LU} to \eqref{(9)}, in view of the 
regularity in \eqref{(1.2)} and \eqref{(9)}, we see 
\begin{equation}\label{(12)}
\Vert \ppp_t\psi(\cdot,t)\Vert_{\HOLST}
\end{equation}
\begin{align*}
\le& C\left( \left\Vert \sum_{j=1}^d (\ppp_tb_j)(\cdot,t)
\ppp_j\psi(\cdot,t)\right\Vert_{\HOLSZ}
+ \Vert (\ppp_tb_0)(\cdot,t)\psi(\cdot,t)\Vert_{\HOLSZ}\right)\\
\le& C\sum_{j=0}^d \Vert b_j(\cdot,t)\Vert_{C^1([0,T];\CONE)}
\Vert \psi(\cdot,t)\Vert_{\HOLST}
\le C_3 \quad \mbox{for $0\le t\le T$.}
\end{align*}
Hence \eqref{(8.8)} and \eqref{(12)} yield
\begin{equation}\label{(13)}
\Vert H_1(\cdot,t,s)\Vert_{\HOLSZ}
\le C_4\sum_{k=0}^1\sum_{j=0}^d
\Vert (\ppp_t^kb_j)(\cdot,t) - (\ppp_t^kb_j)(\cdot,s)\Vert
_{\CONE}.
\end{equation}
Similarly we can prove
\begin{equation}\label{(13a)}
\Vert H_2(\cdot,t,s)\Vert_{\HOLSZ}
\le C_4\sum_{j=0}^d \Vert \ppp_tb_j\Vert_{C([0,T];\CONE)}
\sum_{k=0}^1 \Vert \nabla^k\psi(\cdot,t) - \nabla^k\psi(\cdot,s)
\Vert_{\HOLSZ}.
\end{equation}
Since $\ppp_t^kb_j \in C([0,T];C^1(\ooo{\OOO}))$ for $k=0,1$ and
$0\le j \le d$ by \eqref{(1.2)} and $\nabla^k\psi \in C([0,T];C^{1+\theta}
(\ooo{\OOO}))$ with $k=0,1$ by \eqref{(8)}, from \eqref{(13)} and \eqref{(13a)}
it follows that 
$\lim_{s\to t} \Vert H_{\ell}(\cdot,t,s)\Vert_{\HOLSZ} = 0$ for
$\ell=1,2$.  Thus the verification of \eqref{(11)} is complete.
$\square$

%%=============================================================%%
%% Sample for another appendix section			       %%
%%=============================================================%%

%% \section{Example of another appendix section} \label{secA2}%

%%%%%%%%%%%%%%%%%%%%%%%%%% BACK MATTERS %%%%%%%%%%%%%%%%%%%%%%%%%%%%

\section*{Acknowledgements}
 % Acknowledgments are not compulsory. Where included they should be brief.
 % Grant or contribution numbers may be acknowledged, or help by colleagues. Example:
The second author was supported by Grant-in-Aid 
for Scientific Research Grant-in-Aid (A) 20H00117 of 
Japan Society for the Promotion of Science, 
the National Natural Science Foundation of China
(no. 11771270, 91730303), and the RUDN University 
Strategic Academic Leadership Program.
% \end{acknowledgements}

%%%%%%%%%%%%%%%%%% REFERENCES: %%%%%%%%%%%%%%%%%%%%%%%%%%%%%%%%%%%%%%%%%%%%
%% BibTeX users: please use \bibliographystyle{spmpsci} %% for math. and phys. sci.
%% Non-BibTeX users: please use the model as below !!! %%

%%%% for FCAA - pls. include directly the Refs items here ! %%%
%%%% following STRICTLY the models below %%%%%
%%%% and ARRANGE the items in ALPHABETIC ORDER for authors' family names !!!

 %%%%%%%%%%%%%%%%%%%%%%%%%%%%%%

\end{document}